\documentclass[10pt]{article} 

\title{On the Isomorphism Conjecture in algebraic $K$-Theory} 
\author{Arthur Bartels, Tom Farrell, Lowell Jones, Holger Reich}
\date{July 31, 2001}

\usepackage{amsmath}
\usepackage{amsthm}
\usepackage{amssymb}
\usepackage{amsfonts}  
\usepackage[arrow,matrix,cmtip]{xy} \SelectTips{cm}{}

\newcommand{\cala}{{\cal A}}
\newcommand{\calb}{{\cal B}}
\newcommand{\calc}{{\cal C}}
\newcommand{\cald}{{\cal D}}
\newcommand{\cale}{{\cal E}}
\newcommand{\calf}{{\cal F}}
\newcommand{\calg}{{\cal G}}
\newcommand{\calh}{{\cal H}}

\newcommand{\IE}{{\mathbb E}} 
\newcommand{\IF}{{\mathbb F}}
  
\newcommand{\IH}{{\mathbb H}}

\newcommand{\IK}{{\mathbb K}}

\newcommand{\IN}{{\mathbb N}}

\newcommand{\IR}{{\mathbb R}}

\newcommand{\IZ}{{\mathbb Z}}

\newcommand{\All}{{\mathit{All}}} 
\newcommand{\Cyc}{{\mathit{Cyc}}} 
\newcommand{\VCyc}{{\mathit{VCyc}}} 
\newcommand{\mor}{{\rm mor}} 
\newcommand{\supp}{{\rm supp}} 
\newcommand{\sma}{{\wedge}} 
\newcommand{\wed}{{\bigvee}} 
\newcommand{\einsu}{{\left[ 1, \infty \right)}}
\newcommand{\einsua}{{\left[ 1 , \infty \right]}}
\newcommand{\Or}{{\rm Or}}     

\newcommand{\coker}{{\rm coker}}
\newcommand{\res}{{\rm res}}
\newcommand{\ind}{{\rm ind}}
\newcommand{\dni}{{\rm dni}} 
\newcommand{\inc}{{\rm inc}}

\newcommand{\id}{{\rm id}}
\newcommand{\tr}{{\rm tr}}
\newcommand{\map}{{\rm map}}
\newcommand{\colim}{{\rm colim }}

\newcommand{\M}{{\rm M}}          

\newcommand{\dist}{{\rm dist}}
\newcommand{\cone}{{\rm cone}}
\newcommand{\ch}{{\rm ch}}      
\newcommand{\St}{{\rm St}}      
\newcommand{\Nil}{{\rm Nil}}      
\newcommand{\sh}{{\mathit{sh}}}

\newcommand{\punkt}{\{\mathit{pt}\}}
\newcommand{\e}{\epsilon}
\newcommand{\hra}{\hookrightarrow}

\newcommand{\x}{\times}
\newcommand{\dd}{\partial}

\theoremstyle{plain}
\newtheorem{theorem}{Theorem}[section]
\newtheorem{lemma}[theorem]{Lemma}
\newtheorem{corollary}[theorem]{Corollary}
\newtheorem{proposition}[theorem]{Proposition}
\newtheorem{conjecture}[theorem]{Conjecture}
\newtheorem{addendum}[theorem]{Addendum}

\theoremstyle{definition}
\newtheorem{definition}[theorem]{Definition}
\newtheorem{example}[theorem]{Example}
\newtheorem{notation}[theorem]{Notation}

\theoremstyle{remark}
\newtheorem{remark}[theorem]{Remark}
\newtheorem{claim}[theorem]{Claim}

\begin{document}


\typeout{---------------------- main -----------------------------}

\maketitle

\typeout{--------------------- Inhalt ----------------------------}

\tableofcontents

\typeout{---------------------...und los -------------------------}


\typeout{--------------------- introduction ----------------------}


\section{Introduction}

\subsection{The Isomorphism Conjecture in algebraic $K$-theory}

This paper deals with the 
computation of the algebraic $K$-theory groups 
\[
K_n( R \Gamma ),    \quad n \in \IZ.
\]
Here $R \Gamma$ is a group ring with $R$ an arbitrary associative  
ring with unit and $\Gamma$ is a
discrete group.
Ideally one would like to express $K_n( R \Gamma )$ in terms of 
the group homology of $\Gamma$ and the algebraic $K$-theory of the coefficient 
ring $R$. More precisely there exists a so called assembly map 
\[
A: H_{\ast}( B \Gamma , \IK^{-\infty} R ) \to K_{\ast} ( R \Gamma )
\]
which is believed to be an isomorphism for nice rings and torsion free groups.
The left hand side is the generalized homology theory associated to
the non-connective algebraic $K$-theory spectrum of the ring $R$, i.e.\
on the left we have the homotopy groups of the spectrum $B \Gamma_{+} \sma \IK^{- \infty} R$.
In favorable cases this left hand side is accessible to computations and leads to
very concrete results. Compare Corollary~\ref{maincorollary} below.

Conjectures of this type go back to \cite{Hsiang(1984)}. There are analogous assembly maps and 
conjectures in $L$-theory, see \cite{Quinn(1970)} and \cite{Quinn(1971)}, and for 
the (topological) 
$K$-theory of $C^{\ast}$-algebras \cite{Baum-Connes-Higson(1994)}. 
For a general overview see \cite{Ferry-Ranicki-Rosenberg(1995)}.

There are rings where surjectivity fails already for the infinite cyclic group. 
In these  cases the Bass-Heller-Swan 
splitting \cite{Bass-Heller-Swan(1964)} yields additional summands, the
so called Nil-terms on the right hand side. 
It is also easy to see that in general 
the map can not be surjective if the group is finite or contains torsion.

In \cite{Farrell-Jones(1993)} the second and third author consider a 
modified version of the assembly map building on 
Quinn's ``homology for simplicially stratified fibrations'' \cite{Quinn(1982)}.
The modified assembly map takes the above mentioned problems into account and it is conjectured
in \cite[1.6]{Farrell-Jones(1993)} to be 
an isomorphism for arbitrary groups 
and $R$ the ring of integers. 

The results of this paper suggest that the assumption $R=\IZ$ is not necessary.
The conjecture, if true for a group $\Gamma$, reduces (at least in principle) the computation of
$K_n ( R \Gamma )$ to homological computations  and the computation of the $K$-theory of
the group rings $RH$, where $H$ is finite or contains an infinite  cyclic subgroup of finite 
index.
Roughly speaking the $K$-theory of finite subgroups of $\Gamma$ and the Nil-groups
are built into the left hand side.
The general principle behind such modified assembly maps was clarified 
by Davis and L\"uck in \cite{Davis-Lueck(1998)} building on \cite{Weiss-Williams(1995)}.

We proceed to explain the modified assembly map using their language.
We need the concept of a universal space for a family of 
subgroups which generalizes the universal free $\Gamma$-space $E \Gamma$, the universal
covering of the classifying space $B\Gamma$.

\begin{definition}
A set of subgroups of $\Gamma$ is called a family if it is closed under conjugation
and taking subgroups.
We use the notation 
\[
1, \quad \Cyc, \quad  \VCyc \quad \mbox{  and } \quad \All
\]
for the trivial family, the family
of cyclic subgroups, the family of virtually cyclic subgroups and the family of all
subgroups of a given group $\Gamma$. Recall that a group is called virtually cyclic
if it is finite or contains an infinite cyclic subgroup of finite index.
For a family $\calf$ of subgroups of $\Gamma$ we denote by 
\[
E \Gamma ( \calf )
\]
the universal space among $\Gamma$-spaces with isotropy in $\calf$. 
It is characterized by the universal property that for every $\Gamma$-CW complex $X$ 
whose isotropy groups are all in the family $\calf$ one can find an equivariant continuous 
map $X \to E \Gamma ( \calf )$ which is unique up to equivariant homotopy.
A $\Gamma$-CW complex $E$ is a model for $E \Gamma ( \calf )$ if the fixpoints
$E^H$ are  contractible for $H \in \calf$ and empty otherwise. 
\end{definition}

In \cite{Davis-Lueck(1998)} Davis and L{\"u}ck construct a generalized equivariant homology 
theory $H_n^{\Or \Gamma}( X ; \IK R^{- \infty})$
for $\Gamma$-CW complexes associated to a jazzed up version of the non-connective algebraic 
$K$-theory spectrum $\IK R^{-\infty}$. It has a simple description in terms of balanced
products over the orbit category $\Or \Gamma$. The construction is reviewed in
Section~\ref{identify}. In their formulation the modified assembly map
is simply the map induced by $E \Gamma ( \VCyc ) \to E \Gamma ( \All ) = \punkt$
and the above mentioned conjecture reads as follows.

\begin{conjecture}[Isomorphism Conjecture] \label{isomorphismconjecture} 
$\quad$ \\
The map 
\[
\xymatrix{
H_n^{\Or \Gamma}( E \Gamma ( \VCyc ), \IK R^{-\infty}  ) 
\ar[rr]^-{A_{\VCyc \to \All }} & &
H_n^{\Or \Gamma}( \punkt , \IK R^{-\infty}  ) =K_n ( R \Gamma )
         }
\]
induced by 
\[
E \Gamma( \VCyc ) \to E \Gamma ( \All ) = \punkt
\]
is an isomorphism.
\end{conjecture}

Similar to the classical case there is an Atiyah-Hirzebruch type spectral sequence 
computing the left hand side, see \cite[4.7]{Davis-Lueck(1998)}. 
Further computational tools for the left hand side are developed
in \cite{Lueck-Stamm(2000)}. In favorable cases, e.g.\ if $R$ is a field of characteristic
$0$ the family of virtually cyclic subgroups can be replaced by the smaller
family of finite subgroups. For this family 
a very concrete description of the rationalized left hand side 
is obtained in \cite{Lueck(2000)}. \label{ask Wolfgang!!!}
The conjecture should therefore be seen as a conceptional approach towards a computation
of $K_n( R \Gamma )$.

The conjecture has been verified by the second and third author for 
$n \leq 1$ (and rationally for all  $n$)
in the case where $R=\IZ$ for a large class of groups. 
This class includes: Every discrete cocompact subgroup of any (virtually connected)
Lie group, or any subgroup of one of these \cite{Farrell-Jones(1993)}.
The fundamental group of any complete pinched negatively curved
Riemannian manifold \cite{Farrell-Jones(1998)}. The fundamental group of any complete
non-positively curved locally symmetric space \cite{Farrell-Jones(1998)} 
and the fundamental group of any complete 
$A$-regular non-positively curved Riemannian manifold \cite{Farrell-Jones(1998)}.

In all these cases the restriction to $R = \IZ$ was unavoidable because proofs
went through the translation into questions about $h$-cobordisms \cite{Farrell-Jones(1986)}
or pseudoisotopies \cite{Farrell-Jones(1991)}. This translation only works if $R$ is the ring
of integers.

Some very important results about $K_n( R\Gamma)$ for infinite groups and
arbitrary coefficient rings are: 
The classical Bass-Heller-Swan formula which deals with the infinite cyclic group. 
The fundamental paper \cite{Waldhausen(1978)} of Waldhausen, giving 
exact sequences for amalgamated products and $\mathit{HNN}$-extensions of groups
and \cite{Carlsson-Pedersen(1995)} where 
Carlsson and Pedersen discover 
a big portion of $K_n( R \Gamma )$ for groups arising in non-positive curvature
situations by proving that the classical assembly map $A=A_{1 \to \All}$ is injective.

\subsection{Main Results and Corollaries}

Our aim in this paper is to prove the following two theorems.

\begin{theorem} \label{maintheorem}
Let $\Gamma$ be the fundamental group of a compact Riemannian
manifold with strictly negative sectional curvature. Let $R$
be an associative ring with unit. Then the assembly map
\[
A_{\Cyc \to \All }:H_n^{\Or \Gamma}( E \Gamma ( \Cyc ), \IK R^{-\infty}  ) \to K_n ( R \Gamma )
\]
is an isomorphism for $n \leq 1$.
\end{theorem}

\begin{theorem} \label{maintheoreminjectivity}
The map $A_{\Cyc \to \All }$ in Theorem~\ref{maintheorem} is injective for all $n \in \IZ$.
\end{theorem}

The injectivity statement \ref{maintheoreminjectivity} will be proven in Theorem~\ref{injectivity}.
The remaining surjectivity statement in \ref{maintheorem} is Theorem~\ref{surjectivity} for $n=1$. 
The case $n<1$ is then a consequence of Corollary~\ref{loweriso}.

Note that the negatively curved 
compact Riemannian manifold is a model for the classifying space of 
$\Gamma$. In particular $\Gamma$ contains no torsion and therefore the family of
virtually cyclic subgroups of $\Gamma$ coincides with the 
family of cyclic subgroups \cite[2.6.6.(iii)]{Stamm(1999)}. 
Theorem~\ref{maintheorem} therefore verifies the Isomorphism
Conjecture~\ref{isomorphismconjecture} for this class of groups in dimensions less than $2$.


Recall that $NK_n(R)$ is defined as the cokernel of the split injection $K_n(R) \to K_n( R[t] )$.
Often this group is also denoted $\Nil_{n-1}(R)$ and called the Nil-group of $R$. 
A ring is called (right) regular if it is right Noetherian 
and every finitely generated right $R$-module admits a finite dimensional
projective resolution. For regular rings all Nil-groups vanish \cite{Bass(1972)}.
Theorem~\ref{maintheorem}  has the following consequence.

\begin{corollary} \label{maincorollary}
Let $\Gamma$ be as in Theorem~\ref{maintheorem}.
If for the ring $R$ the groups $NK_n(R)$ vanish for $n\leq 1$, e.g.\ if $R$ is regular, 
then the assembly map
\[
A = A_{1 \to \All}: H_n ( B \Gamma , \IK^{-\infty} R ) \to K_n( R \Gamma )
\]
is an isomorphism for $n \leq 1$.
If moreover $K_{-n}(R)=0$ for $n \geq 1$, e.g.\ if $R$ is regular, we obtain
\begin{eqnarray*}
K_{-n}( R \Gamma )  & = & 0 \mbox{ for } n \geq 1 ,  \\
K_0( R \Gamma )  & = & K_0( R ) \mbox{ and }\\
K_1( R \Gamma )  & = & \Gamma_{ab} \otimes K_0(R) \oplus K_1(R) . 
\end{eqnarray*}
Here $\Gamma_{ab}$ denotes the abelianized group.
\end{corollary}

\begin{proof}
For a group $\Gamma$ and two families of subgroups $\calf \subset \calg$
let $A_{\calf \to \calg}(\Gamma)$ denote the assembly map
induced by $E \Gamma ( \calf ) \to E \Gamma ( \calg )$. 
The vanishing of Nil-terms implies that 
for every cyclic group $C$ the map 
$A_{1 \to \Cyc}(C)$ is an isomorphism in degrees less than $1$. The same holds for 
$A_{\Cyc \to \All}(\Gamma)$ by \ref{maintheorem}.
It follows from the ``transitivity'' proven 
in \cite[2.3]{Lueck-Stamm(2000)} that $A_{1 \to \All}(\Gamma)$
is an isomorphism. Compare also \cite[A.10]{Farrell-Jones(1993)}. This map can be identified with
the ``classical'' assembly map $A$ above. Evaluating the Atiyah-Hirzebruch
spectral sequence for the homology in degrees less than $1$ leads to
the above formulas.
\end{proof}

There is an assembly map for for Nil-groups, compare Subsection~\ref{nilsection}.
The fact that we impose no conditions on the coefficient ring allows us to 
prove that it is an isomorphism under the assumptions 
of~\ref{maintheorem}. As a consequence we obtain for example:
\begin{corollary}  \label{nilcorollary}
Let $\Gamma$ be as in Theorem~\ref{maintheorem}.
If $NK_n (R)=0$ for $n \leq 1$, e.g.\ if $R$ is regular, then 
\[
NK_n( R\Gamma )=0 \mbox{  for  } n \leq 1.
\] 
\end{corollary}

\begin{proof} This follows by combining 
Theorem~\ref{maintheorem} with Proposition~\ref{nilprinciple}. 
\end{proof}

A similar argument leads to the following result.

\begin{corollary}
Let $\Gamma = \Gamma_1 \times  \dots \times  \Gamma_l$, where the groups $\Gamma_i$ 
are as in Theorem~\ref{maintheorem}.
Let $R$ be regular, then the conclusions
of Corollary~\ref{maincorollary} hold for $\Gamma$.
\end{corollary}
\begin{proof}
Observe that for a ring $S$ and groups $G$ and $H$ we have $SG[H]=S[G \times H]$ and similarly 
$SG[t]=S[t]G$. 
Let 
\[
A(G,S): BG_+ \sma \IK^{-\infty}S \to \IK^{-\infty}SG
\]
be the ``classical'' assembly map for 
$G$ with coefficient ring $S$.
If $R$ is regular, then so are
the polynomial rings $R[t_1], \dots ,R[t_1, \dots, t_{l-1}]$ and hence their Nil-groups
vanish. Now use repeatedly Corollary~\ref{maincorollary} and Proposition~\ref{nilprinciple}
to conclude that the assembly maps 
\begin{eqnarray*}
\begin{array}{ccc}
A(\Gamma_1,R),\dots & , A(\Gamma_1, R[t_1, \dots , t_{l-2}]), & A(\Gamma_1,R[t_1,\dots , t_{l-1}]),  \\
A(\Gamma_2 , R\Gamma_1), \dots & , A(\Gamma_2, R\Gamma_1[t_1, \dots , t_{l-2}]), &  \\
\vdots &  & \\
A(\Gamma_l, R[\Gamma_1 \times \dots \times \Gamma_{l-1}])  & &
\end{array}
\end{eqnarray*}
are isomorphisms in degrees less than $2$. The result follows since
\[
A(\Gamma_1 \times \Gamma_2 , R)=
A(\Gamma_2 , R \Gamma_1 ) \circ \id_{{B\Gamma_2}_+} \sma A( \Gamma_1 , R) \mbox{  etc.}
\]
\end{proof}

Finally we would like to mention the following result (Corollary~\ref{lhscomputation}) 
which shows that, for the groups considered in Theorem~\ref{maintheorem},
the Nil-terms split off as direct summands from the left hand side and hence by our main result 
in low dimensions also from the right hand side
of the assembly map $A_{\Cyc \to \All}$.

\begin{proposition} \label{introlhscomputation}
Let $\Gamma$ be as in Theorem~\ref{maintheorem}. We have for all $n \in \IZ$
and every ring $R$
\[
H_{n}^{\Or \Gamma}( E \Gamma ( \Cyc ); \IK R^{-\infty} ) \cong  
 H_{n}(B\Gamma;\IK R^{-\infty}) \oplus \bigoplus_{I} (NK_{n}(R) \oplus NK_{n}(R)).
\]
where the right hand direct sum is indexed over the set $I$ of conjugacy classes
of maximal cyclic subgroups of $\Gamma$.
\end{proposition}

\subsection{A brief outline}

The proof of our main result and hence also this paper fall into three major parts. 
In the first part, consisting of 
Sections~\ref{preparations} to \ref{lower}, we construct a model
for the assembly map which allows us to bring 
geometric assumptions into the picture.
The second part, consisting of Sections~\ref{further} to  \ref{sec_injectivity},
proves injectivity of the assembly map.
The third part, consisting of Sections~\ref{versus} to \ref{sec_surjectivity},
proves the surjectivity statement.
The surjectivity and injectivity parts are completely independent of each other.
So the reader who is only interested in the surjectivity part can skip Sections~\ref{further}
to \ref{sec_injectivity}.

{\bf First part.}
In the first part we construct a functor $\IK^{-\infty} \cald^{\Gamma}(X)$
from $\Gamma$-CW complexes to spectra and prove that it is an
equivariant homology theory, compare Section~\ref{homology}. This means in particular that 
the homotopy groups of this spectrum form an equivariant homology theory in the
ordinary sense, e.g.\ we have homotopy invariance and Mayer-Vietoris sequences.
It is illuminating to consider the
two extreme cases: If  $X$ is a free $\Gamma$-space, e.g.\ $E \Gamma$, then we obtain
up to a shift the (nonequivariant) generalized homology associated to the nonconnective 
algebraic $K$-theory spectrum of the ring $R$ evaluated on the quotient space:
\[
\pi_{n+1} \IK^{-\infty} \cald^{\Gamma}( X) = H_n( X/ \Gamma ; \IK^{-\infty}(R) ).
\]
If $X$ is a point we obtain the algebraic $K$-theory of the group ring $R \Gamma$:
\[
\pi_{n+1} \IK^{-\infty} \cald^{\Gamma}( \punkt )=K_n( R \Gamma).
\]
The functor $\IK^{-\infty} \cald^{\Gamma}( X)$ is in fact the composition of the 
well studied functor $\IK^{-\infty}$ from additive categories to spectra and
the functor $\cald^{\Gamma}(-)$ defined in Section~\ref{functor} below.

More precisely given any $\Gamma$-CW complex we associate to it an 
additive category $\cald^{\Gamma}(X)$, see subsection~\ref{definitions}. 
Its objects are ``geometric modules'' 
over the ``cone'' $X \times \Gamma \times \einsu$. Morphisms are required to satisfy
a control condition near infinity. This is in the spirit of 
\cite{Pedersen-Weibel(1989)} and \cite{Anderson-Connolly-Ferry-Pedersen(1994)}.
For our purposes we have to develop an equivariant version which also digests
$\Gamma$-spaces with infinite isotropy, e.g.\ $E \Gamma ( \Cyc )$. 
Note that such a space can not be locally compact. The technical difficulty
here was to find the right notion of equivariant continous control (see \ref{ccdef}).
Moreover we do the obvious generalizations necessary to deal with non-$\Gamma$-compact
spaces.

Later in Section~\ref{identify} we show that our homology theory coincides
up to a shift with the one constructed by abstract means in \cite{Davis-Lueck(1998)}
and in particular the assembly map can be identified up to a shift with the map induced
on homotopy groups by
\[
\IK^{-\infty} \cald^{\Gamma}(E\Gamma(\calf)) \to \IK^{-\infty} \cald^{\Gamma} ( \punkt ).
\]
To invoke the geometry it is convenient to have our functor depend on a little more
than just $X$. 
We introduce the notion of a resolution of a $\Gamma$-space. This is
a free $\Gamma$-space $\overline{X}$ mapping to $X$ such that every $\Gamma$-compact
set is the image of a $\Gamma$-compact set.
A resolution $p:\overline{X} \to X$ will then be the input for the slight
generalization $\cald^{\Gamma}(\overline{X} ; p)$ and we prove that the $K$-theory of 
$\cald^{\Gamma}(\overline{X};p)$ depends only on  $X$ and not on the chosen resolution 
(see Proposition~\ref{independence-resolution}).

Suppose the universal map $E \Gamma \to E \Gamma( \calf )$ is realized as a resolution, then
a possible model for the assembly map $A_{\calf \to \All}$ is obtained by applying our construction
to the following map of resolutions
\[
\xymatrix{
E \Gamma \ar[d] \ar[r]^{\id} & E \Gamma \ar[d]^{\ast} \\
E \Gamma( \calf ) \ar[r] & \punkt.
         }
\]
Here $\ast$ is our standard notation for the map to a point.
Resolutions should not be confused with Quinn's simplicially stratified fibrations
which serve as the input for the homology theories considered in \cite[Appendix]{Quinn(1982)}
and \cite{Farrell-Jones(1993)}.
Finally in Section~\ref{lower} we use some standard tricks involving cone and suspension
rings to reduce Theorem~\ref{maintheorem} to the case $n=1$. 
This heavily depends on the fact that we have no assumption on our coefficient ring $R$.

{\bf Second part.}
The second  part consisting of Sections~\ref{further} to \ref{sec_injectivity} deals with the proof of the Injectivity 
Theorem~\ref{maintheoreminjectivity}.

The negative curvature assumption
allows us to construct a specific model for the universal map $E \Gamma \to E \Gamma ( \Cyc )$, see
Subsection~\ref{modelforuniversalmap}.
Roughly speaking $E\Gamma$ is realized as a fattened version of the universal covering $\tilde{M}$ 
of the compact negatively curved Riemannian manifold and the map collapses to points those
geodesics which correspond to closed geodesics down  in $M$.

The setwise stabilizer of such a lifted closed geodesic  is an infinite cyclic subgroup of $\Gamma$
and the set of $\Gamma$-orbits of such geodesics is in bijection with the 
the set of conjugacy classes of maximal cyclic subgroups.

Abstractly a lifted closed geodesic is just a real line $\IR$ with the standard
action of the infinite cyclic group $C$. The map
\[
\IK^{-\infty} \cald^{C}( \IR ; \id ) \to \IK^{-\infty} \cald^{C} ( \IR ; \ast )
\]
is a model for the classical assembly map for the infinite cyclic group and hence it is the inclusion
of a split direct summand by the Bass-Heller-Swan splitting.
In the preparatory Section~\ref{further}  we construct an induction isomorphism for our homology theory
which allows us to transport these splittings for the cyclic subgroups to the whole group. This also leads to 
Proposition~\ref{introlhscomputation}.

Since the classical assembly map $A_{1 \to \All}$ is injective by \cite{Carlsson-Pedersen(1995)} 
it remains to produce a map out of the right hand side of the assembly map which detects the 
Nil-groups in the splitting~\ref{introlhscomputation}. 
This is possible because of the geometric description of the assembly map.
In Section~\ref{sec_geometryI} we use the negative curvature assumption to construct for each
maximal cyclic subgroup a partial regain control map: identify $\tilde{M}$
with the normal bundle over a lifted closed geodesic and then push everything linearly
towards the zero section.
This map is not $\Gamma$-equivariant  but it is still $C$-equivariant for the corresponding 
 cyclic subgroup $C$.
Passing to this subgroup of infinite index we unfortunately loose the 
object support condition for our modules which is necessary to identify 
the $K$-theory groups as homology groups. But the resulting modified categories
without the object support condition are still good enough to detect the Nil-groups.

{\bf Third part.}
The remaining Sections~\ref{versus} to \ref{sec_surjectivity} are devoted to the 
proof of the surjectivity of the assembly map $A_{\Cyc \to \All}$ for $n=1$.

In Section~\ref{versus} we give a criterion which guarantees that elements of $K_1(R\Gamma)$ are in
the image of the assembly map (see \ref{smallisintheimage}). 
This criterion is formulated in terms of
$\e$-control and is very much in the spirit of Quinn's description of the 
assembly map for pseudoisotopies (see \cite{Quinn(1982)}) 
which was used in \cite{Farrell-Jones(1993)}. In particular, proving that an element 
is in the image of the assembly map
amounts to finding a suitable controlled representative for it, i.e.\ we have to regain control.
The transition from the description of the assembly map in terms of continuous control 
over a cone as discussed in the first part to the $\e$-control setting closely follows
\cite{Pedersen(2000)} and seems at present only possible for 
$K_1$ but not for higher $K$-theory. 
The next section discusses transfers for modules over a space: given a $\Gamma$-fiber bundle
$\tilde{p} : \tilde{E} \to \tilde{B}$ and a suitable fiber transport we transfer automorphisms 
of modules over $\tilde{B}$
to chain homotopy equivalences of chain complexes of modules over $\tilde{E}$ 
(see \ref{abstracttransfer}). This algebraic construction
mimics a pullback: we replace the fibers $\tilde{p}^{-1}(x)$ by cellular chain 
complexes of the fibers viewed as 
geometric chain complexes  over $\tilde{E}$ and tensor them with the module sitting down at $x$.

>From here on we translate the techniques developed in \cite{Farrell-Jones(1986)}, to
prove vanishing results for Whitehead groups, to our set-up:
Let $M$ be the  compact strictly negatively curved Riemannian manifold. Each fiber of the 
sphere bundle $S \tilde{M}$ of the universal covering can be naturally identified with the 
sphere at infinity and this leads to 
the asymptotic fiber transport for the sphere bundle (see Section~\ref{sec_geometry}) and hence
to a transfer for the bundle $S \tilde{M} \to \tilde{M}$. 
If we apply the geodesic flow on $S\tilde{M}$ to transfered morphisms 
we gain foliated control (see~\ref{transfer}).
The necessary geometric properties of the geodesic flow and the asymptotic fiber 
transport are collected
in Section~\ref{sec_geometry} (see in particular~\ref{nablanice} and \ref{nablacontr}). 
In fact we will work with the hyperbolic enlargement $\IH \tilde{M}$ which comes with a distinguished
direction ``north'' and hence allows us to work with the northern hemisphere-subbundle of the sphere 
bundle
which is a disk bundle and hence homotopically a trivial situation.

A foliated Control Theorem (see~\ref{foliatedcontrol}) 
analogous to Theorem~1.6\ in~\cite{Farrell-Jones(1986)} allows us to 
improve the foliated control we gained
to $\epsilon$-control away from flow lines that cover ``short'' closed geodesics in $M$.
These flow lines are stabilized by cyclic subgroups and collapsing them to 
points introduces cyclic isotropy and naturally leads  to working with $E \Gamma ( \Cyc )$
on the left hand side of the assembly map.

\subsection{Some conventions and notations}

Quite frequently in this paper the interested reader has to digest heavy loaded notation
as for example the definition of $\cald^{\Gamma}(\overline{X};p)$ in Section~\ref{functor}
\[
\cald^{\Gamma} (\overline{X};p)= 
\calc^{\Gamma} ( \overline{X} \times \einsu ; (p \times \id)^{-1} \cale_{\Gamma cc}(X) , 
                                       p_{\overline{X}}^{-1} \calf_{\Gamma c} ( \overline{X} ) )^{\infty}.
\]
Here are some explanations. 
Details can of course be found below in particular in Section~\ref{preparations}.
Our generic symbol for the category of modules and morphisms over the space $X$ is
\[
\calc(X ; \cale , \calf).
\]
Here $\cale$ specifies the morphism support conditions and $\calf$ specifies 
support conditions for the modules. 
The most important ones are the equivariant continuous control condition $\cale_{\Gamma cc}(X)$ 
for $X \times \einsu$ (see \ref{ccdef}) and the $\Gamma$-compact object support
condition $\calf_{\Gamma c}(X)$ for $X$ (see \ref{defFGammac}). Measuring control via a map
is formalized by pulling back support conditions, written e.g.\ as $p^{-1}\cale$. 
Since there is only a limited number
of $p$'s and $q$'s we decided to use the notation $p_X$ for projections with target
$X$ whenever it seemed reasonable.
The upper index $\Gamma$ indicates that we deal with $\Gamma$-invariant 
objects and morphisms (putting it in front saves a lot of brackets). 
The $\infty$-sign indicates ``germs at infinity'', compare \ref{germsatinfty}.

\subsection{Acknowledgments}

The authors would like to thank the ``Sonderforschungsbereich 478 -  Geometrische Strukturen in der
Mathematik'' in M{\"u}nster where this project was initiated for its support. The authors also wish
to thank the National Science Foundation for
its support of this research.  
We would also like to thank Wolfgang L{\"u}ck for answering many questions and in particular
for bringing the two ``teams'' together.
The first author thanks him for drawing his attention towards controlled topology 
and the Isomorphism Conjecture.
The fourth author would also like to thank Erik Pedersen and John Roe for explaining their ideas to him.



\typeout{--------------------- preparations ---------------------}


\section{Preparations} \label{preparations}

\subsection{Non-connective $K$-theory of additive categories}

Let $\cala$ be an additive category. By default we equip an
additive category with the split exact structure or equivalently with the 
structure of a Waldhausen category, where cofibrations are up to isomorphisms
inclusions of direct summands and weak equivalences are isomorphisms. 
A spectrum $\IE$ is a sequence of spaces $\IE_n$ together with structure maps $\sigma_n: \IE_n \to \Omega \IE_{n+1}$
which need not be homotopy equivalences.\label{spectradef}
A map of spectra is a sequence of maps 
which strictly commute with the structure maps in the obvious sense. Equivalence will
always mean stable weak equivalence.

In \cite{Pedersen-Weibel(1985)} the authors
construct a non-connective spectrum whose homotopy groups are the algebraic $K$-groups
of $\cala$ (including the negative $K$-groups).
An alternative construction using Waldhausen's $S_{\bullet}$-construction is described 
in \cite{Cardenas-Pedersen(1997)}. For completeness we give a 
definition in Subsection~\ref{modelforK-infty}.
We denote this spectrum by $\IK^{-\infty}( \cala )$.
In fact $\IK^{-\infty}$ is a functor from additive categories to spectra, i.e.\
an additive and hence exact functor induces a map of spectra.
We will use the following properties of this functor: \label{K-inftyproperties}
\begin{enumerate}
\item
Applied to the category  of finitely generated free $R$-modules we get a
spectrum whose homotopy groups are the Quillen $K$-groups (\cite{Quillen(1973)})
of $R$ in positive degrees and coincide with the negative $K$-groups from \cite{Bass(1968)} in
negative degrees. Hence the construction generalizes work of Gersten \cite{Gersten(1972)}
and Wagoner \cite{Wagoner(1972)}.
\item
An exact functor which is an equivalence of categories induces an equivalence.
\item
If $\cala$ is flasque, i.e.\ there exists an additive functor $S:\calc \to \calc$ together 
with a natural transformation $ID \oplus S \simeq S$, then 
$\IK^{-\infty} ( \cala )$ is contractible.
\item \label{karoubi_fibration}
If $\cala$ is a Karoubi filtration of the category 
$\calb$ (see \cite{Karoubi(1970)} or \cite{Pedersen-Weibel(1989)}), then there is a 
fibration sequence of spectra 
\[
\IK^{-\infty}( \cala ) \to \IK^{-\infty}( \calb ) \to \IK^{-\infty}( \calb/\cala  ). 
\]
\item \label{directed_union}
If $\cala = \bigcup_i \cala_i$ is a directed union of additive subcategories, then the natural map
$\colim_i \IK^{-\infty } \cala_i \to \IK^{-\infty} \cala$ is an equivalence. 
\end{enumerate} 

The functor $\IK^{-\infty}$ does not distinguish between a category and its idempotent
completion, e.g.\ the inclusion of the category of finitely generated free $R$-modules into 
the category of finitely generated projective modules induces an equivalence.

\subsection{The category of modules and morphisms over a space}
\label{modulesandmorphisms}

Let $R$ be an associative ring with unit.
Let $X$ be a set. A {\em general  $R$-module $M$ over $X$} is 
a family of free $R$-modules 
\[
M=(M_x)_{x \in X}
\]
indexed by the points of $X$. We are not interested in general $R$-modules 
but it is convenient to have this terminology. Let $X$ be a topological space.
A general $R$-module $M$ over $X$ is called {\em locally finite} 
if for every compact subset $K \subset X$
the module $\bigoplus_{x \in K} M_x$ is finitely generated.
A morphism $\phi=( \phi_{y,x} )$ from $M$ to $N$ is a family of 
$R$-linear maps 
\[
\phi_{y,x} :M_x \to N_y,
\]
such that
for a fixed $x$ the set of $y$ with $\phi_{y,x}\neq 0$ is finite and
for fixed $y$ the set of $x$ with $\phi_{y,x} \neq 0$ is finite.
Composition of morphisms is matrix multiplication, i.e.\ 
\[
(\phi \circ \psi)_{z,x} = \sum_y \phi_{z,y} \circ \psi_{y,x}.
\]
If a group $\Gamma$ acts on the set $X$ then we get an induced action on general 
modules by $(gM)_x=M_{gx}$. If $\Gamma$ acts by homeomorphisms on the topological space $X$
we get an induced action on locally finite modules. 
A module $M$ over $X$ is {\em $\Gamma$-invariant} if for all $g\in \Gamma$
we have $M_{gx}=M_x$. A morphism $\phi$ between $\Gamma$-invariant modules is {\em $\Gamma$-invariant} 
if  $\phi_{gy,gx} = \phi_{y,x}$. 

\begin{definition} \label{defcgamma}
The category of $\Gamma$-invariant locally finite $R$-modules and $\Gamma$-invariant 
morphisms will be denoted by 
\[
\calc^{\Gamma} ( X ; R ) \mbox{ or by } \calc^{\Gamma} (X)
\]
if the ring is clear from the context.
Given a locally finite $R$-module $M$ we define the support of $M$ as 
\[
\supp(M)=\{ x \in X | M_x \neq 0 \} \subset X.
\]
Given a morphism $\phi$ between modules
we define the support 
\[
\supp( \phi )= \{ (y,x) | \phi_{y,x} \neq 0 \} \subset X \times X.
\]
If $X$ is equipped with a metric $d$ we define the bound of $\phi$
as 
\[
|\phi|=\sup\{ d(x,y)| (x,y) \in \supp( \phi ) \}.
\]
\end{definition}

Given a map $f:X \to Y$ and a general $R$-module $M$ over $X$ we define
$f_{\ast}(M)$ to be the general $R$-module over $Y$ given by 
$f_{\ast}(M)_y =\bigoplus_{x \in f^{-1}(\{y\})} M_x$
for some choice of a direct sum.
In general $M$ locally finite does not imply that $f_{\ast}(M)$ is locally finite.
Given $\phi: M \to N$ there is an obvious induced morphism $f_{\ast}(\phi)$.

\begin{remark} \label{weissremark}
Since there is no canonical choice for the direct sum object 
\[
\bigoplus_{x \in f^{-1}(y)} M_x
\]
the construction is not really functorial in $X$. 
Following Weiss \cite{Weiss(2001)} the problem
can be solved as follows: Let $X$ be a free $\Gamma$-space. Define a general  equivariant
module $M$ over $X$ to be a functor from the partially ordered set of $\Gamma$-invariant
subsets of $X$ to the category of free $R$-modules, such that the natural map
\[
\bigoplus_{\Gamma s \in \Gamma \backslash S}  M( \Gamma s)\to M(S)
\]
induced by the inclusions is an isomorphism.
Every equivariant module $M$ determines an underlying invariant module by $M_x=M(\Gamma x)$.
Morphisms are defined as morphisms of the underlying invariant modules. The support of an 
equivariant module is the support of the underlying invariant module. Sending an equivariant 
module to its underlying invariant module induces an equivalence for all the relevant
categories. This modified construction is now truly functorial by setting $f_{\ast}M(S) = M(f^{-1}(S))$.
To simplify the presentation we will stick to the simple minded definition. 
\end{remark}

\begin{remark} 
In order to apply $K$-theory it is desirable to deal with small categories. To obtain
$\calc^\Gamma(X;R)$ as a small category one has to use a small category of free $R$-modules 
in \ref{weissremark}. This does of course restrict  the dimension of our modules. However,
we will never need modules of uncountable dimension.
\end{remark}

\subsection{Support conditions}

We will now define certain subcategories of 
$\calc^{\Gamma}(X;R)$ by imposing support conditions on objects
and morphisms. We formalize support conditions in the notion of a {\em coarse structure}
following \cite{Higson-Pedersen-Roe(1997)}. \label{coarsestructuredef}
For us a  coarse structure on $X$ consists of a set $\cale$ of subsets of $X \times X$ and 
a set $\calf$ of subsets of $X$ fulfilling certain conditions.
A morphism will be {\em admissible} if there exists an $E \in \cale$ which contains 
its support. An object will be {\em admissible} if there exists an $F \in \calf$
which contains its support.

Let us work out what conditions we need to impose on $\cale$
and $\calf$ in order to really obtain an additive category when we restrict to 
admissible objects and morphisms. 

\begin{enumerate}
\item
Since $\supp(\phi \circ \psi ) \subset \supp( \phi ) \circ \supp( \psi )$
we require that for $E$ and $E^{\prime}$ there exists an $E^{\prime \prime}$ with
$E \circ E^{\prime} \subset E^{\prime \prime }$. Here $\circ$ denotes the 
composition of relations.
\item
Since $\supp( \phi + \psi ) \subset \supp( \phi ) \cup \supp( \psi )$ we require
that for $E$ and $E^{\prime}$ in $\cale$ there exists an $E^{\prime \prime } \in \cale$ such that
$E \cup E^{\prime} \subset E^{\prime \prime }$.
\item 
Since $\supp( \id_M)$ is the diagonal embedding of $\supp ( M )$ in $X \times X$
we require that the diagonal $ \Delta \subset X \times X $ is contained in $\cale$.
\item
Since $\supp( M \oplus N ) = \supp( M ) \cup \supp( N )$ we require that 
for any $F$ and $F^{\prime}$ in $\calf$ there exists an $F^{\prime \prime} \in \calf$
with $F \cup F^{\prime} \subset F^{\prime \prime}$.
\end{enumerate}

It will be convenient to assume in addition that each $E \in \cale$ 
and each $F \in \calf$ is invariant under the $\Gamma$-operation, where $\Gamma$ acts diagonally on $X \times X$.
We will also assume that  each $E \in \cale$ is symmetric
considered as a relation.
If $\cale$ and $\calf$ fulfill these conditions we will denote the corresponding
subcategory of $\calc^{\Gamma}(X)$ consisting of admissible objects and morphisms by
\[
\calc^{\Gamma}( X ; \cale, \calf ).
\]

\vspace{1.5ex}
{\bf Some constructions for support conditions}

Given $(X,\cale , \calf)$
and a map $p:Y \to X$ we can pull back support conditions, i.e.\ form
\[
(Y, p^{-1}\cale,p^{-1}\calf). 
\]
Here $p^{-1} \cale$ consists of all sets
$(p \times p)^{-1}(E)$ with $E \in \cale$ and similar for $p^{-1} \calf$. 
Of course the map $p$ then respects the conditions.
Given two different morphism support conditions $\cale$ and $\cale^{\prime}$
one can impose both conditions, i.e.\ with slight abuse of notation form 
\[
\cale \cap \cale^{\prime} = \{ E \cap E^{\prime}|E\in \cale,E^{\prime}\in \cale^{\prime} \}.
\]
Similarly for object support conditions.
For example the {\em product} of two coarse spaces $(X ; \cale , \calf)$ and 
$(Y ; \cale^{\prime} , \calf^{\prime} )$ is given by 
$( X \times Y, p_X^{-1}\cale \cap p_Y^{-1}\cale^{\prime}, 
p_X^{-1} \calf \cap p_Y^{-1} \calf^{\prime})$.
Choosing $\cale$ or $\calf$ as the power set corresponds to imposing no condition.
In this case we simply omit $\cale$ or $\calf$ from the notation.

\vspace{1.5ex}
{\bf Morphism support conditions}

The most intuitive example of a morphism support condition is the metric morphism control
condition.

\begin{definition}[Metric control]
If  $(X,d)$ is a metric space then $\cale_{d}(X)$ denotes the metric
morphism control condition consisting of all sets of the form
\[
E_{\alpha}= \{ (x,y) | d(x,y)\leq \alpha \} \subset X \times X,
\]
where $\alpha$ is a positive real number.
If $X$ is a $\Gamma$-space we will always assume that the metric is 
$\Gamma$-invariant.
\end{definition}
 
The most important example of a morphism support condition for us is
the equivariant continuous control support condition $\cale_{\Gamma cc}$ on 
a $\Gamma$-space of the form
$X \times \einsu$. It generalizes the continuous control condition from
\cite{Anderson-Connolly-Ferry-Pedersen(1994)} to an equivariant setting. 
Recall that $\Gamma_x$ denotes the isotropy subgroup of the 
point $x$.

\begin{definition}[Equivariant continuous control]\label{ccdef}
Let $X$ be a topological space with a continuous $\Gamma$-action.
A subset $E \subset (X\times \einsu )^{\times 2}$ is in $\cale_{\Gamma cc}(X)$
if the following holds:
\begin{enumerate}
\item \label{ccdefcc}
For every $x\in X$, every $\Gamma_{x}$-invariant open neighborhood $U$ of $(x, \infty)$ in
$X \times \einsua$ there exists
a $\Gamma_x$-invariant  open neighborhood $V \subset U$ of $(x,\infty)$ such that
\[
(U ^c \times V)   \cap E = \emptyset
\mbox{ and } 
(V \times U^c )\cap E = \emptyset.
\]
Here $U^c$ denotes the complement of $U$ in $X \times \einsua$.
\item \label{ccdefmetric}
We require that $p_{\einsu} \times p_{\einsu} ( E ) \in \cale_d( \einsu )$, where $p_{\einsu}$
is the projection onto $\einsu$ and $d$ is the 
standard Euclidean metric on $\einsu$.
\item \label{ccdefequ}
$E$ is symmetric and invariant under the diagonal operation of $\Gamma$.
\end{enumerate}
\end{definition}

One checks that $\cale_{\Gamma cc}(X)$ is closed under finite unions and 
composition of relations and contains the diagonal of $X$. 
The following Lemma gives an idea about the nature of this control condition in the case where 
$X$ is a free $\Gamma$-space.

\begin{lemma} \label{convergence}
Let $X$ be a free $\Gamma$-space. Let $E$ be in $\cale_{\Gamma cc}(X)$.
Suppose that $((x_n , t_n) ,( x_n^{\prime} , t_n^{\prime} ))$ is a sequence of pairs
in $E$. If $(x_n , t_n)$ converges to $\overline{x}=(x, \infty)$, then also
$(x_n^{\prime}, t_n^{\prime} )$ converges to $\overline{x}$.
\end{lemma}

\begin{proof}
Note that in the free case all neighborhoods of $\overline{x}$ are $\Gamma_x$-invariant.
Choose a neighborhood $U_{\overline{x}}$ for $\overline{x}$ in 
$X \times \left[ 1, \infty \right]$. By assumption there exists
a neighborhood $V_{\overline{x}} \subset U_{\overline{x}}$ of $\overline{x}$ such 
that no pair in $E$ lies in 
$V_{\overline{x}}\times U_{\overline{x}}^c$. In particular since $(x_n , t_n)$
lies eventually in $V_{\overline{x}}$ we see that $(x_n^{\prime} , t_n^{\prime} )$ lies 
eventually in $U_{\overline{x}}$.
\end{proof}

{\bf Warning:} The analogous statement is wrong in the non-free case if we have infinite isotropy:
Let the infinite cyclic group $C$ act on $\IR^2$ by shifting in the 
$x$-direction. Let $X$ be the $C$-space obtained from $\IR^2$ by collapsing the
$x$-axis.
Now consider $x'_n=(n,1/n) \in X$. There is an open neighborhood of $(0,0)\in X$
such that $x'_n \notin U$ for all $n$. However, there is no such $C$-invariant
neighborhood. In particular there is an $E \in \cale_{\Gamma cc} (X)$
with $(((0,0),n),(x'_n,n)) \in E$ but $x'_n$ does not converge to $(0,0)$.

\vspace{1.5ex}
{\bf Object support conditions}

To later obtain constructions which are functorial for arbitrary
maps we will impose compactness conditions on objects, compare \ref{funclocfin} on
page \pageref{funclocfin}.

\begin{definition}[$\Gamma$-compact object support] \label{defFGammac}
Let $X$ be a $\Gamma$-space. A subset 
$F \subset X$ is in $\calf_{\Gamma c }(X)$ if it is $\Gamma$-compact, i.e.\ if it is 
of the form $\Gamma K$ for some compact set $K\subset X$. 
\end{definition}

Most of the time we will
pull back this $\Gamma$-compact support condition to $X \times \einsu$ via the standard
projection $p_X:X\times \einsu \to X$, i.e.\ we will work with $p_{X}^{-1} \calf_{\Gamma c} (X)$.

The following lemma explains in which way the category of $\Gamma$-invariant objects and
morphisms is related to the group ring $R \Gamma$, 
cf.~\cite[1.5]{Carlsson-Pedersen(1995)}.

\begin{lemma} \label{XversusGamma}
Let $X$ be a free $\Gamma$-space. The natural projection $p:X \to X/\Gamma$
induces an equivalence of categories
\[
\calc^{\Gamma}(X; \calf_{\Gamma c} ,R) \to \calc(X/\Gamma ; \calf_{c} , R\Gamma ).
\]
For any space $Y$ the map $q:Y \to \punkt$ induces an equivalence
\[
\calc(Y ; \calf_c ) \to \calc( \punkt ) .
\]
Consequently $\calc^{\Gamma}( X ; \calf_{\Gamma c} , R)$ is 
equivalent to the category of finitely generated free $R\Gamma$-modules.
\end{lemma}
\begin{proof}
Define a $\Gamma$-operation on $(p_{\ast}M)_{\Gamma x } = \bigoplus _{gx \in \Gamma x} M_{gx}$ 
by $h(m_{gx})=m_{hgx}$ for $(m_{gx}) \in \bigoplus M_{gx}$ and $h \in \Gamma$. 
The $R$-linear map $(p_{\ast}\phi)_{\Gamma y , \Gamma x}$ 
is given by the matrix $(\phi_{hy,gx})$. The fact that
$\phi$ is $\Gamma$-invariant translates into 
$(p_{\ast}\phi)_{\Gamma y , \Gamma x}$ being an $R\Gamma$-module
homomorphism. The functor $p_{\ast}$ is  full (surjective on morphism sets) 
and faithful (injective on morphism sets) and for every object in
$\calc( X/\Gamma ; \calf_{c} , R \Gamma)$ we find an isomorphic object of the 
form $p_{\ast}M$.
\end{proof}

\subsection{Germs at infinity} \label{germsatinfty}

Given $\calc^{\Gamma}(X \times \einsu ; \cale , \calf )$ we want to define
the category of germs at infinity 
\[
\calc^{\Gamma}( X \times \einsu ; \cale , \calf )^{\infty}.
\]
Intuitively we want to ignore everything that
happens in a finite region with respect to the $\einsu$-coordinate.
Objects remain the same but the morphism sets in the new category are defined as
quotients of the old ones:
\[
\mor_{\calc^{\infty}}(M,N) = \mor_{\calc}(M,N)/\sim .
\]
Here morphisms $f$ and $g$ are identified if the difference $f-g$ factors through 
some object $L$ with $\supp( L) \subset X \times \left[ 1 , n \right]$ for some finite $n$.
To really obtain a category we need some mild conditions on $\cale$ and $\calf$ that are satisfied
in all our examples. Namely we want that 
 $\calc^{\Gamma}( X \times \einsu ; \cale , \calf )^{\infty}$ is a Karoubi quotient
of $\calc^{\Gamma}( X \times \einsu ; \cale , \calf )$ 
and in particular an additive category (compare the beginning of section~\ref{coarsetools}).
As explained in \ref{standardkaroubi} there is a fibration sequence of 
spectra
\[
\IK^{-\infty} \calc^{\Gamma}( X ; \cale , \calf ) \to 
\IK^{-\infty} \calc^{\Gamma}( X \times \einsu ; \cale , \calf )
\to \IK^{-\infty} \calc^{\Gamma}( X \times \einsu ; \cale , \calf )^{\infty}.
\]
This sequence is used frequently to obtain results about the germ category.

\subsection{A model for the  non-connective algebraic $K$-theory spectrum}
\label{modelforK-infty}

For completeness we now give the definition of the model for the non-connective algebraic 
$K$-theory spectrum of an additive category which we will use throughout 
the paper. For more details see \cite{Pedersen-Weibel(1985)} and \cite{Cardenas-Pedersen(1997)}. 
Note that the definitions in the previous subsections make sense for any (small) additive category 
instead of the category of finitely generated free modules.
In particular for any additive category $\cala$ one can form $\calc( \IR ; \cale_d ; \cala)$
where $d$ denotes the standard metric on $\IR$.
Let $|iS_{\bullet} \cala|$ be  the realization of the nerve of the category of isomorphisms
in Waldhausen's $S_{\bullet}$-construction applied to $\cala$. Set 
$K( \cala )=\Omega |iS_{\bullet} \cala |$ and define the $n$-th space of the spectrum $\IK^{-\infty}\cala$
via an $n$-fold iteration as
\[
\IK^{-\infty} \cala _n = 
K( \calc( \IR ; \cale_d ;\calc( \IR ; \cale_d ; \dots \calc( \IR ; \cale_d ; \cala ) \dots ))),
\] 
where there are $n$ occurrences of $\IR$. The 
structure maps are constructed from the square associated to 
the decomposition 
$\IR = (-\infty , 0 ] \cup [0 , \infty)$ using the distinguished null-homotopies induced
by the Eilenberg-Swindle on the two half-lines. Compare the last pages in \cite{Cardenas-Pedersen(1997)}.
Note that in general the structure maps are not weak equivalences.


\typeout{--------------------- functoriality ---------------------}


\section{The functor to additive categories} \label{functor}

\subsection{Resolutions}

For the isomorphism conjecture it is essential to deal with non-free $\Gamma$-spaces.
But the translation from free $R\Gamma$-modules to $R$-modules over a $\Gamma$-space 
as in Lemma~\ref{XversusGamma} works only in the free case.
Fortunately we can separate the space which leads to the morphism support conditions 
from the actual space on which the modules live and we can 
always arrange that the latter is free. The following definition generalizes an idea from
\cite{Carlsson-Pedersen-Roe(2001)}.

\begin{definition}
A {\em resolution} of the $\Gamma$-space $X$ is a free $\Gamma$-space $\overline{X}$ together 
with an equivariant continuous map $p:\overline{X} \to X$ which fulfills the following condition:
for every $\Gamma$-compact set $\Gamma K \subset X$ we can find a 
$\Gamma$-compact set $\Gamma \overline{K} \subset \overline{X}$ such that 
$p(\Gamma  \overline{K} )=\Gamma K$. In addition we require the $\Gamma$-action on $\overline{X}$
to be properly discontinuous and $\Gamma \backslash \overline{X}$ to be Hausdorff.
(This is in particular fulfilled whenever $\overline{X}$ is a free $\Gamma$-CW complex.) 
\end{definition}

Note that $p$ is necessarily surjective.
It is immediately clear that the standard projection $p_X:X \times \Gamma \to X$ is a resolution of $X$.
We call this the {\em standard resolution of $X$}. For a space $Y$ we denote the unique map
$Y \to \punkt$ by $\ast$. Note that this is a resolution if $Y$ is a free $\Gamma$-CW complex.

Suppose we are given a map $f:X \to Y$ between $\Gamma$-spaces. Let $p:\overline{X} \to X$
and $q:\overline{Y} \to Y$ be resolutions. A map of resolutions covering $f$ is of course
\label{mapofresolutions}
a continuous equivariant map $\overline{f}: \overline{X} \to \overline{Y}$ such that the square 
\[
\xymatrix{
\overline{X} \ar[d]_{p} \ar[r]^{\overline{f}} & \overline{Y} \ar[d]^{q} \\
X \ar[r]^f & Y
         }
\]
commutes. Note that the standard resolution $p_X: X \times \Gamma \to X$ is functorial in $X$.

\subsection{The category $\cald^\Gamma(X)$} \label{definitions}

The following definitions are crucial for all that follows.
Given any $\Gamma$-space $X$ and a resolution $p:\overline{X} \to X$ we define 
\[
\calc^{\Gamma} (\overline{X};p)= 
\calc^{\Gamma} ( \overline{X} \times \einsu ; (p \times \id)^{-1} \cale_{\Gamma cc}(X) , 
                                       p_{\overline{X}}^{-1} \calf_{\Gamma c} ( \overline{X} ) ) 
\]
In words: we take the the equivariant continuous control morphism support condition 
on $X \times \einsu$ and pull it 
back to $\overline{X} \times \einsu$. 
Objects live over sets of the form $\Gamma K \times \einsu$ with $K \subset \overline{X}$ compact.
Moreover all modules and morphisms are $\Gamma$-invariant. 

Finally we pass to germs at infinity and introduce the notation
\[
\cald^{\Gamma}(\overline{X};p)=\calc^{\Gamma}( \overline{X};p )^{\infty}.
\]
If we work with the standard resolution $p_X : X \times \Gamma \to X$ we write even shorter
\[
\cald^{\Gamma} (X) = \cald^{\Gamma} ( X \times \Gamma ; p_X ).
\]
Our aim is now to prove that maps of resolutions induce functors between the corresponding 
categories. In particular we will see that $\cald^{\Gamma} ( X \times \Gamma ; p_X )$ is functorial
in $X$. 
Towards the end of this section we will show that up to an equivalence $\cald^{\Gamma}(\overline{X};p)$
does not depend on the chosen resolution of $X$. 
After some preparations in Section~\ref{coarsetools} we show in Section~\ref{homology} that
$\IK^{-\infty} \cald^{\Gamma} ( X )$ is an equivariant homology theory.

\subsection{Functoriality}

Let us first discuss abstractly what we have to do.
Let $f:X \to Y$ be given. Recall that for a general $R$-module $M$ over $X$ we defined $f_{\ast}(M)$
by $f_{\ast}(M)_y =\bigoplus_{x \in f^{-1}(y)} M_x$.
To obtain a functor 
$f_{\ast}:\calc(X; \cale , \calf ) \to \calc( Y ; \cale^{\prime} , \calf^{\prime} )$
we have to check the following: \label{abstractfunctoriality}
\begin{enumerate}
\item \label{funclocfin}
For every compact $K \subset Y$ and every $F \in \calf$ 
the set $f^{-1}( K  )\cap F$ is contained in a  compact set.
This implies that $f_{\ast}(M)$ is again locally finite.
\item
The map $f$ respects support conditions on morphisms, i.e.\
for every $E \in \cale$ there exists an $E^{\prime} \in \cale^{\prime}$ such that
$(f \times f)( E ) \subset E^{\prime}$. Sometimes one can only prove the weaker statement
that for every $E\in \cale$ and every $F \in \calf$ there exists an $E^{\prime} \in \cale^{\prime}$
with
$(f \times f)( E \cap F \times F ) \subset E^{\prime}$. This is also sufficient.
\item
The map $f$ respects support conditions on objects, i.e.\
for every $F \in \calf$ there exists an $F^{\prime} \in \calf^{\prime}$
such that $f(F) \subset F^{\prime}$.
\end{enumerate}

We will now prove that our construction is functorial with respect to maps
of resolution. The proof depends on the
slice theorem for $\Gamma$-CW complexes.
We state it at the end of this section.

\begin{proposition}
Let $f:X \to Y$ be an equivariant map between $\Gamma$-CW complexes. 
Let $\overline{f}$ be a map between resolutions $p:\overline{X} \to X$
and $q: \overline{Y} \to Y$ covering $f$, then $\overline{f}$
induces a functor
\[
\overline{f}_{\ast} : \cald^{\Gamma} ( \overline{X} ; p ) \to \cald^{\Gamma} ( \overline{Y}; q ).
\]
\end{proposition}

\begin{proof}
Let $M$ be an object of $\cald^{\Gamma}( \overline{X} ;p)$. The support condition
implies that there exist $\Gamma$-compact subcomplexes 
$\overline{X_0} \subset \overline{X}$ and $X_0 \subset X$ with 
$p(\overline{X_0}) \subset X_0$ such that $M$ is in the image of the inclusion
$\cald^{\Gamma}( \overline{X_0} ; p|_{\overline{X_0}}) 
\to \cald^{\Gamma}( \overline{X};p)$. Therefore we may  assume that $\overline{X}$
and $X$ are $\Gamma$-compact. We have to check the conditions (i), (ii) and (iii) above.
By definition $\Gamma$ acts properly discontinuously on $\overline{Y}$ and $\Gamma\backslash\overline{Y}$
is Hausdorff. This can be used to deduce (i).
The second condition  follows from
Lemma~\ref{cclemma} and (iii) is obvious.
\end{proof}

\begin{lemma} \label{cclemma}
Let $f:X \to Y$ be an equivariant continuous  map 
between $\Gamma$-CW complexes. Assume that $X$ is $\Gamma$-compact. Then
$f \times \id_{\einsu}$ respects the equivariant continuous control
morphism support conditions.
\end{lemma}

\begin{proof}
Let $E\in \cale_{\Gamma cc}$ be given.
That $(f \times \id_{\einsu })^{\times 2} (E)$ satisfies the metric 
condition \ref{ccdefmetric} in Definition \ref{ccdef} and 
is $\Gamma$-invariant and symmetric is clear. Given $y \in Y$, a $\Gamma_y$-invariant 
neighborhood $U$ of $y$ and an $r \geq 1$ we have to find a $\Gamma_y$-invariant 
neighborhood $V \subset U$ of $y$ and an $R \geq r$ such that the 
condition \ref{ccdef} \ref{ccdefcc} is satisfied for $U \times (r,\infty]$ and $V \times (R,\infty]$. 
Suppose such a neighborhood $V$
together with an $R$ does not exist. Let $V^n$ be a descending sequence of open
slice neighborhoods of $y$ as in Proposition~\ref{slice} below. It follows that 
there exists a sequence $(x_n , t_n , x_n^{\prime} , t_n^{\prime}) \in E$ such that 
for all $n \in \IN$ we have 
\begin{eqnarray*} 
(f(x_n),t_n) & \in & V^n \times (n , \infty ) \\
\mbox{ but } \quad 
(f(x_n^{\prime}) , t^{\prime}_n) & \notin &  U \times ( r, \infty ).
\end{eqnarray*}
Note that $t_n$ tends to infinity.
Since $X$ is $\Gamma$-compact there exists a sequence $g_n \in \Gamma$ such that 
$g_n x_n$ converges to some $x \in X$. For all $n \geq m$ we have
\[
f(g_n x_n)= g_n f(x_n) \in g_n V^n \subset g_n V^m \subset \Gamma V^m
\]
and hence for the limit we have $f(x) \in \overline{ \Gamma V^m}$. Since this holds for every
$m$ the last sentence in Proposition~\ref{slice} gives us
\[
f(x) \in \bigcap_{m \geq 1} \overline{ \Gamma V^m }= \Gamma \{ y \}.
\]
In particular $gf(x)=y$ for some $g \in \Gamma$ and $gg_n f(x_n)$ tends to $gf(x)=y$.
For $n$ large enough we have therefore $gg_n f(x) \in V^1$. Since $f(x) \in V^1$
it follows from \ref{slice} \ref{sliceslice} that $gg_n \in \Gamma_y$ for large $n$.
Since $U$ and $V^n$ are $\Gamma_y$-invariant we conclude that for large $n$ 
\begin{eqnarray} \label{eqn1}
(f(gg_n x_n) , t_n ) & \in & V^n \times (n , \infty ) \\
\mbox{ but } \quad (f(gg_n x_n^{\prime}) , t_n^{\prime}) & \notin & U \times (r , \infty).\label{eqn2}
\end{eqnarray}
Let us now use that $E \in \cale_{\Gamma cc} (X)$. 
There exists a $\Gamma_{gx}$-invariant
neighborhood $W$ of $gx$ and an $R \geq r$ such that 
\begin{eqnarray}\label{eqn3}
E \cap \left(  ( W \times (R , \infty )) \times (f^{-1}(U) \times (r, \infty) )^c \right) = \emptyset
\end{eqnarray}
Since $gg_n x_n $ tends to $gx$ and $t_n$ tends to infinity there exists an $n_0$
such that $(gg_{n_0} x_{n_0} , t_{n_0} ) \in W \times (R, \infty)$. Since $E$ is $\Gamma$-invariant
$(x_{n_0} , t_{n_0} , x_{n_0}^{\prime} , t_{n_0}^{\prime} ) \in E$ implies
$(gg_{n_0} x_{n_0} , t_{n_0} , gg_{n_0} x_{n_0}^{\prime} , t_{n_0}^{\prime} ) \in E$.
Now (\ref{eqn3}) implies $(gg_{n_0} x_{n_0}^{\prime} , t_{n_0}^{\prime} ) \in f^{-1}(U) \times (r, \infty)$.
This contradicts (\ref{eqn2}).
\end{proof}

\begin{proposition}[Slice theorem]\label{slice}
Let $X$ be a $\Gamma$-CW complex. For every $x\in X$ there exists a
sequence of open neighborhoods $V^{k}$, $k \in \IN$ \ with the following properties.
\begin{enumerate}
\item \label{sliceinvariant}
Each $V^k$ is an open $\Gamma_x$-invariant neighborhood of $x$.
\item \label{sliceslice}
The obvious map 
\[
 \Gamma \times_{\Gamma_x} V^k \to \Gamma V^k
\]
is a homeomorphism. \\

(A neighborhood with these properties is called an open slice around $x$. In particular 
for such a slice we have that $g V^k \cap V^k = V^k$ for every $g \in \Gamma_x$ and 
$gV^k \cap V^k = \emptyset$ for $g \notin \Gamma_x$.)
\item \label{slicedescending}
The sequence is descending, i.e.\
$V^1 \supset V^2 \supset V^3 \dots $ and 
\[ 
\bigcap_{k \geq 1} \overline{\Gamma V^k} = \Gamma \{ x \}.
\]
\end{enumerate}
\end{proposition}

\begin{proof} 
The proof  is a slight modification of the proof of Theorem~1.37 in \cite{Lueck(1989)}.
The compact isotropy groups required there are not necessary in the case of a discrete
group.
\end{proof}

\subsection{Independence of the chosen resolution}

We will now prove that up to an equivalence of categories $\cald^{\Gamma}( \overline{X} ; p )$
only depends on $X$ and not on the resolution $p: \overline{X} \to X$.

\begin{proposition}  \label{independence-resolution}
Let $p:\overline{X} \to X$ and $p^{\prime}:\overline{X}^{\prime} \to X$ be two resolutions
of $X$.
If in the diagram
\[
\xymatrix{
Z \ar[d]_{q^{\prime}} \ar[r]^q & \overline{X} \ar[d]^p \\
\overline{X}^{\prime} \ar[r]_{p^{\prime}} & X
         }
\]
the space $Z$ is defined as the pullback of $\overline{X} \to X \leftarrow \overline{X}^{\prime}$,
then $p \circ q = p^{\prime} \circ q^{\prime}: Z \to X$ is a resolution of $X$.
The maps $q$ and $q^{\prime}$ induce equivalences of categories
\[
\xymatrix{
\cald^{\Gamma}( \overline{X} ; p)& \ar[l]_{q_{\ast}} 
\cald^{\Gamma}( Z ; p\circ q ) \ar[r]^{q^{\prime}_{\ast}} &
\cald^{\Gamma}( \overline{X}^{\prime} ; p^{\prime} ).
         }
\]
\end{proposition}

\begin{proof}
We first show that $p \circ q: Z \to X$ is a resolution. Recall that 
$Z = (p \x p^{\prime})^{-1} \left(\{(x,x)|x\in X\} \right) \subset \overline{X} \x \overline{X}^{\prime}$.
Let $K$ be a compact subset of $X$.
By assumption we find compact subsets $\overline{K} \subset \overline{X}$ and 
$\overline{K}^{\prime} \subset \overline{X}^{\prime}$ with 
$p(\Gamma \overline{K})=p^{\prime}(\Gamma \overline{K}^{\prime})=\Gamma K$.
Now $q^{-1}(\overline{K}) \cap q^{\prime -1}(\overline{K}^{\prime}) 
= \overline{K} \times \overline{K}^{\prime} \cap Z \subset \overline{X} \times \overline{X}^{\prime}$
is the compact set we are looking for. 
It remains to be verified that the functor induced by $q$ is an equivalence. 
It is faithful, which is always the case even if we map 
everything to a point. It is full 
because we pulled back the morphism support condition
via $q$. Every object in the target category is isomorphic to an image object
under the functor: Simply choose a preimage for every $\Gamma$-orbit. The resolution
condition assures that one can do this in such a way that the choice respects the 
support condition.
\end{proof}


\typeout{---------------------coarse ---------------------}


\section{Some coarse tools} \label{coarsetools}

In this section we briefly collect the more or less obvious generalizations
of statements about the passage from coarse spaces to spectra which can be found for example
in \cite{Higson-Pedersen-Roe(1997)}. All proofs are straightforward generalizations
of the proofs given there. 

Let $(X, \cale , \calf ) $ be a coarse $\Gamma$-space, compare page~\pageref{coarsestructuredef}.
Given a $\Gamma$-invariant 
subset $A \subset X$ and an $E \in \cale$ define the {\em $E$-neighborhood of $A$ in $X$} 
(or the {\em $E$-thickening}) as \label{Ethickening}
\[
A^E = \{ x \in X| \mbox{ there exists } a \in A \mbox{ with } (a,x) \in E \} .
\] 
Note that $A^E$ is $\Gamma$-invariant. Let
\[
\calf \langle A \rangle_\cale = \{ (F \cap A)^E \cap F | F \in \calf, E \in \cale \}.
\]
We write $\calf \langle A \rangle$ if $\cale$ is clear  from the context.
Note that $\calc^\Gamma(X;\cale,\calf \langle A \rangle ))$ is a full subcategory of 
$\calc^\Gamma(X;\cale,\calf)$. 
The notion of a Karoubi filtration
goes back to \cite{Karoubi(1970)} and is also explained in \cite{Cardenas-Pedersen(1997)}.

\begin{proposition}[Coarse Pair] \label{standardkaroubi}
Let $(X ; \cale, \calf )$ be a coarse $\Gamma$-space and $A$ be $\Gamma$-invariant subset.
Let $i:A \to X$ be the inclusion.
\begin{enumerate}
\item
The inclusion $\calc^{\Gamma}(X; \cale , \calf  \langle A \rangle ) \to \calc^{\Gamma} ( X ; \cale , \calf )$
is a Karoubi filtration. Let $\calc^{\Gamma} ( X,A; \cale , \calf)$ denote 
the corresponding quotient category.
\item
The inclusion 
$\calc^{\Gamma}( A ; i^{-1}\cale , i^{-1}\calf )  \to \calc^{\Gamma} ( X ; \cale , \calf  \langle A \rangle)$
is an equivalence of categories.
\item
There is a fibration sequence of spectra
\[
\IK^{-\infty}\calc^{\Gamma}( A ;  i^{-1}\cale , i^{-1}\calf)  
\to  \IK^{-\infty} \calc^{\Gamma} ( X ; \cale , \calf)
\to \IK^{-\infty} \calc^{\Gamma} ( X,A ; \cale , \calf).
\]
\end{enumerate}
\end{proposition}

\begin{proof}
The first assertion is immediate from the definitions. It is straightforward to see that the
functor in (ii) is full and faithful. Also every object in the target category is isomorphic to
an image object. This implies (ii). 
Finally, (iii) is a consequence of property 
\ref{karoubi_fibration} on page~\pageref{K-inftyproperties}.
\end{proof}

Suppose the space is $\overline{X} \times \einsu$ and $A=\overline{X} \times \{ 1 \}$,
where $p:\overline{X} \to X$ is a resolution. Consider $\cale = (p \x \id)^{-1}\cale_{\Gamma cc}$
and $\calf = p_{\overline{X}}^{-1}\calf_{\Gamma c}(\overline{X})$.
Then $\calf \langle A \rangle$ and $\calf' = \{ (\overline{X} \x [1,n]) \cap F | F \in \calf, n \geq 1 \}$
are equivalent in the following sense: every $S \in \calf \langle A \rangle$ is contained in some 
$S' \in \calf'$ and vice versa.   
We are therefore in the 
situation of Subsection~\ref{germsatinfty} and we write 
$\calc^{\Gamma}( \overline{X} \times \einsu; \cale , \calf)^{\infty}$
instead of $\calc^{\Gamma}( \overline{X} \times \einsu, \overline{X} \times \{1\} ; \cale,\calf)$.

The next proposition is the equivariant version of the coarse Mayer-Vietoris principle:

\begin{proposition}[Coarse Mayer-Vietoris] \label{CMV}
Let $A$ and $B$ be  $\Gamma$-invariant subsets of the coarse $\Gamma$-space 
$(X, \cale , \calf)$ such that $A \cup B = X$. 
Suppose that the triple $(X,A,B)$ is coarsely excisive in the following sense:
for every $E \in \cale$ and $F \in \calf$ 
there exist $E' \in \cale$ and $F' \in \calf$ such that 
$(A \cap F)^E \cap (B \cap F)^E \cap F \subset (A \cap B \cap F')^{E'} \cap F'$.
Then the following diagram induces a  
homotopy pull-back square after applying $\IK^{-\infty}$.
\[
\xymatrix{
\calc^{\Gamma}( A \cap B ; j^{-1}\cale , j^{-1}\calf ) \ar[r] \ar[d] &
\calc^{\Gamma}( A ; i_A^{-1}\cale , i_A^{-1}\calf ) \ar[d] \\
\calc^{\Gamma}( B ; i_B^{-1}\cale , i_B^{-1}\calf ) \ar[r] &
\calc^{\Gamma}( X ; \cale , \calf ).
         }
\]
Here $i_A$, $i_B$ and $j$ are the obvious inclusions.
\end{proposition}

\begin{proof}
This can be proven exactly as \cite[9.3]{Higson-Pedersen-Roe(1997)}:
One checks that the inclusion 
$\calc(A,A \cap B;i_A^{-1}\cale,i_A^{-1}\calf) \to \calc(X,B;\cale,\calf)$
is an equivalence of categories and uses \ref{standardkaroubi}. 
\end{proof}

The next proposition gives a criterion for a category of the type
$\calc^{\Gamma}(X ; \cale , \calf )$ to be flasque.

\begin{proposition}[Eilenberg Swindle]  \label{eilenbergswindle}
The category $\calc^{\Gamma}( X ; \cale , \calf )$ is flasque if
there exists a $\Gamma$-equivariant self map $s:X \to X$ with the following
properties:
\begin{enumerate}
\item
For every compact $K \subset X$, every $F \in \calf$ and every $n \geq 0$
the set $(s^n)^{-1}(K) \cap F$ is compact and eventually empty.
\item
For every $E \in \cale$ and $F \in \calf$ there exists an $E^{\prime}$
with 
\[
\bigcup_{n \geq 1} (s \times s)^n( E \cap F \times F) \subset E^{\prime}.
\]
\item
For every $F \in \calf$ there exists an $F^{\prime} \in \calf$ with
$\bigcup_{n \geq 1} s^n( F ) \subset F^{\prime}$.
\item
For every $F \in \calf $ the set $\{ (x,s(x)) | x \in \bigcup_{n \geq 0} s^n( F ) \}$
is contained in some $E \in \cale$.
\end{enumerate}
In many applications (ii)-(iv) follow because for every $E \in \cale$ one has 
$(s \times s) (E) \subset E$ and for every $F \in \calf$ one has $s(F) \subset F$.
\end{proposition}

\begin{proof}
There is a functor $T=\bigoplus_{n \geq 1} (s^n)_{\ast}$ to general modules over
$X$ and  a natural isomorphism of functors $\tau: ID \oplus T \to T$.
One only needs to assure that every $T(M)$ is again locally finite
and that $T(M)$, $T(\phi)$ and $\tau_M$ respect the support conditions.
The first condition is responsible for locally finiteness.
The second checks that $\supp( T(M) )$ is admissible. The third
that $\supp( T( \phi ))$ is admissible and the last condition deals
with $\supp( \tau_M )$.
\end{proof}

\begin{example}[Swindle towards infinity] \label{swindletoinfty}
Consider $\einsu$ with the metric coarse structure $\cale_d$. Then 
$\calc(\einsu;\cale_d)$ is flasque by Proposition~\ref{eilenbergswindle}. Indeed,
$s(t) = t+1$ satisfies the above hypothesis.
\end{example}


\typeout{---------------------homology ---------------------}


\section{An equivariant homology theory}
\label{homology}

In this section we will show that the functor $\IK^{-\infty} \cald^{\Gamma}( - )$ from
the category of $\Gamma$-CW complexes and equivariant continuous maps to the category of 
spectra and maps (compare page~\pageref{spectradef}) is a generalized equivariant 
homology theory in the following sense: It satisfies
\begin{enumerate}
\item
{\bf Homotopy invariance.} \label{homotopyinvariance}
An equivariant homotopy equivalence $f:X \to Y$ induces an equivalence of spectra
\[
\IK^{-\infty} \cald^{\Gamma}(X) \to \IK^{- \infty} \cald^{\Gamma} (Y).
\]
\item
{\bf Pushout-pullback property.} \label{pushoutpullback}
A homotopy pushout square in the category of $\Gamma$-CW complexes
induces a homotopy pullback square of spectra.
\item
{\bf Directed unions.} \label{colimproperty} 
Suppose the $\Gamma$-CW complex $X$ is the directed union of subcomplexes $X_i$, then 
the natural map
\[
\colim_{i \in I} \IK^{-\infty} \cald^{\Gamma} ( X_i ) \to \IK^{-\infty} \cald^{\Gamma}( X)
\]
is an equivalence.
\end{enumerate}

The last two properties may be replaced by the maybe more familiar 
\begin{enumerate}
\item[(ii)']
{\bf Exactness.}
A cofibration sequence of $\Gamma$-CW complexes induces a homotopy
fibration sequence of spectra.
\item[(iii)']
{\bf Disjoint union axiom.}  \label{disjointunion}
Given a  family $X_i$, $i \in I$ of $\Gamma$-CW complexes 
the natural map 
\[
\bigvee_{i \in I} \IK^{-\infty} \cald^{\Gamma} ( X_i) 
\to \IK^{-\infty} \cald^{\Gamma} ( \coprod_{i \in I} X_i )
\]
is an equivalence of spectra. 
\end{enumerate}

Consequently the homotopy groups 
\[
K_n( \cald^{\Gamma}(X) ) = \pi_n ( \IK^{-\infty} \cald^{\Gamma} ( X ) )
\]
give a generalized equivariant homology theory in the usual sense.
Compare for example Definition 6.7
on page~144 in \cite{tomDieck(1987)}.
Note that the exactness axiom  corresponds to a ``soft'' excision in comparison to a
``hard'' excision statement involving arbitrary closed subsets. Compare the discussion in \cite{Weiss(2001)}.
In fact it is not clear whether an analogue of Proposition~\ref{CMVassumption} holds without
the extra security zone $U$.
Note also that we have no properness restrictions or the like on maps.
We have a usual homology theory, not a locally finite one.

Let us briefly discuss relative versions. 
\begin{notation} \label{pairnotation}
Let $p : \overline{X} \to X$ be a resolution and 
$\overline{A} \subset \overline{X}$ be $\Gamma$-invariant. Let
\[
  \cald^\Gamma(\overline{X},\overline{A};p) =
  \calc^\Gamma(\overline{X} \x \einsu ,\overline{A} \x \einsu;
               (p \times \id)^{-1} \cale_{\Gamma cc}(X) ,
                p_{\overline{X}}^{-1} \calf_{\Gamma c} ( \overline{X} ) ),
\]
be the Karoubi quotient in the notation of Proposition~\ref{standardkaroubi}(i).
It is the quotient of $\cald^{\Gamma}(\overline{X};p)$ by the thickened version
of $\cald^{\Gamma}(\overline{A};p|_{\overline{A}} )$, compare~\ref{standardkaroubi}(ii).
Of course we write 
\[
\cald^{\Gamma}(X,A)
\]
if we work with the standard resolution.
\end{notation}

Given a homology theory $\IH(-)$ in the sense above one can always extend it to 
a functor on pairs by defining $\IH(X,A)$ to be the functor evaluated on the 
cone of the inclusion $i:A \to X$.
It follows from Proposition~\ref{standardkaroubi} and the properties
of a homology theory that the two possible interpretations
of the third term in the following fibration sequence are  
naturally equivalent.
\[
\IK^{-\infty} \cald^\Gamma(A) \to 
\IK^{-\infty} \cald^\Gamma(X) \to \IK^{-\infty} \cald^\Gamma(X,A).
\]
\begin{remark} \label{resolvedversion}
Note that combined with Proposition~\ref{independence-resolution}
we obtain ``resolved'' versions of all the statements above for arbitrary instead of standard resolutions. 
For example
a map of resolutions $(\overline{f},f):(\overline{X},p) \to (\overline{Y},q)$ where
$f:X \to Y$ is an equivariant  homotopy equivalence induces an
equivalence
\begin{eqnarray} \label{homotopyplusindependence}
\xymatrix{
\cald^{\Gamma}( \overline{X};p) \ar[r]^{\simeq} & \cald^{\Gamma}( \overline{Y};q) .
         }
\end{eqnarray}
\end{remark}

The proof of properties (i) to (iii) depends on results from Section~\ref{coarsetools}. 
In order to apply the coarse Mayer-Vietoris principle in the continuous control setting 
we will need the next two results.

\begin{proposition} \label{CMVassumption}
Let $X$ be a $\Gamma$-CW complex and $A$, $B \subset X$ be two
$\Gamma$-invariant closed subsets, where $A$ is $\Gamma$-compact. Assume
that there exists a $\Gamma$-invariant
open neighborhood $U$ of $A\cap B$ in $X$ such that
$U$ is homeomorphic to $A \cap B \times (-1,1)$, where $\Gamma$ acts
trivially on $(-1,1)$. Moreover,
$U \cap A$ should correspond to $A \cap B \x [0,1)$ and $U \cap B$ to $A
\cap B \x (-1,0]$.
Let $E$ be in $\cale_{\Gamma cc}(X)$,
then there exists an $E^{\prime} \in \cale_{\Gamma cc}(X)$, such that
\[
( A \times \einsu )^E \cap  ( B \times \einsu )^E \subset (A \cap B
\times \einsu )^{E^{\prime}}.
\]
\end{proposition}

\begin{proof}
Let $Z_A = (A \x \einsu)^E - (A \x \einsu)$ and $Z_B = (B \x \einsu)^E -
(B \x \einsu)$.
For $0<\e<1$ let $U_\e$ be the neighborhood of $A \cap B$ corresponding
to $(A \cap B) \x (-\e,\e)$.
Then our assumptions imply the following: 
for any $0<\e<1$ there is $t_\e$ such that
$Z_A \cap X \x [t_\e,\infty),Z_B \cap X \x [t_\e,\infty) \subset U_\e \x
\einsu$.
This allows us to define for $t \geq t_{1/2}$ the function
\[
\e(t)= \min \{ \epsilon | (Z_A \cup Z_B) \cap X \x [t,\infty) \subset
U_\e \x \einsu \}.
\]
For $z \in A \cap B$ let $L_\e(z) \subset U$ correspond to $\{ z \} \x
(-\e,\e)$.
Now set
\begin{eqnarray*}
E'=\{ ((x,t),(y,t)) &  | &  t<t_{1/2}
 \mbox{ or }  x=y  \\ & \mbox{ or } & x \in A \cap B \mbox{ and } y \in
L_{\e(t)}(x) \\
& \mbox{ or } &  y \in A \cap B \mbox{ and } x \in L_{\e(t)}(y) \}.
\end{eqnarray*}
It is not hard to see that this $E'$ satisfies our claim.
\end{proof}

\begin{lemma} \label{i-1lemma}
If $i:X \to Y$ is the inclusion of a closed $\Gamma$-invariant subspace, then 
$\cale_{\Gamma cc}(X) = i^{-1} \cale_{\Gamma cc} ( Y)$.
\end{lemma}

\begin{proof}
This follows straightforward from the definitions.
\end{proof}

We can now prove homotopy invariance using the coarse Mayer-Vietoris principle and an Eilenberg swindle.

\begin{proposition}  \label{homotopyinvarianceCor}
The functor $\IK^{-\infty} \cald^{\Gamma} ( - )$ is homotopy invariant.
\end{proposition}

\begin{proof} \label{read!!!H}
For $i=0,1$ let $X_i = [i,1] \x X$ and  $Z_i = [i,\infty) \x X$. 
Let $\overline{X}, \overline{X_i}$ and $\overline{Z_i}$ be standard resolutions. 
Denote by $p$ the projections $\overline{Z_i} \x \einsu \to Z_i \x \einsu$ and
$\overline{X_i} \x \einsu \to X_i \x \einsu$. 
Denote by $q$ the projections $\overline{Z_i} \x \einsu \to \overline{X}$ and
$\overline{X_i} \x \einsu \to \overline{X}$. 
Let 
\begin{eqnarray*}
 \cala(X_i)  & 
 = & 
 \calc^\Gamma(\overline{X_i} \x \einsu;p^{-1}\cale_{\Gamma cc}(X_i),q^{-1}\calf_{\Gamma c}(\overline{X_i}))
 \\
 \cala(Z_i) &  
 = & 
 \calc^\Gamma(\overline{Z_i} \x \einsu;p^{-1}\cale_{\Gamma cc}(Z_i),q^{-1}\calf_{\Gamma c}(\overline{Z_i})) 
\end{eqnarray*}
Using  \ref{CMVassumption} the conditions of \ref{CMV} can be checked.
Therefore (and by \ref{i-1lemma}) 
the  following square induces a homotopy pull-back square in $K$-theory
\[
\xymatrix
{
\cala(X_1) \ar[r] \ar[d] 
& 
\cala(X_0) \ar[d]
\\
\cala(Z_1) \ar[r]
&
\cala(Z_0).
\\
}
\]
The map $(s,x,\gamma,r) \mapsto (s+1/r,x,\gamma,r)$ induces Eilenberg swindles on $\cala(Z_0)$ and $\cala(Z_1)$.
(We want to point out that this makes use of the metric control in the $\einsu$-direction 
from Definition~\ref{ccdef} \ref{ccdefmetric}.) By  \ref{eilenbergswindle} this implies, that $\cala(X_0) \to \cala(X_1)$
is an equivalence in $K$-theory.
Note that $\cald^\Gamma(X_i)$ is obtained from $\cala(X_i)$ by taking germs at infinity. 
Homotopy invariance follows now from the fibration sequence in Subsection \ref{germsatinfty}.
\end{proof}

Property~\ref{pushoutpullback} also depends on \ref{CMV}.

\begin{proposition} \label{pushoutpullbackProp}
The functor $\IK^{-\infty} \cald^{\Gamma} ( - )$ applied to a homotopy
pushout diagram yields a homotopy-pullback diagram of spectra.
\end{proposition}

\begin{proof}
Suppose $X_1 \leftarrow X_0 \to X_1$ is given. 
By homotopy invariance it  
is sufficient to prove that the double mapping cylinder corresponding to
this diagram yields a homotopy pullback square. The triple consisting of 
the double mapping cylinder together with the individual mapping cylinders 
corresponding to $X_0 \to X_1$ respectively $X_0 \to X_2$ fulfill 
the assumptions of \ref{CMVassumption}. Using this we can verify
the condition in the Coarse Mayer-Vietoris principle~\ref{CMV}, which we can hence apply.
Lemma~\ref{i-1lemma} finishes the proof.
\end{proof}

It remains to consider the behavior of our functor applied to directed unions.

\begin{proposition}
The functor $\IK^{-\infty}\cald^{\Gamma}(-)$ satisfies (iii) on page \pageref{colimproperty}.
\end{proposition}

\begin{proof}
Suppose the $\Gamma$-CW complex $X$ is the directed union of subcomplexes $X_i$.
Observe that $\cald^{\Gamma}( X_i)$ is a full additive subcategory of $\cald^{\Gamma}(X)$.
Because of the object support condition $p_{X \times \Gamma}^{-1} (\calf_{\Gamma c})$ it follows
that $\cald^{\Gamma}(X)$ is the directed union of the subcategories $\cald^{\Gamma}(X_i)$.
Now apply property \ref{directed_union} from page~\pageref{K-inftyproperties}.
\end{proof}


\typeout{---------------------identify ---------------------}


\section{Identifying the assembly map} \label{identify}

An ordinary homology theory is determined by its behavior on a point.
In the equivariant situation the basic building blocks of a space are not points but
orbits. We will review from \cite{Davis-Lueck(1998)} that similarly
to the ordinary case
an equivariant homology theory is determined by its behavior on orbits.
This will be used to identify our version of the assembly map with other existing
assembly maps in the literature.

\subsection{Review of the Davis-L{\"u}ck construction}

Let $\Or \Gamma$-denote the orbit category, i.e.\ the category whose
objects are the $\Gamma$-sets of the form $\Gamma/H$ and whose morphisms
are $\Gamma$-maps. Note that objects are in bijection with the subgroups of $H$.
Following \cite{Davis-Lueck(1998)} we will
call a functor from the orbit category into the category of spectra (and strict maps)
an $\Or \Gamma$-spectrum. We use a question mark to indicate the place where
objects and morphisms are plugged into the functor. 
Associated to each $\Or \Gamma$-spectrum $\IK(?)$
is an equivariant homology theory in the sense of Section~\ref{homology}.
It is constructed as follows: 
every $\Gamma$-space $X$ gives rise to a contravariant $\Or \Gamma$-space 
$X^{?}=\map_{\Gamma}( \Gamma/ ? , X )$ and we can form the balanced
smash product over the orbit category between a contravariant $\Or \Gamma$-space
and a covariant $\Or \Gamma$-spectrum to obtain an ordinary
spectrum, e.g.\
\[
X^{?}_{+} \sma_{\Or \Gamma} \IK (?)= 
\left( \bigvee_{\Gamma/H \in \Or \Gamma} X^H_{+} \sma \IK ( H ) \right)/ \sim.
\]
Compare \cite[p.237]{Davis-Lueck(1998)}. 
This construction is functorial in $X$ and satisfies the properties listed
in Section~\ref{homology}. To stress the homological behavior we write
the homotopy groups of the spectrum as
\[
H_n^{\Or \Gamma}( X ; \IK ) = \pi_n( X ^{?}_{+} \sma_{\Or \Gamma} \IK( ? ) ).
\]
This notation was already used in the introduction.

Let us now recall the construction of the algebraic $K$-Theory $\Or \Gamma$-spectrum 
from \cite[Section~2]{Davis-Lueck(1998)}. We will denote it by $\IK R^{-\infty}$.
Given a $\Gamma$-set $S$ we denote by $\overline{S}$ the transport category
associated to $S$, i.e.\ the category whose objects are points in $S$
and with $\mor_{\overline{S}}(s,t)=\{ g \in \Gamma | gs=t \}$.
The transport category is a groupoid, i.e.\ every morphism is an isomorphism.
Given any small category $\calc$ we can form the associated $R$-linear category
$R\calc$ with the same objects and new morphism set 
$\mor_{R\calc}(c,d)=R \mor_{\calc}(c,d)$ the free $R$-module
generated by the old morphism set. Finally we turn such an $R$-linear category
into an additive category, i.e.\ we artificially introduce finite sums
(or products). The resulting category is denoted $R\calc_{\oplus}$.
The Davis-L{\"u}ck functor is now given as
\begin{eqnarray*}
\IK   R^{-\infty}: \Or \Gamma & \to & Spectra \\
\Gamma/H & \to & \IK^{-\infty} ( R \overline{\Gamma/H}_{\oplus}) .
\end{eqnarray*}
Note that $\IK R^{-\infty}$ and $\IK^{-\infty}( R )$ are different objects.

\subsection{Comparing the theories}
\label{comparing}

In Section~\ref{functor} we constructed the functor $\IK^{-\infty} \cald^{\Gamma}( ? )$
which sends $\Gamma$-CW complexes to spectra. In Section~\ref{homology}
we proved that it is in  fact a generalized equivariant homology theory. 
We now want to identify this homology theory (up to a shift) with the one associated to the 
$\Or \Gamma$-spectrum $\IK R^{-\infty}$ as above. In particular we will see
that our assembly map coincides with the one constructed in \cite{Davis-Lueck(1998)}.

\begin{proposition}
Suppose $\IF$ is a functor from $\Gamma$-CW complexes to spectra
which is homotopy invariant, exact and satisfies the disjoint union axiom, compare Section~\ref{homology}.
Suppose there is a zig-zag of equivalences of $\Or \Gamma$-spectra
\[
\xymatrix{
\IK R^{-\infty} \ar[r]^-{\simeq} & \IF_1 & \dots \ar[l]_-{\simeq}  \ar[r]^-{\simeq} & \IF_n=\IF  |_{\Or \Gamma}.
         }
\]
Then for every $\Gamma$-CW complex $X$ we have a zig-zag of equivalences of spectra
\[
\xymatrix{
X_{+}^{?} \sma_{\Or \Gamma} \IK R^{-\infty} ( ? ) \ar[dr]^-{\simeq} & \dots \ar[d]_-{\simeq} \ar[dr]^-{\simeq} &
\IF^{\%}( X ) \ar[d]^-{B_{X}}_-{\simeq}  \ar[dr]^-{A_{X}}_-{\simeq} \\
& X_{+}^{?} \sma_{\Or \Gamma} \IF_1 (?) &
X_{+}^{?} \sma_{\Or \Gamma} \IF|_{\Or \Gamma}(?) &
\IF ( X ) . }
\]
The diagram is natural in $X$.
\end{proposition}
\begin{proof}
The functor $\IF^{\%}$ and the natural transformations $A$ and $B$ are constructed in
\cite[6.3]{Davis-Lueck(1998)}. It is shown there that $B$ is always an equivalence and 
that $A$ is an equivalence if $\IF$ satisfies the assumptions stated in the theorem.
\end{proof}

We will now construct a zig-zag of equivalences between the restriction of 
our (shifted) functor $\Omega \IK^{-\infty} \cald^{\Gamma} ( ? )$ to the orbit category and the
Davis-L{\"u}ck functor $\IK R^{-\infty}( ? )$. 

Let $S$ be a discrete $\Gamma$ set.
Let $\cale_{\Delta}=\{ \Delta \}$ be the morphism support
condition on $S$ consisting only of the diagonal $\Delta$, i.e.\ morphisms do not move things
at all. 
On $S \times \einsu$ we define the discrete morphism support condition $\cale_{dis}$ as 
$p_{S}^{-1} \cale_{\Delta} \cap p_{\einsu}^{-1} \cale_d$, i.e.\ things do not move in the $S$-direction
and we have the standard metric control structure in the $\einsu$-direction.
As usual projections are indexed by their targets.
Recall that we write $\calc( \einsu  ; \cale_d ,\cala)$ if in the definition of
$\calc( \einsu ; \cale_d )$ we work with an additive category $\cala$ 
instead of  the category of finitely generated modules.

Now define the following functors from the orbit category to additive categories:
\begin{eqnarray*}
\cald_1( S ) &  = & 
\calc^{\Gamma}( S  \times \Gamma ; p_{S}^{-1}\cale_{\Delta} , \calf_{\Gamma c} ) \\
\cald_2( S ) & = &
\calc( \left[ 1 , \infty \right) ; \cale_d ; \cald_1( S ) )^{\infty} \\
\cald_3( S ) &  = &
\calc^{\Gamma}( S  \times \Gamma \times \einsu ; 
p_{S \times \einsu}^{-1}\cale_{dis} , p_{S \times \Gamma}^{-1} \calf_{\Gamma c} )^{\infty} 
\end{eqnarray*}
and recall the definition
\[
\cald^{\Gamma} ( S )   = 
\calc^{\Gamma}( S  \times \Gamma \times \einsu ; 
p_{S \times \einsu}^{-1}\cale_{\Gamma cc} , p_{S \times \Gamma}^{-1} \calf_{\Gamma c} )^{\infty}. 
\]

For a spectrum $\IE$ define $\Omega \IE$ to be the spectrum with 
$n$-th space $(\Omega \IE)_n=\Omega(\IE_n)$ and define $\sh \IE$ with $\sh \IE_0= \Omega \IE_0$ and $\sh \IE_n = \IE_{n-1}$. 
\begin{proposition} \label{zigzag}
We have $\cald_3( ? )=\cald^{\Gamma} ( ? )$. There is  a natural transformation
$R \overline{( ? )}_{\oplus}  \to   \cald_1( ? )$ which is objectwise an equivalence and there
is a zig-zag of equivalences of $\Or \Gamma$-spectra:
\[
\xymatrix{
\IK^{-\infty} \cald_1 ( ? )  \ar[r]^-{\simeq} & 
\sh \IK^{-\infty} \cald_2 ( ? )  & 
\sh \IK^{-\infty} \cald_3 ( ? ) \ar[l]_-{\simeq} \ar[d]^-{\simeq} \\
\IK R^{-\infty} (?) \ar[u]_-{\simeq}  &  &
\Omega \IK^{-\infty} \cald_3 ( ? ) = \Omega \IK^{-\infty} \cald^{\Gamma}( ?) 
         }
\]
\end{proposition}
\begin{proof}
(i) The equality $\cald_3=\cald^{\Gamma}$ is discussed in Remark~\ref{disremark}.

(ii) The equivalence  $R \overline{(?)}_{\oplus} \to \cald_1( ? )$ is the crucial step.
Let $S$ be a discrete $\Gamma$-set.
We describe a functor from the subcategory $\overline{S}$ of 
$R \overline{S}_{\oplus}$ to 
$\cald_1(S )= \calc^{\Gamma}( S \times \Gamma ; p_S^{-1}\cale_{\Delta}, \calf_{\Gamma c},R)$.
An object $s$ in $\overline{S}$ is sent to the equivariant module $M(s)$
with $M(s)_x = R$ if $x \in \Gamma (s,e)$ and $M(s)_x=0$ otherwise. Here $\Gamma (s,e)$
is the orbit of $(s,e)$ under the diagonal operation. A morphism 
$k \in \mor_{\overline{S}}(s,s^{\prime})$, i.e.\ a $k \in \Gamma$ with
$ks=s^{\prime}H$ is sent to the unique $\Gamma$-invariant morphism $M(k)$
with $\M(k)_{(s,k^{-1})(s,e)} = \id_R$. Note that 
$p_1\times p_1(\supp M(k) )=\Gamma  p_1 \times p_1((s,k^{-1},s,e))=\Gamma(s,s)$ lies in
the diagonal, i.e.\ fulfills the discrete morphism support condition.
By the universal properties of 
the linearization $\overline{S} \to R \overline{S}$
and of $R \overline{S} \to R \overline{S}_{\oplus}$ the functor
extends to $R \overline{S}_{\oplus}$. It induces a bijection on morphism
sets. It remains to convince oneself that in every isomorphism
class of modules in $\cald_1(S)$ there is a module of the form 
$M_{s_1} \oplus \ldots \oplus M_{s_n}$.
It is straightforward to check that the functor $M=M^S$ is natural in $S$, i.e.\ 
given an equivariant map $f:S \to T$ we have an equality 
$M^T \circ R\overline{f}_{\oplus} = (f \times \id_{\Gamma})_{\ast} \circ M^S $.

(iii) From the explicit description of the spectrum $\IK^{-\infty}$ in \ref{modelforK-infty}
it follows that the map $\IK^{-\infty} \cald_1 \to \sh \IK^{-\infty} \cald_2$ can be constructed
from the functor on the right in the ``germs at infinity sequence'', see Subsection~\ref{germsatinfty}.
\[
\calc( \left( -\infty , 1 \right] ; \cale_d ; \cald_1 ) 
\to  
\calc( \IR ; \cale_d ; \cald_1) 
\to 
\cald_2 = \calc( \einsu  ; \cale_d ; \cald_1 )^{\infty}.
\]
Since the left hand term admits an Eilenberg swindle the map induces an equivalence.

(iv) Consider the diagram
\[
\xymatrix{
\cald_1 \ar[r] \ar[d]^-{=} & 
\calc^{\Gamma}( S  \times \Gamma \times \einsu ; 
p_{S \times \einsu}^{-1}\cale_{dis} , p_{S \times \Gamma}^{-1} \calf_{\Gamma c} ) \ar[r] \ar[d]_F &
\cald_3 \ar[d] \\
\cald_1 \ar[r] & 
\calc( \einsu ; \cale_d ;  \cald_1 ) \ar[r] & 
\cald_2
         }
\]
Here the functor $F$ is given by considering a module over $S \times \Gamma \times \einsu$ as 
an object in $\calc( \einsu ; \cale_d , \cald_1 )$. Note that 
an object in $\calc( \einsu ; \cale_d , \cald_1 )$ considered as a module need not 
have support in $p_{S \times \Gamma}^{-1} \calf_{\Gamma c}$. So there is no obvious functor in the 
reverse direction. The rows are both germs at infinity sequences (compare Subsection~\ref{germsatinfty}) and 
the middle terms both admit an Eilenberg swindle. This gives the equivalence
$\sh \IK^{-\infty} \cald_3 \to \sh \IK^{-\infty} \cald_2$.

(v) Finally: for every spectrum $\IE$ there is a natural equivalence $\sh \IE \to \Omega \IE$.
\end{proof}

In particular 
we can now identify our assembly map with the one
from \cite{Davis-Lueck(1998)}.

\begin{corollary} \label{compareassemblymaps}
The map 
\[
K_{\ast + 1} \cald^{\Gamma}( E\Gamma ( \calf )) \to K_{\ast + 1} \cald^{\Gamma}( \punkt )
\]
is a model for the assembly map
\[
A_{\calf \to \All}: H_{\ast}^{\Or \Gamma} ( E\Gamma ( \calf ) ;\IK R^{-\infty} ) 
\to H_{\ast}^{\Or \Gamma} ( \punkt ;\IK R^{-\infty} ) = K_{\ast} ( R \Gamma ) .
\]
\end{corollary}


\typeout{---------------------lower ---------------------}


\section{Lower $K$-theory and Nil-groups} \label{lower}

In this section we would like to explore the fact that we do not impose any assumptions 
on our coefficient ring.

\subsection{The passage to lower $K$-theory}

To prove the isomorphism in our main Theorem~\ref{maintheorem} 
we can restrict ourselves 
to the case $n=1$, i.e. 
it suffices to prove that  
\[
H_1^{\Or \Gamma }( E \Gamma ( \calf ) ; \IK R^{-\infty} ) \to K_1( R \Gamma )
\]
is an isomorphism. This is a consequence of the fact that we
have no assumptions on the ring $R$.

Let $R$ be a ring with unit. Let $M_f( R) $ be the ring of all row and column-finite infinite
matrices over $R$, i.e.\ all matrices $(r_{ij})_{i,j \in \IN}$  such that 
for fixed $i$ only finitely many $r_{ij}$ are nonzero and similar  for fixed $j$.
Let the cone ring  $CR$ be the subring of $M_f(R)$ generated by all matrices $(r_{ij})$ for which
the set $\{ r_{ij} | i,j \in \IN \}$ is finite. Let $I(R)$ be the ideal of $CR$
consisting of all finite matrices, i.e. matrices $(r_{ij})$ such that there exists an $N$
with $r_{ij}=0$ if $i \geq N$ or $j\geq N$. Finally define the suspension ring
$\Sigma R$ as the quotient $CR/ I(R)$.

\begin{proposition} \label{suspensionfibr}
There exists a homotopy fibration sequence of $\Or \Gamma$-spectra
\[
\IK R^{-\infty}  \to \IK CR^{-\infty}  \to \IK \Sigma R^{-\infty}
\]
and $\IK CR^{-\infty}$ is a contractible $\Or \Gamma$-spectrum.
\end{proposition}

\begin{proof}
The Davis-L{\"u}ck functor $\IK R^{-\infty}$ is obviously functorial in the ring $R$.
It is sufficient to show that evaluated at any object $\Gamma/H$
we have a fibration sequence of ordinary spectra. Recall that 
$\IK R^{-\infty} ( \Gamma / H )= \IK^{-\infty} ( R \overline{ \Gamma / H }_{\oplus} )$.
Since the transport category $\overline{\Gamma / H }$ is a connected groupoid,
the inclusion of the full subcategory $\overline{eH}$ on the object $eH$ induces an equivalence.
This subcategory can be identified with the group $H$. It  is therefore sufficient to show that
associated to the group rings $R H$, $(CR) H$ and $(\Sigma R) H$ we have a fibration
sequence on the level of $K$-theory spectra. But in fact $(CR) H= C (RH)$ and $(\Sigma R)H= \Sigma (R H )$.
So we have reduced things to the case of rings, i.e. we have to prove that 
\[
\IK^{-\infty}( R_{\oplus} ) \to \IK^{-\infty}( CR_{\oplus} ) \to \IK^{-\infty}( \Sigma R_{\oplus} ) 
\]
is a fibration sequence and that $\IK^{-\infty}( CR_{\oplus} )$ is contractible.
This is a standard fact. Compare e.g.\ \cite{Wagoner(1972)}.
\end{proof}

\begin{remark}
There are several possible definitions of a cone and suspension ring leading to
a fibration sequence for ordinary algebraic $K$-theory.
Our choice has the advantage that $C(RH)=(CR)H$ and $\Sigma(RH)=(\Sigma R)H$
which was  the essential trick in the proof.
\end{remark}

\begin{corollary} \label{loweriso}
If the assembly map induces an isomorphism in the $n$-th homotopy group
for arbitrary coefficient rings then it also induces an isomorphism for all $i \leq n$.
\end{corollary}

\begin{proof}
By \ref{suspensionfibr} we obtain from
\[
X^{?}_{+} \sma_{\Or \Gamma} \IK  R^{-\infty}  ( ? ) \to
X^{?}_{+} \sma_{\Or \Gamma} \IK  CR^{-\infty}  ( ? ) \to
X^{?}_{+} \sma_{\Or \Gamma} \IK  \Sigma R^{-\infty}  ( ? ) 
\]
a long exact sequence of homotopy groups. Moreover, the homotopy groups of the middle term vanish.
The sequence is natural in $X$. Hence inserting the map $E \Gamma (\calf ) \to \punkt$ 
we can identify the assembly map 
for $R$ in dimension $i$
with the assembly map for $\Sigma R$ in dimension $i+1$. 
\end{proof}

\subsection{Assembly for Nil-groups} \label{nilsection}

We will now study the assembly map for $NK_*(R)$. Let $\IN \IK^{-\infty} R$ denote the
homotopy cofiber of the map of spectra $\IK^{-\infty} R \to \IK^{-\infty} R[t]$ induced by the inclusion of
$R$ into the polynomial ring over $R$. Thus, the homotopy groups of $\IN \IK^{-\infty} R$ are given by $NK_*(R)$.
This construction is functorial in the ring $R$ and we obtain an $\Or\Gamma$-spectrum $\IN \IK R^{-\infty}$.
For every family $\calf$ of subgroups of $\Gamma$ we have therefore an assembly map
\[
\mathit{NA}_{\calf \to \All} : 
 H^{\Or\Gamma}_*(E\Gamma(\calf);\IN \IK R^{-\infty} ) \to 
H^{\Or\Gamma}_*( \punkt; \IN \IK R^{-\infty}) = NK_*(R\Gamma). 
\]

\begin{proposition} \label{nilprinciple}
Suppose that the assembly map 
\[
H^{\Or\Gamma}_* (E\Gamma(\calf);\IK S^{-\infty}) \to 
H^{\Or\Gamma}_* (\punkt ; \IK S^{-\infty}) 
= K_*(S\Gamma)
\]
is an isomorphism for $* \leq n$ and $S=R,R[t]$. Then for the ring $R$ and $* \leq n$
the assembly map $\mathit{NA}_{\calf \to \All}$ is also an isomorphism.
In particular, if $\calf = 1$ and $R$ is regular then $NK_*(R\Gamma)=0$ for $* \leq n$.
\end{proposition}

\begin{proof}
The first statement is a consequence of the long exact sequence of homotopy groups obtained from
\[
X^{?}_{+} \sma_{\Or \Gamma} \IK  R^{-\infty}  ( ? ) \to
X^{?}_{+} \sma_{\Or \Gamma} \IK  R[t]^{-\infty}  ( ? ) \to
X^{?}_{+} \sma_{\Or \Gamma} \IN \IK R^{-\infty}  ( ? ). 
\]
The second follows from the fact that $NK_*(R)=0$ for regular rings.
\end{proof}



\typeout{---------------------further ---------------------}


\section{Induction and restriction} \label{further}

We collect some further constructions and formal properties related to the functor 
$\cald^{\Gamma}(\overline{X}, p)$ which we will need in the proof of injectivity 
in Section~\ref{injectivity}.

Let $p:\overline{X} \to X$ be a resolution. Every morphism $\phi$ in
$\calc^\Gamma(\overline{X};p)$ is given as a matrix $\phi_{(x,t),(y,s)}$. 
The following lemma is a straightforward consequence of our definitions. In particular, it is 
important here that modules in $\cald^\Gamma(\overline{X};p)$ have support contained in 
$\Gamma K \x \einsu$ for some compact $K \subset \overline{X}$.

\begin{lemma} \label{separate}
Every morphism in $\cald^\Gamma(\overline{X};p)$ can be represented by a morphism 
$\phi$ in $\calc^\Gamma(\overline{X};p)$ such that $\phi_{(x,t),(y,s)}=0$ whenever $p(x)$ and $p(y)$ 
lie in different components of $X$. 
\end{lemma}

\begin{remark} \label{disremark}
The above observation is particularly useful in the following situation:
let $p:\overline{X} \to X$ be a resolution where $X$ is a discrete space.
Let $\cale_{\Delta}=\{ \Delta \}$ be the morphism support
condition consisting only of the diagonal $\Delta$, i.e.\ morphisms do not move things
at all. As usual projections are indexed by their targets.
On $X \times \einsu$ we define the discrete morphism support condition $\cale_{dis}$ as 
$p_{X}^{-1} \cale_{\Delta} \cap p_{\einsu}^{-1} \cale_d$, i.e.\ things do not move in the $X$-direction
and we have the standard metric control structure in the $\einsu$-direction.
Now Lemma~\ref{separate} implies that we can add this to our support conditions
in $\cald^\Gamma(\overline{X};p)$ without changing our category. In symbols: If $X$ is discrete then
\[ 
\cald^\Gamma(\overline{X};p) =  
\calc^{\Gamma} ( \overline{X} \times \einsu ; (p \times \id)^{-1} \cale_{dis} , 
        p_{\overline{X}}^{-1} \calf_{\Gamma c} ( \overline{X} ) )^\infty.
\]  
\end{remark}

\begin{proposition}[Induction] \label{inddni}
Let $H \subset \Gamma$ be a subgroup.
Let $X$ be an $H$-space and let $p:\overline{X} \to X$ be a resolution.
Consider the induced resolution $\id \times_H p: \Gamma \times_H \overline{X} \to \Gamma \times_H X$.
There are inverse isomorphisms of categories
\[
\xymatrix{
\cald^{H}( \overline{X}; p ) \ar@<0.5ex>[r]^-{\ind} & 
\cald^{\Gamma}(\Gamma \times_H \overline{X} ; \id \times_H  p ) \ar@<0.5ex>[l]^-{\dni}.
         }
\]
\end{proposition}

\begin{proof} 
Let $i:X=H \times_H X \to \Gamma \times_H X$ denote the inclusion. 
For an invariant  module $M$ over $X \times \einsu$ define $(\ind M )_{\left[ g,x \right]} = M_x $.
For a morphism $\phi$ set
$(\ind \phi)_{\left[g^{\prime} , x^{\prime} \right], \left[g , x \right] } 
  = \phi_{g^{-1}g^{\prime}x^{\prime},x}$ 
if $g^{-1}g^{\prime} \in H$ and otherwise 
set $(\ind \phi)_{\left[g^{\prime} , x^{\prime} \right], \left[g , x \right] }=0$.
This defines a functor to general modules over $\Gamma \times_H X \times \einsu$.
It is straightforward to check that this also respects the support conditions
and gives a functor 
$\ind:\cald^{H}( \overline{X}; p ) \to 
 \cald^{\Gamma}(\Gamma \times_H \overline{X} ; \id \times_H  p )$.
Define an inverse by $(\dni M)_x = M_{i(x)}$ and $(\dni \phi)_{y,x} = \phi_{i(y), i(x)}$.
This does not give a functor on general modules. However, Lemma~\ref{separate} implies 
that we do get a functor 
$\dni: \cald^{\Gamma}(\Gamma \times_H \overline{X} ; \id \times_H  p )
 \to \cald^{H}( \overline{X}; p )$.
That $\dni \circ \ind$ is the identity functor is immediate from the above definitions.
Lemma~\ref{separate} implies that $\ind \circ \dni$ is also the identity. 
\end{proof}

If $H \subset \Gamma$ is a subgroup of finite index, then 
there is also an obvious restriction functor
\[
\res : \cald^{\Gamma} (\overline{X};p) \to \cald^{H} (\overline{X};p) .
\]
But the compactness conditions for objects prevents such a functor in general.
Nevertheless we will in  
Section~\ref{injectivity} need a kind of restriction functor for subgroups
which are not of finite index. We achieve this by simply dropping the object
support condition in the target category.
This leads to some technical difficulties since 
$\Gamma$-compactness arguments are not available anymore.
Fortunately in the situations we will consider $\overline{X}$ always 
comes  equipped with an invariant metric
and we can solve the problems by 
imposing  a weak metric condition on morphisms over $\overline{X} \times \einsu$.

\begin{definition} \label{deflf}
Let $p:\overline{X} \to X$ be a resolution of $\Gamma$-spaces.
Suppose $\overline{X}$ is equipped with an invariant metric $d$. 
Let $d_{prod}$ be a product metric on $\overline{X} \times \einsu$. Define
\begin{eqnarray*}
 \cald^{\Gamma}_{lf}(X;p) &   = &   \calc^{\Gamma}( \overline{X} \times \einsu; 
       p_{X \times \einsu}^{-1} \cale_{\Gamma cc})^{\infty}     \\
\tilde{\cald}^{\Gamma}(X;p) &  = & \calc^{\Gamma}( \overline{X} \times \einsu; 
   p_{X \times \einsu}^{-1} \cale_{\Gamma cc} \cap \cale_{d_{prod}},
       \calf_{\Gamma c} ( \overline{X} ))^{\infty}         \\
\tilde{\cald}^{\Gamma}_{lf}(X;p) &  = & 
  \calc^{\Gamma}( \overline{X} \times \einsu; 
       p_{X \times \einsu}^{-1} \cale_{\Gamma cc} \cap \cale_{d_{prod}})^{\infty}
\end{eqnarray*}
Short: the ``lf'' drops the $\Gamma$-compact object support condition and the 
twiddle adds an additional $\cale_{d_{prod}}$ condition, i.e.\ each morphism is 
bounded in the product metric. We use the obvious analogous notation
\[
\cald_{lf}^{\Gamma}(\overline{X},\overline{A};p), \quad \tilde{\cald}^{\Gamma}(\overline{X},\overline{A};p)
\quad  \mbox{and} \quad \tilde{\cald}_{lf}^{\Gamma}(\overline{X},\overline{A};p)
\]
for pairs. Compare~\ref{pairnotation}.
\end{definition}

\begin{remark} \label{pairseq}
It can be deduced from  Proposition~\ref{standardkaroubi} that there are ``twiddle'' and ``lf''
analogues of the pair fibration sequences, e.g.\ there is a fibration sequence
\[
\IK^{-\infty} \tilde{\cald}_{lf}^\Gamma(\overline{A};p|_{\overline{A}}) \to 
\IK^{-\infty} \tilde{\cald}_{lf}^\Gamma(\overline{X};p) \to 
\IK^{-\infty} \tilde{\cald}_{lf}^\Gamma(\overline{X},\overline{A};p).
\]
\end{remark}

\begin{remark} \label{resremark}
One could now hope for a restriction functor  
$res : \tilde{\cald}_{lf}^\Gamma(\overline{X};p) \to \tilde{\cald}_{lf}^H(\overline{X};p)$
for arbitrary subgroups $H \subset  \Gamma$. But recall that our definition of 
equivariant continuous control (compare \ref{ccdef} \ref{ccdefcc}) used isotropy invariant neighborhoods
in $X$, and of course the isotropy groups may change under restriction. 
However, a sufficient condition for such a restriction is that $X$ is a free $\Gamma$-space or  
discrete.
\end{remark}

\begin{remark} \label{lfsepremark}
There is a version of Lemma~\ref{separate} for morphisms in
$\tilde{\cald}^\Gamma_{lf}(\overline{X};p)$ under certain compactness
conditions. In particular, consider a disjoint union 
$\overline{X} = \overline{X}_0 \cup \overline{X}_1$, where the $p(\overline{X}_i)$ consist
of different components of ${X}$ and $\overline{X}_0$ is $\Gamma$-compact.
Then 
\[
\tilde{\cald}_{lf}^\Gamma(\overline{X};p) = 
\tilde{\cald}_{lf}^\Gamma(\overline{X}_0;p) \oplus 
\tilde{\cald}_{lf}^\Gamma(\overline{X}_1;p).
\]
\end{remark}

Note that for $\Gamma$-compact spaces $\overline{X}$ we have of course
$\cald^{\Gamma}_{lf}(\overline{X};p)=\cald^{\Gamma}(\overline{X};p)$. The following 
Lemma tells us that in some cases it is  irrelevant whether
or not we add the product metric condition.

\begin{lemma} \label{tildedrop}
Let $X$ be $\Gamma$-space and let $p:X \to \pi_0(X)$ be the quotient map
which assigns to a point its path component.
\begin{enumerate}
\item
We have $\tilde{\cald}^{\Gamma}(X;\id) = \cald^{\Gamma} ( X; \id)$.
\item
The inclusions
$\tilde{\cald}^{\Gamma}( X; \ast) \to  \cald^{\Gamma}(X;\ast)$   and  
$\tilde{\cald}^{\Gamma}( X; p)  \to  \cald^{\Gamma}(X;p)$ 
induce isomorphisms on the level of $K$-theory.
\end{enumerate}
\end{lemma}

\begin{proof}
(i) This is a consequence of the following (easily checked) fact: 
Every morphism $\phi$ which fulfills the $\cale_{\Gamma cc} (X)$ condition
is already bounded with respect to the product metric. 

(ii) The map of coarse pairs 
\begin{eqnarray*}
(X \times \einsu , X \times \{ 1 \} ,
(\ast \times \id)^{-1} \cale_{\Gamma cc} ( \punkt ) \cap \cale_{d_{prod}} ) \\
\to 
(X \times \einsu , X \times \{ 1 \} ,
(\ast \times \id)^{-1} \cale_{\Gamma cc} ( \punkt ) )
\end{eqnarray*}
leads to a comparison map between the corresponding fibration sequences in $K$-theory, 
compare \ref{standardkaroubi}.
The left hand comparison map is an equality already on the level of categories.
The middle terms of the fibration sequences are both contractible by the usual Eilenberg
swindle towards infinity, cf.\ Example~\ref{swindletoinfty}.
Hence the right hand map induces an  equivalence
as claimed. 

To treat the second map one first uses Remark~\ref{disremark} to observe that 
$(p \times \id_{\einsu})^{-1} \cale_{\Gamma cc}$ can be replaced by 
$(p\times \id_{\einsu} )^{-1} \cale_{dis}$. Now 
$(X \times \einsu ,( p\times \id_{\einsu} )^{-1} \cale_{dis} )$ admits an Eilenberg swindle
and one can argue as above.
\end{proof}


\typeout{---------------------geometryI---------------------}


\section{Geometric preparations needed for proving  injectivity}
\label{sec_geometryI}

Let $M$ be an $n+1$-dimensional compact Riemannian manifold with 
strictly negative sectional curvature. Let $\Gamma = \pi_1(M)$.
Let $\pi:\tilde{M} \to M$ denote the universal cover and 
let $d$ denote the distance on $\tilde{M}$.
In general there is no $\Gamma$-equivariant contracting self map on $\tilde{M}$.
In this section we will discuss contracting maps that are invariant under an infinite cyclic subgroup.
This will be important in Section~\ref{sec_injectivity}. 

Let $f : \IR \to \tilde{M}$ be a geodesic such that the setwise stabilizer of the image
$C = \{\gamma \in \Gamma|\gamma f(\IR) = f(\IR)\}$ is nontrivial (and hence infinite cyclic) 
and let $g =f(\IR)$. Let $Exp: T\tilde{M} \to \tilde{M} $ be the exponential map; i.e.\ 
if $v \in T_x \tilde{M}$, then $Exp(v)=\alpha_v(1)$, where $\alpha_v$
is the unique geodesic in $\tilde{M}$ such that $\dot{\alpha}_v(0)=v$.
Let $E$ be the subset of $T\tilde{M}$ consisting of all vectors $v$
whose foot $\alpha_v(0)\in g$ and such that $v \perp g$ at $\alpha_v(0)$.
Note that $E \to g$, $v \to \alpha_v(0)$ is a $C$-equivariant $\IR^n$-vector
bundle with structure group $O(n)$, where $n+1=\dim \tilde{M}$. 
A choice of a trivialisation allows us to identify $E$ with $g \x \IR^{n}$.
Let 
\[
\psi : g \x \IR^n \to \tilde{M}
\]
be the restriction of $Exp$ under this identification.
Note that $\psi$ is a $C$-equivariant diffeomorphism where $C$ acts by translation 
on $g$ and orthogonally on $\IR^{n}$. 
We want to study the following $C$-equivariant map:
\begin{eqnarray} \label{defineq}
q: \tilde{M} \times \einsu & \to & \tilde{M} \times \einsu  \\
(\psi(x,v) , t ) & \mapsto & (\psi(x, v/t ), t)   \nonumber
\end{eqnarray}
Our aim is to prove that this map is contracting and therefore allows to regain
some control, cf.\ Lemma~\ref{q-map}.

\begin{proposition}
\label{coolingtower}
There is a constant $K$
such that
\[
d(\psi(x,v/r'),\psi(y,w/r'))
\leq
\frac{2}{r} \cdot d(\psi(x,v),\psi(y,w))
\]
provided 
$|v|,|w| \geq  r \cdot
 (K + d(\psi(x,v),\psi(y,w)))$ and $r' \geq r \geq 1$. 
Here $|\cdot|$ denotes the Euclidean norm on $\IR^n$.
\end{proposition}

The proof depends on the following lemmata.

\begin{lemma} \label{constK}
Let $p:\tilde{M} \to g$ be the orthogonal projection.
There is a constant $K$ such that the following holds:
if $X \in T_{\psi(x,V)}(\tilde{M})$ with $|V| \geq K$ then
$\|dp(X)\| \leq \|X\| / 2$. Here $\| . \|$ denotes the norm induced
by the Riemannian metric on the tangent space.
\end{lemma}

\begin{proof}
This is \cite[Lemma 1.2]{Farrell-Jones(1987b)}.
\end{proof}

\begin{lemma} \label{claims}
Let $u,u' \in E$ and $\phi(t) = d(\alpha_u(t),\alpha_{u'}(t))$. Then 
\begin{enumerate}
\item \label{claim2}
$\phi(s) \leq \phi(t)$ for $0\leq s \leq t$,
\item \label{claim1} 
$\phi(t) \geq (\phi(1)-\phi(0)) t + \phi(0)$ for $t \geq 1$. 
\end{enumerate}
\end{lemma} 

\begin{proof}
It is well known that $\phi$ is a convex function, 
cf.~\cite[p.4]{Ballmann-Gromov-Schroeder(1985)}.
This implies our claims by elementary arguments.
\end{proof}

\begin{proof}[Proof of \ref{coolingtower}.]
By \ref{claims} \ref{claim2} we may assume $r=r'$.
Set $V=\frac{1}{r}v$ and $W=\frac{1}{r} w$, then
\ref{claims} \ref{claim2} shows that 
\[
d( \psi( x , V ), \psi( y , W) ) 
 \leq d( \psi (x,v) , \psi( x, w))
\]
and hence that 
\[
|V|,\mbox{ }|W| 
 \geq K + d( \psi(x, V) , \psi(y, W) ).
\]
Consequently the geodesic segment connecting $\psi(x, V)$ to $\psi(y , W)$
never gets closer than $K$ to $g$. Therefore
\[
d(x,y) \leq \frac{1}{2} d( \psi(x, V), \psi(y, W) )
\]
by \ref{constK}. 
By applying \ref{claims} \ref{claim1} with $u=(x,V)$, $u'=(y,W)$ and $t=r$
together with this  inequality and observing that $\alpha_u(r) = \psi(x , v )$
and $\alpha_{u'}(r)=\psi(y, w)$ we see that 
\[
d( \psi( x, v ) , \psi (y, w )) \geq (\phi(1)-\phi(0)) r 
 \geq \frac{1}{2} d(\psi(x,V),\psi(y,W))r.
\]
Multiplying this last inequality by $\frac{2}{r}$ yields the one posited
in the Proposition.
\end{proof}

\begin{remark}  \label{vierecksungleichung}
For $(x,v)$ and $(y,w)$ in $g \times \IR^n$ we have 
\[
|w| \leq d ( \psi(x,v),\psi(y,w) ) + |v |
\]
because $|v|=d ( \psi(x,v) , g)$, $|w|=d ( \psi(y,w) , g)$
and 
\[
d ( \psi(y,w) , g )\leq 
d(\psi(y,w) , \psi(x,v)) + d( \psi(x,v),g).
\]
\end{remark}

We can now give an exact statement about the 
control gained by the map $q$, cf.~(\ref{defineq}). 
Proposition~\ref{coolingtower} will enter the proof of injectivity through the 
following lemma.
Recall the ``lf'' and $\tilde{\cald}$-notation from Definition~\ref{deflf}.
Let $p_{\tilde{M}/g}: \tilde{M} \to \tilde{M}/g$ be the quotient map
which collapses $g$ to a point.

\begin{lemma} \label{q-map}
The map $q$ is $C$-invariant and induces functors
\begin{eqnarray*}
\tilde{\cald}^{\Gamma}( \tilde{M} ; \id ) & \to & \tilde{\cald}^{C}_{lf} (\tilde{M} ; \id ) \\
\tilde{\cald}^{\Gamma}( \tilde{M} ; \ast) & \to & \tilde{\cald}^{C}_{lf} (\tilde{M} ; p_{\tilde{M} / g}) .
\end{eqnarray*}
\end{lemma}
\begin{proof}
We have to prove the following: 
\begin{enumerate}
\item
If $E \in \cale_{\Gamma cc}( \tilde{M} ) \cap \cale_{d_{prod}}(\tilde{M} \times \einsu )$, then 
$(q \times q)(E) \in \cale_{C cc}( \tilde{M} )$.
\item
If $E \in \cale_{d_{prod}}(\tilde{M} \times \einsu )$, then 
$(q \times q)(E) \in p_{\tilde{M}/g}^{-1}\cale_{Ccc}(\tilde{M}/g)$.
\end{enumerate}
(i) Note that since we have a free $\Gamma$-action all neighborhoods in the 
definition of equivariant continuous control are ordinary neighborhoods.
Let $U$ be an open neighborhood of $\psi(x,v) \in \tilde{M}$ and let $r \geq 1$. Suppose there
exists no $V \subset U \times (r, \infty)$ as in the continuous control
condition~\ref{ccdef}~\ref{ccdefcc}. Then there exists a sequence
$((\psi(x_n, v_n),t_n),(\psi(y_n,w_n),s_n)) \in E$ such that 
\begin{eqnarray}
(\psi(x_n,v_n/t_n),t_n) &  \in & B_{1/n}(\psi(x,v)) \times ( n,\infty )  \label{stern1} \\
\mbox{ and } 
(\psi(y_n,w_n/s_n),s_n) &  \notin & U \times (r,\infty).  \label{stern2}
\end{eqnarray}

Since $\psi$ is a homeomorphism (\ref{stern1}) implies that $x_n$ tends to $x$,
$v_n/t_n$ tends to $v$ and $t_n$ tends to $\infty$. Using  the $\cale_{d_{prod}}$-condition
and \ref{ccdef}~\ref{ccdefmetric}
there exists an $\alpha >0$ such that 
\begin{eqnarray}
|s_n - t_n| \mbox{ , }  d( \psi(x_n,v_n) , \psi(y_n,w_n) ) < \alpha . \label{stern3}
\end{eqnarray}
Hence $s_n$ tends to $\infty$ and $v_n/s_n$ tends to $v$.
Let $K$ be the constant from Proposition~\ref{coolingtower} and set
$\tilde{r}_n = \min \{ | v_n | , | w_n | \} \cdot (K + \alpha)^{-1}$.
By considering a suitable subsequence it will be sufficient to consider the 
case where $\tilde{r}_n$ tends to infinity (Case I) or the case where
it is bounded (Case II).

{\bf Case I}: Let $r_n= \min \{ s_n , \tilde{r}_n \}$, then $r_n$ tends to
infinity. Eventually $s_n \geq r_n \geq 1$ and 
\[
\min \{ |v_n | , | w_n | \} = \tilde{r}_n ( K + \alpha ) 
\geq r_n ( K + d ( \psi(x_n , v_n) , \psi(y_n , w_n ) ).
\]
Applying Proposition~\ref{coolingtower} gives that 
\[
d ( \psi(x_n , v_n/s_n ) , \psi(y_n , w_n/s_n ) ) < \frac{2}{r_n} \alpha
\]
which tends to zero. Hence $\lim \psi( y_n , w_n/s_n) = \lim \psi( x_n , v_n/s_n ) = \psi( x,v) $
which contradicts (\ref{stern2}).

{\bf Case II}:
Let $\tilde{r}_n$ be bounded. Then also $\min \{ | v_n |, |w_n | \}$  is bounded.
Now (\ref{stern3}) and \ref{vierecksungleichung} imply that $|v_n | \leq \alpha + |w_n |$.
Passing to a subsequence we may assume therefore that
$v_n$ tends to $\tilde{v}$  and hence $\psi(x_n , v_n) \to \psi(x, \tilde{v} )$ for some 
$\tilde{v}$. Using (\ref{stern3}) this implies 
$d ( \psi(x ,\tilde{v} ) , \psi ( y_n , w_n ) ) < \alpha + 1$ for large $n$ and by a compactness 
argument (passing to another subsequence) we have $\psi(y_n , w_n ) \to \psi( y, \tilde{w} )$
for some $\tilde{w}$. Now since $E \in \cale_{\Gamma cc}$ Lemma~\ref{convergence} implies that 
$\psi(x, \tilde{v} ) = \psi( y, \tilde{w} )$.
Hence both $\psi(y_n , w_n/s_n )$ and $\psi ( x_n , v_n/t_n)$ tend to 
$\psi( x,0)$. This contradicts (\ref{stern1}) or (\ref{stern2}).

(ii) Let $z=p_{\tilde{M}/g}(g)$ be the special point in $\tilde{M}/g$.
Let $p_{\tilde{M}/g}(\psi(x,v)) \in \tilde{M}/g$ be given. Consider first the case where $\psi(x,v) \notin g$.
While trying to find a $V \subset U \times (r , \infty)$ satisfying the 
equivariant continuous control condition we can always replace $U$ by the smaller open
neighborhood $U - z$. From there on we can proceed exactly as in (i).
Note that $E \in \cale_{\Gamma cc}$  was only used in Case II, but Case II implied that $|v_n|$
is bounded and therefore $v=0$.
In Case I we only used (\ref{stern3}) which follows from the $\cale_{d_{prod}}$ condition.

Now consider the case $v=0$, i.e.\ $\psi(x,v) \in g$. Let $U^{\prime}$ be a $C$-invariant
neighbourhood of $z$. Let $U = p_{\tilde{M}/g}^{-1}(U^{\prime})$. Since 
$C$-acts orthogonally on the $\IR^n$-factor we can assume that $U=\psi(g \times B_{\epsilon}(0))$
for some $\epsilon >0$. Let $r >0$. We proceed as usual by contradiction:
Assume that there exists a sequence $((\psi(x_n,v_n),t_n),(\psi(y_n,w_n),s_n)) \in E$ such that
\begin{eqnarray}
(\psi(x_n , v_n/t_n),t_n ) & \in & \psi(g \times B_{1/n}(0)) \times (n, \infty ) \label{stern5} \\
\mbox{ and } (\psi(y_n , w_n/s_n ) , s_n ) & \notin & 
\psi( g \times B_{\epsilon}(0) ) \times (r , \infty)  \label{stern6}
\end{eqnarray}

Let $\alpha > 0$ satisfy (\ref{stern3}) and. As above we get
$t_n$, $s_n \to \infty$ and $v_n/s_n \to 0$. If $|w_n|$ is bounded we have
$(y_n , w_n/s_n ) \in g \times B_{\epsilon}(0) $ for large $n$ contradicting (\ref{stern6}). 
Therefore and by \ref{vierecksungleichung}
we may assume that $|w_n|$ and  $|v_n|$ tend to infinity.
If we define $r_n$ as in (i) we see that it also tends to infinity. 
Now applying Proposition~\ref{coolingtower} we get
\begin{eqnarray*}
| w_n / s_n | &  =  & d ( g , \psi(y_n, w_n/s_n ) ) \\
& \leq &
d ( g , \psi( x_n , v_n/s_n ) )  
+ d (\psi(x_n , v_n/s_n ) , \psi( y_n , w_n/s_n )) \\
& \leq & 
| v_n /s_n | +  \frac{2}{r_n} \alpha.
\end{eqnarray*}
The right hand side tends to zero. This contradicts (\ref{stern6}).
\end{proof}


\typeout{--------------------- injectivity ---------------------}


\section{Injectivity}  
\label{sec_injectivity}
In this section we will prove the following theorem.

\begin{theorem} \label{injectivity}
Let $\Gamma$ be the fundamental group of a compact Riemannian manifold with
strictly negative sectional curvature.
Then for every $n \in \IZ$ and every associative ring $R$ the assembly map
\[
\xymatrix{
H_n^{\Or \Gamma}( E \Gamma ( \Cyc ) ; \IK R^{-\infty} ) 
\ar[rr]^-{A_{\Cyc \to \All}} & &  K_n( R \Gamma ) 
         }
\]
is injective.
\end{theorem}

We will first use the negative curvature assumption to construct a specific model $p:E \to E(\Cyc)$
for the universal map  $E\Gamma \to E \Gamma ( \Cyc )$. 
This will allow us to interpret
the above assembly map as the map induced by the map of 
resolutions (compare page \pageref{mapofresolutions})
\[
\xymatrix{
E \ar[d]_p \ar[r]^{\id} & E \ar[d]^p  \\
E(\Cyc) \ar[r] & \punkt .
         }
\]
We then prove that the left hand side of the assembly map 
is a direct sum of $H_n(B\Gamma;\IK R^{-\infty})$
and Nil-terms (see Corollary~\ref{lhscomputation}).
Since the classical assembly map $A_{1 \to \All}$ is known to be injective in our case 
by \cite{Carlsson-Pedersen(1995)} it remains to detect the Nil-terms in the image.
Using again the negative curvature assumption 
(in particular Proposition~\ref{coolingtower})
we construct maps out of $K_n(R\Gamma)$ which in fact do detect the Nil-terms.

\subsection{A model for the universal map $E \Gamma \to E \Gamma ( \Cyc )$}
\label{modelforuniversalmap}

Let $M$ be our compact Riemannian manifold of strictly negative sectional curvature.
Let $\pi: \tilde{M} \to M$ be its universal covering and let $\Gamma= \pi_1 (M)$ be its 
fundamental group.

\begin{lemma} 
\label{fixonlyone}
Let $f,f':\IR \to \tilde{M}$ be two different geodesics, i.e.\ 
$f(\IR) \neq f'(\IR)$. Suppose that an isometry $\gamma$ 
of $\tilde{M}$ stabilizes (setwise) $f(\IR)$ and $f'(\IR)$. Then 
$\gamma$ does not act freely  on $\tilde{M}$.
\end{lemma}

\begin{proof}
Since 
$\tilde{M}$ has negative curvature, there are unique
$t,t'$ such that 
\[
\mathit{dist}_{\tilde{M}}(f(t),f'(t'))=
 \mathit{dist}_{\tilde{M}}(f(\IR),f'(\IR)),
\]
cf.~\cite[p.8]{Ballmann-Gromov-Schroeder(1985)}.
Uniqueness implies that $f(t)$ and $f'(t')$ are fix under $\gamma$. 
\end{proof}

Let $\Cyc$ denote the family of cyclic subgroups of $\Gamma$.
Choose a set of representatives $\{ C_i | i \in I \} \subset \Cyc$ for the 
conjugacy classes of maximal infinite cyclic subgroups, i.e\ every maximal infinite
cyclic subgroup is conjugate to a unique $C_i$. Let $E= \tilde{M} \times \IR^2$ and 
let $p_{\tilde{M}}: E \to \tilde{M}$ be the standard projection.

\begin{lemma}
\label{geodesics}
There are embeddings $f_i:\IR \to  E$ satisfying the following:
\begin{enumerate}
\item
Each $p_{\tilde{M}} \circ f_i$ is a geodesic whose image in $M$ is closed.
\item
The setwise stabilizer of $f_i(\IR)$ is $C_i$,\ i.e.\ $C_i = 
 \{\gamma \in \Gamma|\gamma f_i(\IR)=f_i(\IR)\}$.
\item \label{unique}
Each $\gamma_0 \neq 1 \in \Gamma$ stabilizes (setwise) a unique set
of the form $\gamma f_i(\IR)$ (with $i \in I,\gamma \in \Gamma$).
\item
We have $f_i(\IR) \cap \gamma f_j(\IR) = \emptyset$ unless
$\gamma \in C_i$ and $i=j$.
There is an isomorphism of $\Gamma$-spaces
\[
\coprod_{i \in I} \Gamma \times_{C_i} f_i( \IR ) \to X= \bigcup_{i \in I} \Gamma f_i( \IR  ) \subset E.
\]
Moreover $X = \bigcup_{i \in I } \Gamma f_i(\IR)$ is an embedded
submanifold of $E$ and therefore the inclusion is a cofibration.
\end{enumerate}
\end{lemma} 

\begin{proof}
Let $\gamma_i$ be a generator of $C_i$. Then
$\gamma_i$ can be represented (up to free homotopy) 
by a closed geodesic ${f'_i} : S^1 \to {M}$ (\cite[I.3.16]{Bridson-Haefliger(1999)}). 
Choose embeddings 
$h_i : S^1 \to \IR^{2}$, such that
$\mathit{dist}(h_i(S^1), h_j(S^1)) > 1$ for
$i \neq j$.  
Then the embeddings
$f'_i \times h_i : S^1 \to M \times \IR^{2}$ 
can be lifted to the required
embeddings $f_i : \IR \to E$.
Uniqueness in (iii) is a consequence of Lemma~\ref{fixonlyone}.
\end{proof}

One should think of the $f_i$ as perturbed copies of lifted closed geodesics
avoiding intersections and self intersections. Let now $E(\Cyc)$ be the 
quotient of $E$ where every path component of 
$X= \bigcup_{i \in I } \Gamma \gamma f_i( \IR )$ is collapsed to a point
and let $p:E \to E( \Cyc )$ be the natural map, i.e.\ we have by definition
the following push-out diagram:
\[
\xymatrix{
X \ar[d] \ar[r]^-{\mathit{inc}} & E \ar[d]_p \\
\pi_{0}(X) \ar[r] & E(\Cyc ).
         }
\]

\begin{proposition}
The map $p:E \to E(\Cyc)$ is a model for the universal map 
$E \Gamma \to E \Gamma (\Cyc)$.
\end{proposition}
\begin{proof} Using Lemma~\ref{geodesics} (iii) and (iv) 
one checks that $E(\Cyc)^H$ is
contractible if $H \subset \Gamma$ is infinite 
cyclic or trivial and empty otherwise.  
\end{proof}


\subsection{Splitting off the Nil-terms}

Let $X_i= \Gamma f_i(\IR)$ be the $\Gamma$-orbit of $f_i(\IR)$ in $E$. Consider the following
diagram
\begin{eqnarray} \label{seconddiagram}
\xymatrix{
f_i( \IR ) \ar[r] \ar[d]^{p|_{g_i}} & 
\Gamma \times_{C_i} f_i( \IR ) = X_i  \ar[r] \ar[d]^{p|_{X_i}} & 
X = \coprod X_i \ar[r] \ar[d]_{p|_{X}} \ar@{}[dr]|{\mbox{{\tiny  h-pushout}}} & 
E \ar[d]^{p} \\
\punkt \ar[r] & 
\Gamma \times_{C_i} \punkt = \pi_0( X_i )  \ar[r] & 
\pi_0( X ) \ar[r] & E(\Cyc).
         }
\end{eqnarray}

Replacing the left 
hand arrow in (\ref{seconddiagram}) by a map of resolutions one obtains a functor
\[
\cald^{C_i}(f_i(\IR) ; \id) \to 
\cald^{C_i}( f_i( \IR ) ; p|_{f_i( \IR )} ).
\]
Since $f_i(\IR )$ is a copy of the real line with the standard action
of the infinite cyclic group  this functor induces a model for the classical
assembly map $A_{1 \to \All}$ of the infinite cyclic group on the level of $K$-theory.
Therefore by the Bass-Heller-Swan decomposition
we have a split cofibration sequence of spectra 
\[
\IK^{-\infty} \cald^{C_i}(f_i(\IR) ; \id) \to \IK^{-\infty} \cald^{C_i}( f_i( \IR ) ; p|_{f_i( \IR )} ) \to
\IN_i.
\]
Here $\IN_i$ is by definition the homotopy-cofiber of the assembly map whose homotopy groups
are known to be
\[
\pi_n( \IN_i)   = NK_n( R ) \oplus NK_n( R ).
\]
These Nil-groups were defined before \ref{nilcorollary}.
Resolving the left hand square in (\ref{seconddiagram}) and using the induction isomorphisms
from \ref{inddni}
we see that the above leads to a split homotopy-cofibration sequence
\[
\IK^{-\infty} \cald^{\Gamma} ( X_i ; \id )  \to
\IK^{-\infty} \cald^{\Gamma} ( X_i ; p|_{X_i} )\to \IN_i.
\]
The disjoint union axiom (compare page~\pageref{disjointunion}) 
gives us the split homotopy-cofibration sequence
\[
\IK^{-\infty} \cald^{\Gamma} ( X ; \id )  \to
\IK^{-\infty} \cald^{\Gamma} ( X ; p|_{X} ) \to \wed_{i \in I} \IN_i.
\]
The right hand square in (\ref{seconddiagram})
is a homotopy push-out square since $X \to E$ is an embedded submanifold and therefore 
a $\Gamma$-cofibration. Our homology theory turns this into a homotopy-pullback
or equivalently a homotopy-pushout square and we hence obtain the following
split homotopy-cofibration sequence:
\[
\IK^{-\infty} \cald^{\Gamma} ( E; \id )  \to
\IK^{-\infty} \cald^{\Gamma} ( E; p ) \to \wed_{i \in I} \IN_i.
\]
In particular we have the following consequence:
\begin{proposition}  \label{lhscomputation}
There exists a splitting
\[
H_{n}^{\Or \Gamma}( E \Gamma ( \Cyc ); \IK R^{-\infty} ) \cong  
 H_{n}(B\Gamma;\IK R^{-\infty}) \oplus \bigoplus_i NK_i (R) \oplus NK_i (R) .
\] 
\end{proposition}


\subsection{Detecting the Nil-terms}

We will now finish the proof of Theorem~\ref{injectivity}.
Below we will construct for each $n \in \IZ$ the  following commutative diagram:
\begin{eqnarray} \label{maindiagram}
\xymatrix{
K_n \cald^{\Gamma} (X_i; \id ) \ar[d] \ar[r] & 
K_n \cald^{\Gamma} (X_i ; p|_{X_i} ) \ar[d] \ar[r] & 
N_i  \ar[d] \\
K_n \cald^{\Gamma} (E; \id) \ar[r]^{A_{1 \to \Cyc}} \ar[dd]  \ar[dr]_{A_{1 \to \All}}  &
 K_n \cald^{\Gamma} (E;p)  \ar[d]^{A_{\Cyc \to \All}}  \ar[r] &   
\bigoplus_{k \in I}  N_k  \ar[dd] \\
&  K_n \cald^{\Gamma} (E; \ast )  \ar[d] &   \\
K_n \tilde{\cald}^{C_j}_{lf} ( \tilde{M} ; \id ) \ar[r] & 
K_n \tilde{\cald}^{C_j}_{lf} ( \tilde{M} ;p_{\tilde{M}/{g_j} }) \ar[r] & 
\tilde{N}_j
           }
\end{eqnarray}
Here the $n$'th homotopy group of $\IN_i$ is denoted $N_i$ and $\tilde{N_j}$ is by definition the cokernel
of the lower left hand arrow.
For the definition of $\tilde{\cald}^{\Gamma}_{lf}(\overline{X};p)$ compare~\ref{deflf}.

We have the following facts about this diagram:
\begin{enumerate}
\item
By the results from the previous Subsection the first two rows
are short split exact sequences. 
\item
We will show below in Lemma~\ref{iji} that the right  hand vertical map $N_i \to \tilde{N}_j$ 
is the zero map if $i\neq j$ and
it is injective if $i=j$. 
\item
The map $A_{1 \to \All}$ is injective since this is 
a model for the ``classical'' assembly map for $\Gamma$ and injectivity is proven
in \cite{Carlsson-Pedersen(1995)}.
\end{enumerate}
These facts imply the injectivity of $A_{\Cyc \to \All}$ and hence Theorem~\ref{injectivity}.

Let us now construct the diagram.
The non-obvious maps in the diagram are the vertical arrows in the lower left hand
square. To construct them we will need the geometry discussed in the 
previous section. Let $p_{\tilde{M}/g_i}: \tilde{M} \to \tilde{M}/g_i$ 
be the quotient map obtained by collapsing the geodesic 
$g_i=p_{\tilde{M}}(f_i(\IR ))$ to a point. Using the identification 
$\tilde{M} \cong g_i \times \IR^N$ provided 
by \ref{coolingtower} we define a map $q_i$ as follows
\begin{eqnarray*}
q_i: \tilde{M} \times \einsu = g_i \times \IR^N \times \einsu  & \to &
\tilde{M} \times \einsu = g_i \times \IR^N \times \einsu \\
(x,v,r) & \mapsto & (x,v/{r},r) 
\end{eqnarray*}
Recall from Definition~\ref{deflf} that a twiddle adds a product-metric condition on 
morphisms and ``lf'' drops the $\Gamma$-compact object support condition.
Let $p_{\tilde{M}}:E=\tilde{M} \times \IR^2 \to \tilde{M}$ be the standard projection.
By Lemma~\ref{q-map} the map $q_j$ induces the functors with the same name in the following 
commutative diagram.
\[
\xymatrix{
\tilde{\cald}^{\Gamma} ( X_i;\id ) \ar[r]^-{\inc} \ar[d]_{F_{N}} &
\tilde{\cald}^{\Gamma} (E ; \id) \ar[r]^{p_{\tilde{M}}} \ar[d] &
\tilde{\cald}^{\Gamma} (\tilde{M} ;\id) \ar[r]^-{q_j}  \ar[d] &
\tilde{\cald}^{C_j}_{lf} (\tilde{M} ;\id) \ar[d]_{F_{\tilde{N}}} \\
\tilde{\cald}^{\Gamma} (X_i ; p|_{X_i}) \ar[r]^-{\inc} &
\tilde{\cald}^{\Gamma} (E ; \ast) \ar[r]^{p_{\tilde{M}}}  &
\tilde{\cald}^{\Gamma} (\tilde{M} ; \ast) \ar[r]^{q_j} & 
\tilde{\cald}^{C_j}_{lf} (\tilde{M}; p_{\tilde{M}/ g_j} )
         }
\]
Here all vertical functors are induced by identities. By Lemma~\ref{tildedrop}
we can drop some of the twiddles and we obtain the lower left hand side of (\ref{maindiagram}).
(Note that vertical arrows here are horizontal in (\ref{maindiagram}) and vice versa.)

By definition $N_i$ and $\tilde{N}_j$ are the cokernels of the maps induced by the functors $F_{N}$ and 
$F_{\tilde{N}}$:
\begin{eqnarray} 
 N_i &  = & \coker(K_n(F_{N}))  \label{ai} \\ 
 \tilde{N}_j &  = & \coker(K_n(F_{\tilde{N}})). \label{bi}
\end{eqnarray}
To study the induced map between these cokernels we may replace the above diagram with the following. 
(For the existence of the arrows labeled $\res$ see Remark~\ref{resremark}.)
\[
\xymatrix{
\tilde{\cald}^{\Gamma} ( X_i;\id ) \ar[r]^-{\res} \ar[d]_{F_{N}} &
\tilde{\cald}_{lf}^{C_j} (X_i ; \id) \ar[r]^-{q_j|_{X_i}} \ar[d]^F &
\tilde{\cald}^{C_j}_{lf} (\tilde{M} ;\id) \ar[d]_{F_{\tilde{N}}} \\
\tilde{\cald}^{\Gamma} (X_i ; p|_{X_i}) \ar[r]^-{\res} &
\tilde{\cald}_{lf}^{C_j} (X_i; p|_{X_i} ) \ar[r]^-{q_j|_{X_i}} &
\tilde{\cald}^{C_j}_{lf} (\tilde{M}; p_{\tilde{M}/ g_j} )
         }
\]
Studying the functor $F$ we obtain:

\begin{lemma} \label{iji}
$\mbox{   }$
\begin{enumerate}
\item \label{ij}
For $i \neq j$ the map $N_i \to \tilde{N}_j$ is the zero map.
\item \label{ii} 
If $i=j$ then $N_i \to \tilde{N}_i$ is injective. 
\end{enumerate}
\end{lemma}

\begin{proof}
(i) 
If $i \neq j$ we have an isomorphism of $C_j$-spaces $X_i \cong \Gamma/C_i \times \IR$, 
where $C_j$ acts trivially on $\IR$. Under this identification $p|_{X_i}$ is the projection 
onto the first factor. Define the $C_j$-equivariant contraction
\begin{eqnarray*}
X_i \times [1,\infty) = \Gamma/C_i \times \IR \times [1,\infty)
 & \to & X_i \times [1,\infty) \\
(\gamma C_i,t,r) & \mapsto & (\gamma C_i,t/r,r).
\end{eqnarray*}
This map induces a functor 
$G:\tilde{\cald}_{lf}^{C_i} (X_i ; p|_{X_i} ) \to \tilde{\cald}_{lf}^{C_i} ( X_i; \id )$ which splits 
$F$, i.e.\ $F \circ G \simeq ID$ and in particular the cokernel of the map induced by $F$ vanishes.

(ii) 
Let $Y_i := X_i - f_i(\IR)$. By Remark~\ref{lfsepremark} there is a decomposition of the functor 
$F$ as 
\[
\xymatrix{
\tilde{\cald}^{C_i} ( f_i(\IR) ;\id ) 
\oplus
\tilde{\cald}_{lf}^{C_i} (Y_i ; \id) 
\ar[rr]^-{\left( 
\begin{array}{cc} F_1 & 0\\ 0 & F_2 \end{array}
           \right) }  & &
\tilde{\cald}^{C_i} ( f_i(\IR) ; p|_{f_i (\IR)} ) \oplus
\tilde{\cald}_{lf}^{C_i} ( Y_i ; p|_{Y_i} ) .
        }
\]
Using the induction isomorphisms from \ref{inddni}  we see that 
the map $F_{N}$ is a retract of $F$. It will therefore be sufficient to prove
that that the cokernel of the map $K_n(F_2)$ is trivial
and that the map from $\coker(K_n(F_1))$ to $\tilde{N}_i=\coker(K_n(F_{\tilde{N}}))$ is injective.

Note that $Y_i \cong  (\Gamma - C_i) \times_{C_i} \IR$. Since the  left action
of $C_i$ on $(\Gamma - C_i)/C_i$ is free $Y_i$ is non canonically
isomorphic as a $C_i$-space to $(\Gamma - C_i) / C_i \times \IR$ where
$\IR$ is equipped with the trivial $C_i$-action, see Proposition~\ref{geodesics} \ref{unique}.
Hence similar as in the previous proof we can split $F_2$.

Now note that $f_i(\IR)\cong g_i$ and consider the diagram
\begin{eqnarray} \label{nocheindiagram}
\xymatrix
{
 \tilde{\cald}^{C_i}(g_i;\id) \ar[r] \ar[d]^{F_1} &
 \tilde{\cald}^{C_i}_{lf}(\tilde{M};\id)  \ar[r] \ar[d]^{F_{\tilde{N}}} &
 \tilde{\cald}^{C_i}_{lf}(\tilde{M}, g_i ;\id) 
               \ar[d]^{F_3}  
 \\
 \tilde{\cald}^{C_i}(g_i ;*) \ar[r] &
 \tilde{\cald}^{C_i}_{lf}(\tilde{M};p_{\tilde{M}/g_i}) \ar[r] &
 \tilde{\cald}^{C_i}_{lf}(\tilde{M}, g_i ; p_{\tilde{M}/g_i} ) 
}
\end{eqnarray}
We claim that $F_3$ induces an isomorphism in $K$-theory. 
Using the fibrations of spectra (see \ref{pairseq}) induced by the diagram this will imply that the 
map $\coker(K_n(F_1)) \to \coker(K_n(F_{\tilde{N}}))$ is injective. 
To verify the claim one replaces $g_i$ with a $C_i$-invariant tubular neighborhood $T \supset g_i$
in (\ref{nocheindiagram}). By homotopy invariance (\ref{homotopyinvarianceCor}) 
\[
\tilde{\cald}_{lf}^{C_i} ( g_i ; \id) =  \cald^{C_i} (g_i ; \id ) \to
\cald^{C_i} ( T ; \id ) =  \tilde{\cald}_{lf}^{C_i} ( T ; \id )  
\]
induces an  isomorphism in $K$-theory. (The equalities use Lemma~\ref{tildedrop}.)
Similarly  we get an isomorphism with $p=p_{\tilde{M}/g_i}$ instead of the identities. 
Using the various fibrations (see \ref{pairseq}) we conclude  that it is sufficient to prove that 
$\cald_{lf}^{C_i} ( \tilde{M} , T ; \id) \to \cald_{lf}^{C_i} ( \tilde{M} , T ; p )$
is an equivalence of categories. This is seen as follows: on both sides every object is isomorphic 
to one with support in $(\tilde{M}-T) \times \einsu$ and for morphisms between such objects 
the control conditions agree. 
\end{proof}




\typeout{---------------------versus ---------------------}


\section{Bounded versus $\e$-control}  \label{versus}

So far we always considered the whole category of morphisms
over a space and investigated the functorial properties
of such a construction. In the proof of surjectivity we
will shift our attention to individual morphisms in the category.

In this section we will develop a criterion which  decides that  an element
is in the image of the assembly map if it satisfies certain control conditions.
We start by developing some more language concerning control of morphisms.
 
\subsection{Control}

Below and in Section~\ref{squeezing}
it is important to use the right notion of ``control over a subset''.
For a module $M$ over $X$ and a subset $B \subset X$ we denote by $M|_B$ the 
largest submodule of $M$ with support in $B$. Let 
$i_B$ and $p_B$ the inclusion of respectively the projection onto $M|_B$.
For a morphism $\phi$ we then have the restriction $\phi|_B= p_B \circ \phi \circ i_B$.
Also recall the $E$-thickening of a subset from page~\pageref{Ethickening}.

\begin{definition}[Control] \label{controlover}
Let $X$ be a metric space or 
let $E \subset X \times X$ be a symmetric neighborhood of the diagonal.
Let $B \subset X$. Let $\phi$ and $\psi$ be morphisms in $\calc(X)$.
\begin{enumerate}
\item
A morphism $\phi$ is $\alpha$-{\em{controlled}} if $| \phi | \leq \alpha$.
A morphism $\phi$ is $E$-{\em{controlled}} if $\supp( \phi ) \subset E$.
\item
An automorphism $\phi$ in $\calc(X)$ is called an $\alpha$-{\em automorphism}
if $\phi$ and $\phi^{-1}$ are $\alpha$-controlled, it is called an $E$-{\em automorphism}
if $\phi$ and $\phi^{-1}$ are $E$-controlled.
\item
A morphism $\phi$ is $\alpha$-{\em controlled
over} $B \subset X$ if for every $x \in X$  and $y \in B$ with $d(x,y) > \alpha$
we have $\phi_{x,y}= \phi_{y,x}=0$.
A morphism $\phi$ is said to be $E$-{\em controlled over} $B$ if 
\[
\supp( \phi \circ i_B ) \subset E \quad  \mbox{ and } \quad \supp( p_B \circ \phi  ) \subset E
\]
or equivalently $\phi \circ i_B = i_{B^E} \circ p_{B^E}\circ  \phi \circ i_B$ and
$p_B \circ \phi = p_B \circ \phi \circ i_{B^E} \circ p_{B^E}$.
\end{enumerate}
\end{definition}

{\bf Warning:} Control over $B$ is not a local notion, i.e.\ we can not 
compute the control of $\phi$  over $B$ from the knowledge of $\phi|_B$.

\subsection{An explicit description of $K_2$}

>From the germs at infinity fibration (cf.\ Subsection~\ref{germsatinfty}) and an Eilenberg swindle (cf.~\ref{swindletoinfty}) 
we know abstractly that there is an isomorphism
\[
\xymatrix{
K_2 \cald^{\Gamma}(X, \ast ) \ar[r]^{\simeq} &  K_1 \calc^{\Gamma}( X ).
         }
\] 
But in order to translate control conditions we need a more explicit
description of such an isomorphism. We  follow \cite{Pedersen(2000)}.
In particular compare Section~5 of that paper. 

Let us  briefly recall
how elements in $K_2$ of an additive category $\cala$ can be described.
Compare~\cite[p.360]{Pedersen(2000)} and \cite[p.4]{Weibel(1981)}. 
The description is completely analogous to the definition 
of $K_2$ of a ring as the kernel of the homomorphism $\St(R) \to E(R)$.
An automorphism $e:M\to M$ in $\cala$ is called elementary if there
is a decomposition $M = M_0 \oplus  \dots \oplus M_k$ such that 
with respect to this decomposition $e$ has only one nonzero off diagonal entry
and identities on the diagonal. A deformation $\eta = (e_1, \dots ,e_l)$ 
is a finite sequence of elementary automorphisms with respect to a fixed decomposition
of a fixed module.
One usually thinks  of such a sequence as a formal product. 
We denote by $\overline{\eta}$ the automorphism
associated to $\eta$, i.e.\ we have a map
\[
\eta = ( e_1, \dots , e_l) \mapsto \overline{\eta} = \prod e_i.
\]
An element 
in $K_2 \cala  $ is represented by a deformation $\eta$ with $\overline{\eta}=\id$.
For an elementary automorphism $e=e_{\phi}$ whose off diagonal term is $\phi$
let $e|_A$ be the same elementary automorphism with $\phi$ replaced by $\phi|_A$.
For a deformation $\eta=(e_1, \dots , e_l)$ define $\eta|_A=(e_1|_A, \dots ,e_l|_A)$.
The control notions have their obvious analogues for deformations.

\begin{definition} \label{deformationcontrol}
Let $X$ be a metric space or let $E \subset X \times X$ be a symmetric neighborhood of the diagonal.

A deformation
$\eta=( e_1, \dots , e_l )$ is  $\alpha$-{\em controlled} if $\sum |e_i| \leq \alpha$.
The deformation is $E$-{\em controlled} if all partial products of the $e_i$
are  $E$-controlled. 
\end{definition}

We now review a construction from \cite{Pedersen(2000)} 
which shows how a sequence $(\phi_i,\eta_i)$ of automorphisms and deformations 
over a $\Gamma$-space $X$ gives rise to an element 
$a(\phi_i,\eta_i) \in K_2(\cald^\Gamma(X;*))$.

\begin{definition}[Bounded sequences] 
Let $\{ \phi_i : M_i \to M_i | i \in \IN \}$ be a sequence of
automorphisms in $\calc^{\Gamma}(X)$. Suppose $M_{i+1} = M_i \oplus L_i$ 
and there are deformations $\eta_i$
such that $(\phi_i \oplus \id_{L_i}) \overline{\eta}_{i+1} = \phi_{i+1}$. 
We will say $(\phi_i,\eta_i)$ is a {\em bounded sequence} if 
there is an $\alpha > 0$ such that each $\phi_i$ is an $\alpha$-automorphism 
and each $\eta_i$ is $\alpha$-controlled.
\end{definition}

Assume that $X$ is $\Gamma$-compact.
Use the inclusion $X \x \{i\} \hra X \x [1,\infty)$ to consider $M_i$ as a 
module over $X \x [1,\infty)$. In this way $M = \bigoplus_i M_i$ is a module
over $X \x [1,\infty)$. Now define 
$\overline{\eta}_{\mathit{ev}},\overline{\eta}_{\mathit{odd}},\overline{\phi}_{\mathit{ev}}$ and 
$\overline{\phi}_{\mathit{odd}}$ as indicated in the following diagram:
\[
\xymatrix
{
 M \ar[d]^{\overline{\phi}_{\mathit{even}}} & &
 M_1 \ar[d]^\id & 
 M_2 \ar[d]^{\phi_2} & 
 M_3 \ar[d]^{\phi_2^{-1} \oplus \id} & 
 M_4 \ar[d]^{\phi_4} & 
 M_5 \ar[d]^{\phi_4^{-1} \oplus \id} &
 \dots
 \\
 \cdot \ar[d]^{\overline{\eta}_{\mathit{odd}}} & &
 \cdot \ar[d]^\id &
 \cdot \ar[d]^{\id} &
 \cdot \ar[d]^{\overline{\eta}_3^{-1}} &
 \cdot \ar[d]^{\id} &
 \cdot \ar[d]^{\overline{\eta}_5^{-1}} &
 \dots
 \\
 \cdot \ar[d]^{\overline{\phi}_{\mathit{odd}}} & &
 \cdot \ar[d]^{\phi_1} &
 \cdot \ar[d]^{\phi_1^{-1} \oplus \id} &
 \cdot \ar[d]^{\phi_3} &
 \cdot \ar[d]^{\phi_3^{-1} \oplus \id} &
 \cdot \ar[d]^{\phi_5} &
 \dots 
 \\
 \cdot \ar[d]^{\overline{\eta}_{\mathit{even}}} & &
 \cdot \ar[d]^{\id}  &
 \cdot \ar[d]^{\overline{\eta}_2^{-1}} &
 \cdot \ar[d]^{\id} &
 \cdot \ar[d]^{\overline{\eta}_4^{-1}} &
 \cdot \ar[d]^{\id} &
 \dots  
 \\
 M & & M_1 & M_2 & M_3 & M_4 & M_5 & \dots
 \\
}  
\]
Recall that the matrix
\[
\left( \begin{array}{cc}
           \psi  & 0 \\
           0     & \psi^{-1}  
       \end{array} \right)
\]
can be written as the following product
of elementary matrices
\[
\left( \begin{array}{cc}
           1     & 1 \\
           0     & 1         
       \end{array} \right)
\left( \begin{array}{cc}
           1     & 0 \\
           -1    & 1         
       \end{array} \right)
\left( \begin{array}{cc}
           1     & 1 \\
           0     & 1         
       \end{array} \right)
\left( \begin{array}{cc}
           1     & 0 \\
           \psi  & 1         
       \end{array} \right)
\left( \begin{array}{cc}
           1     & -\psi^{-1} \\
           0     & 1         
       \end{array} \right)
\left( \begin{array}{cc}
           1     & 0 \\
           \psi  & 1         
       \end{array} \right).
\]
This identity can be used to factor 
$\overline{\phi}_{\mathit{ev}}$ and $\overline{\phi}_{\mathit{odd}}$ into products of 
elementary automorphisms. Clearly, $\overline{\eta}_{\mathit{ev}}$ and 
$\overline{\eta}_{\mathit{odd}}$ are also products of elementary automorphisms. 
(These products are only finite if the number 
of elementary automorphisms in the deformation $\eta_i$ is 
uniformly bounded. However, by inserting more
terms in our sequence one can easily arrange that each $\eta_i$ itself is an 
elementary automorphism.) Note that 
\[
 \overline{\eta}_{\mathit{even}} \circ 
       \overline{\phi}_{\mathit{odd}} \circ
 \overline{\eta}_{\mathit{odd}} \circ \overline{\phi}_{\mathit{even}}
 = \phi_1 \oplus \id \oplus 
           \id \oplus \dots.
\]
In $\cald^\Gamma(X;*)$ this is equivalent to $\id_M$. 
Hence we have constructed a deformation whose associated
automorphisms is $\id_M$. This defines $a(\phi_i,\eta_i) \in K_2 \cald^\Gamma(X;*)$.

Recall that for $\Gamma$-compact and free $X$ we have by definition 
\[
\cald^{\Gamma}( X ; \ast) = \calc^{\Gamma}( X \times \einsu; p_{\einsu}^{-1} \cale_d )^{\infty}
\]
where $\cale_d$ is the standard metric control condition on $\einsu$.
Consider the map 
\begin{eqnarray*}
P : K_2 \cald^\Gamma(X;*) & \to & K_1 \calc^{\Gamma}(X) \\
\eta=(e_1 , \dots , e_l) & \mapsto & (p_{X})_{\ast} \left( \overline{\eta}|_{X \times \left[ 1, r \right]} \right) = ( p_X)_{\ast} \left(
\left( \prod e_i \right)|_{X \times \left[ 1, r \right]} \right)
\end{eqnarray*}
We may assume that each $e_i$ is represented
by a morphism which is already elementary before taking the germs at infinity.
On the right hand side we restrict to $X \times \left[ 1, r \right]$
for some large enough $r$ which depends on the element and then 
project down to $X$ via $p_X$. Compare \cite[5.1]{Pedersen(2000)}.

\begin{lemma} \label{pedersensmap}
The map $P$ is an isomorphism.
\end{lemma}

\begin{proof} 
Using Steinberg relations, enlarging $r$ or making a different choice for the 
representatives of the $e_i$ in $\calc^{\Gamma}(X \times \einsu ; p^{-1}_{\einsu} \cale_d )$ does not change $P(\eta )$.
We have a well defined homomorphism. 
Surjectivity is clear since $P(a(\phi_i,\eta_i)) = [\phi_1]$.
Suppose $( \prod e_i)|_{[1,r]}$ is a product of elementary matrices. Complete these 
elementary matrices via the identity
to elementary matrices $f_1, \dots,f_m$ in $\calc^{\Gamma}(X \times \einsu ; p^{-1}_{\einsu} \cale_d )$. In 
the germ category each $f_i$ represents the identity
and hence $a=[(e_1, \dots , e_l)]=[(e_1,\dots,e_l,f_m^{-1},\dots f_1^{-1})]$ in
$K_2\calc^{\Gamma}(X \times \einsu ; p^{-1}_{\einsu} \cale_d )^{\infty}$. Now 
$\prod e_i \prod f_j^{-1}=\id$ in $\calc^{\Gamma}(X \times \einsu ; p^{-1}_{\einsu} \cale_d )$ by construction
and hence $a$ lies in the image of 
\[
K_2 \calc^{\Gamma}( X \times \einsu ; p^{-1}_{\einsu} \cale_d  ) \to 
 K_2 \calc^{\Gamma} (X \times \einsu ; p^{-1}_{\einsu} \cale_d )^{\infty}.
\]
But $K_2 \calc^{\Gamma} ( X \times \einsu ; p^{-1}_{\einsu} \cale_d ) =0$ by the usual Eilenberg swindle \ref{swindletoinfty}
and hence the map is injective.
\end{proof}

\begin{remark} It is conceivable that $P$ coincides with the boundary map
in the long exact sequence derived
from the germs at infinity fibration sequence in Subsection~\ref{germsatinfty} but we do not need this.
\end{remark}

\subsection{Deciding that an element is in the image}

Let $X$ be a Riemannian manifold.
Let $\Gamma$ act freely, properly discontinuously 
and cocompactly by isometries on $X$. Let $A$ be closed  
$\Gamma$-invariant submanifolds of $X$. Let $X(A)$ be the space obtained from
$X$ by collapsing each path component of $A$ to a point, i.e.\ we have a 
pushout diagram
\[
\xymatrix{
A \ar[d] \ar[r] & X \ar[d]^{p_{X(A)}} \\
\pi_0 (A) \ar[r] & X(A).
         }
\]
We now develop a criterion which decides that an element is in the image
of the map $K_2 \cald^{\Gamma} ( X , p_{X(A)} ) \to K_2 \cald^{\Gamma} (X , \ast)$.

\begin{definition}[Controlled sequence]  \mbox{ }
\label{controlseq}
\begin{enumerate}
\item A morphism $\phi\in \calc^{\Gamma}(X)$ is said to be $p_{X(A)}$-{\em separating}
if for $x$, $y$ from different path components of $A$ we have $\phi_{x,y}=0$.
\item
Let $(\phi_i, \eta_i )$ be a bounded sequence, where $\eta_i = (e_1^i, \dots , e_{l_i}^i)$. We will
say that $(\phi_i, \eta_i)$ is a $p_{X(A)}$-{\em controlled sequence} if the following
holds:
\begin{enumerate}
\item
All $\phi_i$, $\phi_i^{-1}$ and $e_j^i$ are $p_{X(A)}$-separating.
\item
There are $\epsilon_i>0$ such that
$\phi_i$, $\phi_i^{-1}$ and $\eta_i$ are $\epsilon_i$-controlled over $X-A$ and $\lim \e_i=0$.
\end{enumerate}
\end{enumerate}
\end{definition}

\begin{lemma}
\label{seqintheimage}
Let $(\phi_i,\eta_i)$ be a $p_{X(A)}$-controlled sequence,
then 
$P^{-1}([\phi_1])$ is in the image 
of the assembly map
\[
K_2 \cald^{\Gamma}( X , p_{X(A)} ) \to K_2 \cald^{\Gamma} ( X , \ast ).
\]
\end{lemma}

\begin{proof} 
By definition a 
$p_{X(A)}$-controlled sequence $(\phi_i,\eta_i)$ gives an element 
$\bar{a}(\phi_i,\eta_i)$ in $K_2(\cald^\Gamma(X;p_{X(A)}))$ 
which maps to 
$a(\phi_i,\eta_i) = P^{-1}( \left[ \phi_1 \right] )
 \in K_2(\cald^\Gamma(X;*))$ under
the map induced by 
$\cald^\Gamma(X;p_{X(A)}) \to \cald^\Gamma(X;*)$.
\end{proof}

The following criterion which will be used in the proof of Theorem~\ref{surjectivity}
is a consequence of the Squeezing Theorem~\ref{quinnsqueezing}.

\begin{proposition}
\label{smallisintheimage} 
There is an $\e_0 > 0$ such that the following holds: Let
$\phi$ be an $p_{X(A)}$-separating automorphism over $X$, 
such that $\phi$ and $\phi^{-1}$ are $\e_0$-controlled
over $X-A$. Then $[\phi]$ is in the image of the composition
\[
\xymatrix
{
 K_2(\cald^{\Gamma} (X,p_{ X(A) })) \ar[r] &
 K_2(\cald^{\Gamma} (X,*)) \ar[r]^\cong_P &
 K_1(\calc^{\Gamma} (X)).
}   
\]
\end{proposition}
\label{phi bounded ????}

\begin{proof}
Let $A \subset U_1 \subset U_2 \subset \dots \subset U$ 
be a sequence of tubular neighborhoods of $A$. Let
$p_{X(U)} : X \to X(U)$ be the map obtained by collapsing
each path component of $U$ to a point. We can now use  
\ref{quinnsqueezing} to produce from
$\phi$ a $p_{X(U)}$-controlled sequence. 
The point is that we have to give up control over $U_n$ in the 
$n$-th squeezing step, but never outside of $U$. 
Note that by \ref{homotopyinvarianceCor} we have 
$K_2(\cald^\Gamma(X,p_{X(A)})) \cong 
 K_2(\cald^\Gamma(X,p_{X(U)}))$.
The claim follows now from \ref{seqintheimage}.
\end{proof}

\begin{addendum} \label{smalladdendum}
The same statement holds if we replace the Riemannian metric
on the $\Gamma$-compact space $X$ by another equivariant metric $d$ which generates the same
topology.
\end{addendum}
\begin{proof} By compactness there is for given $\e_0$ an $\delta$ such that $\delta$-control with respect to $d$ implies
$\e_0$-control with respect to the Riemannian metric.
\end{proof}


\typeout{---------------------transfer ---------------------}


\section{The transfer}
\label{sec_transfer}

Before we construct the transfer we will briefly recall some facts about the torsion of a 
self-homotopy equivalence of a chain complex. For a detailed account
see for example  Chapter 12 in \cite{Lueck(1989)}, in particular Example 12.17.

\subsection{Torsion of a self-homotopy equivalence}
\label{sec_torsion}

Let $\cala$ be a small additive category. 
We denote by $\ch \cala$ the category of bounded chain complexes in $\cala$. 
For a complex $C$ define the objects $C_\mathit{ev} = \oplus_{n} C_{2n}$ and 
$C_\mathit{odd} = \oplus_{n} C_{2n+1}$. 
The cone $\cone(f)$ of a self homotopy equivalence
$f : C \to C$ is the contractible chain complex given by
\[
(\cone(f), \dd_{\cone(f)}) = 
   (\Sigma^{1}C \oplus C,
    \left( \begin{array}{cc}
              -\dd_C & 0 \\ f & \dd_C
           \end{array}
    \right)).
\] 
Let $\gamma$ be a chain contraction and let $\Phi$ denote the canonical isomorphism
$\cone(f)_\mathit{ev} \cong \cone(f)_\mathit{odd}$. 
The torsion of $f$ is defined as
\[
\tau(f) := [\Phi \circ 
     (\dd_{cone(f)} + \gamma)
          |_{\cone{(f)}_\mathit{odd}}]
     \in K_1(\cala).
\]
It is independent of the choice of $\gamma$ and depends only on the chain
homotopy class of $f$. It has the logarithmic property $\tau(f \circ g)= \tau(f) + \tau(g)$
and $\tau(\id)=0$ and a ladder diagram of self-homotopy equivalences leads
to the usual additivity relation.\label{standardproperties}

Now suppose that $\cala \subset \calc(X)$ is a category of modules over some space $X$.
Later we will need to control the support of the automorphisms representing
$\tau(f)$ and its inverse. We therefore need an explicit description
of the contraction $\gamma$.

Let $g$ be a chain homotopy inverse of $f$, let $s$ be a chain homotopy from $gf$ to $\id_C$,
i.e.\ $gf-\id_C=\dd_Cs + s \dd_C$. Let $t$ be a chain homotopy from $\id_C$ to $fg$. Then
we have a chain homotopy 
\[
\left( \begin{array}{cc}
         s & g \\ 0 & t
       \end{array}
\right) 
\mbox{ from }
\left( \begin{array}{cc}
         \id_C & 0 \\ tf+fs & \id_C
       \end{array}
\right)
\]   
to the zero map. Consequently
\begin{equation}
\label{explicit}
\gamma = \left( \begin{array}{cc}
         s & g \\ 0 & t
       \end{array}
\right)
\left( \begin{array}{cc}
         \id_C & 0 \\ -tf-fs & \id_C
       \end{array}
\right)
\end{equation}
is a chain contraction for $\cone(f)$. Moreover if $2N+2$ exceeds the dimension of $C$, then
an explicit inverse for $\dd_C + \gamma$ is given by
$(\dd_C + \gamma)(\id_C - \gamma^2+\gamma^4-\dots \pm \gamma^{2N})$.
\begin{remark} \label{estimateprinciple}
In particular the support of $(\dd_{cone(f)} + \gamma)|_{\cone{(f)}_\mathit{odd}}$ in $X \times X$
and its inverse can be obtained from the support of $f$, $g$, $s$ and $t$ by a finite
number of the operations ``union'' and ``composition of relations''.
\end{remark}
 
Given a bounded chain complex of finitely generated free $\IZ$-modules $C$ and 
a module $M$ in $\cala$ it makes sense to form the tensor product $M \otimes_{\IZ} C$
given by $(M \otimes_{\IZ} C)_x=M_x \otimes_{\IZ} C$.

\begin{lemma} \label{eulermultabstract}
Let $\varphi$ be an automorphism in $\cala$. Let 
$C$ be a finite complex of finitely generated free $\IZ$-modules, then 
\[
\tau( \varphi \otimes \id_C) = \chi( C ) \cdot \left[ \varphi \right]   \in  K_1 \cala
\]
\end{lemma}
\begin{proof}
This follows from the definition of $\tau$.
\end{proof}

\subsection{The transfer} \label{subsec_transfer}

Let $p : E \to B$ be a fiber bundle of Riemannian manifolds.
Let $\pi : \tilde{B} \to B$ be a regular cover with group of
deck transformations $\Gamma$.
Define $\tilde{E}$ as the pullback 
\[
\xymatrix{
\tilde{E} \ar[d]_{\pi_E} \ar[r]^{\tilde{p}} & \tilde{B} \ar[d]^{\pi} \\
E \ar[r]_{p} & B.
         }
\]
To construct the desired  transfer 
we need further structure.
For $x\in \tilde{B}$ let $\tilde{E}_x$ denote the fiber $\tilde{p}^{-1}( x  )$.
Suppose that we are given a fiber transport $\nabla$, i.e.\ a family of maps 
$\nabla_{y,x}: \tilde{E}_x \to \tilde{E}_y$ which is 
\begin{enumerate}
\item $\Gamma$-{\em invariant}, i.e.\ for all $\gamma \in \Gamma$ we have 
      $\gamma \circ \nabla_{y,x} \circ \gamma^{-1}
            = \nabla_{\gamma y,\gamma x}$.
\item $\nabla$ is {\em functorial}, i.e.\
      $\nabla_{x,x} = \id_{\tilde{E}_x}$ and 
      $\nabla_{z,y} \circ \nabla_{y,x} 
       = \nabla_{z,x}$.
\item And $\nabla$ induces a {\em homotopically trivial action}:
      $\nabla_{\gamma x,x}$ is homotopic to the deck transformation
      $\gamma : \tilde{E}_x 
       \to \tilde{E}_{\gamma x}$.   
\end{enumerate}

For every $q \in B$ we choose a triangulation $T_q$ of the fiber $E_q=p^{-1}(q)$. 
Let $C(T_q)$ denote the corresponding  cellular $\IZ$-chain complex.
For $x \in \tilde{B}$ the triangulation $T_{\pi(x)}$
induces a triangulation of $\tilde{E}_x$ under the 
canonical identification $\tilde{E}_x = E_{\pi(x)}$.   
For $x,y \in \tilde{B}$ let $\nabla^T_{y,x}$ be a cellular approximation 
of $\nabla_{y,x}$ with respect to these triangulations.
Note that $\nabla^T$ is not functorial anymore but since $\nabla^T_{y,x}$ is 
homotopic to $\nabla_{y,x}$ we know that $\nabla^T_{z,y} \circ \nabla^T_{y,x}$
is homotopic to $\nabla^T_{z,x}$. We can and will 
assume that  $\nabla^T$  is again $\Gamma$-invariant and that $\nabla^T_{x,x}=\id_{\tilde{E}_x}$.

\begin{proposition} \label{abstracttransfer}
There exists a transfer depending on $\nabla^T$ 
\[
tr_{\nabla^T}: \calc^{\Gamma}( \tilde{B} , \calf_{\Gamma c}) \to 
\ch \calc^{\Gamma}( \tilde{E} , \calf_{\Gamma c } )
\]
which need not be a functor but has the following properties.
\begin{enumerate}
\item It is functorial up to homotopy, i.e.\ $\tr_{\nabla^T} ( \varphi \circ \psi )$
is chain homotopic to $\tr_{\nabla^{T}} ( \varphi ) \circ \tr_{\nabla^{T}} ( \psi )$
and $\tr_{\nabla^T}(\id)=\id$. In particular it sends automorphisms to self-homotopy equivalences.
\item
Torsion and transfer induce a well defined map
\[
\tau \circ \tr_{\nabla^T}:K_1 \calc^{\Gamma}( \tilde{ B } , \calf_{\Gamma c} )  \to 
K_1 \calc^{\Gamma}( \tilde{E} , \calf_{\Gamma c}).
\]
\item \label{eulermult}
For an automorphism $\varphi$ in $\calc^{\Gamma}( \tilde{B} , \calf_{\Gamma c} )$ we have
in $ K_1 \calc^{\Gamma}( \tilde{B} , \calf_{\Gamma c} )$
\[
\tilde{p}_{\ast} \tau(\tr_{\nabla^T}(\varphi)) =  \chi(E_q) \cdot [\varphi],
\] 
where $\chi( E_q )$ is the Euler characteristic of the fiber.
\end{enumerate}
\end{proposition}
\begin{proof} 

Let $M$ be an object in $\calc^{\Gamma} (\tilde{B} , \calf_{\Gamma c})$.
We construct the chain complex $\tr_{\nabla^T}(M)$ over $\calc^{\Gamma} (\tilde{E} ,  \calf_{\Gamma c})$
as follows. Define 
\begin{equation}
\label{transferM}
\tr(M)= \bigoplus_{x \in \tilde{B}} 
          M_{x} \otimes_{\IZ} C(T_{\pi(x)}).
\end{equation}
The differential is given by 
\[ 
\partial_{\tr(M)} = 
            \bigoplus_{x \in \tilde{B}} 
          \id_{M_{x}} 
             \otimes \partial_{C(T_{\pi(x)})}.
\]
Using the barycenters of the triangulations we turn $\tr_{\nabla^T}(M)_n$ into 
a module over $\tilde{E}$: Let $e \in \tilde{E}_x$.
If there is an $n$-simplex $\sigma$ in $T_{\pi(x)}$ having
$\pi_E(e)$ as its barycenter then
\[
(\tr(M)_n)_e = M_{\tilde{p}(e)} 
                   \otimes_{\IZ} \IZ \sigma,
\]
otherwise $(\tr(M)_n)_e = 0$. Since we used the triangulation of $E_{\pi(x)}$
instead of $\tilde{E}_x$ we obtain a $\Gamma$-invariant object.

In order to lift morphisms we need $\nabla^T$. Let $\varphi = (\varphi_{x,y})$ be a morphism in
$\calc^{\Gamma} (\tilde{B} , \calf_{\Gamma c})$.
We define $\tr_{\nabla^T}(\varphi)$ 
as a matrix with respect to the direct sum decomposition in (\ref{transferM}):
\[
\tr_{\nabla^T}(\varphi)_{y,x} =
    \varphi_{y,x} \otimes 
           C(\nabla^T_{y,x}).
\]
Here we consider $\nabla^T_{y,x}$ as a map $E_{\pi(x)} \to E_{\pi(y)}$ 
using the identifications $\tilde{E}_x=E_{\pi(x)}$  
and $C(\nabla^T_{y,x}): C(T_{\pi(x)}) \to C(T_{\pi{y}})$ is the induced map
of chain complexes.
Since $\nabla^T$ is functorial up to homotopy we get that  
$\tr_\nabla(\varphi \circ \psi)$ is always chain homotopic to 
$\tr_\nabla(\varphi) \circ \tr_\nabla(\psi)$. 

(ii) follows from the standard properties of $\tau$, see page~\pageref{standardproperties}.

Let us now prove \ref{eulermult}.
Fix a point $x \in \tilde{B}$.
By Lemma~\ref{XversusGamma} the inclusion $\calc^{\Gamma}( \Gamma x) \to 
\calc^{\Gamma}( \tilde{B} , \calf_{\Gamma c} )$ is an equivalence
(and both categories are equivalent to the category of finitely generated free
$R\Gamma$-modules). We may therefore assume that we have an automorphism
$\varphi= (\varphi_{hx,gx}): M \to M$ with $\supp M \subset \Gamma x$.
Since $\nabla$ induces a homotopically trivial action we know that $\nabla_{hx, gx}$ and therefore
$\nabla_{hx, gx}^T$ is homotopic to the deck transformation given by $h^{-1}g$.
It follows that $\tr_{\nabla^T}( \varphi )_{hx , gx}$ is homotopic to 
$\varphi_{hx,gx} \otimes \id_{C(T_{\pi(x)})}$. This implies that $p_{\ast} \tau \tr( \varphi )
= \tau p_{\ast} \tr ( \varphi ) = \tau( \varphi \otimes \id_{C(T_{\pi(x)})})$.
The result follows by Lemma \ref{eulermultabstract}.
\end{proof}

\subsection{Transfer and flow}

We still assume the situation described in the beginning of the last subsection.
Let $\Phi_t$ be a $\Gamma$-invariant flow on $\tilde{E}$. 
The flow lines of $\Phi_t$ induce a foliation 
with one dimensional leaves on $\tilde{E}$. 
We will need the concept of foliated or $(\alpha,\delta)$-control.
Usually $\alpha$ will be large and $\delta$ will be small.
 
\begin{definition}[$(\alpha,\delta)$-control]
\label{alphadeltacontrol}
Let $\alpha,\delta > 0$ and $v,w \in \tilde{E}$. 
The {\em pair} $(v,w)$ is
said to be $(\alpha,\delta)$-{\em controlled} if there is a 
flow line $L$ and a path $\omega : [0,1] \to L$ of 
arclength at most $\alpha$, such that $dist_{\tilde{E}}(v,\omega(0)) < \delta$
and $dist_{\tilde{E}}(\omega(1),w) < \delta$. 
Let $E_{\alpha,\delta} \subset \tilde{E} \times \tilde{E}$ be the set of all
$(\alpha, \delta)$-controlled pairs. This is a symmetric neighborhood
of the diagonal and the definition of $E$-control, $E$-control over a subset etc. from 
\ref{controlover} and \ref{deformationcontrol} apply.
For brevity we usually write  $(\alpha, \delta)$-controlled instead of $E_{\alpha,\delta}$-controlled etc.. 
\end{definition}

\begin{definition}[$\Phi_t$ contracts $\nabla$]
\label{foliated_contracted}
We say that the flow $\Phi_t$ contracts the fiber transport $\nabla$
if there is a function $f$ such that for every $\alpha>0$ there exists an assignment $t \mapsto \delta_t>0$ for $t \in \einsu$
with $\lim_{t \to \infty} \delta_t = 0$ such that the following holds:
\begin{quote}
For all $x,y \in \tilde{B}$ with 
$dist_{\tilde{B}}(x,y) < \alpha$ and every $e \in \tilde{E}_y$ 
the pair 
\[
(\Phi_t(\nabla_{x,y}(e)),\Phi_t(e))
\] is
$(f(\alpha),\delta_t)$-controlled.
\end{quote}
\end{definition}  
Later in the application the function $f$ will in fact be of the form $f(\alpha)=D_1 \alpha + D_2$ with constants
$D_1$, $D_2 >1$. Compare~\ref{HeinzImHof} and~\ref{Ballmann}.

\begin{proposition}[Gaining foliated control]
\label{transfer}
Suppose we are in the situation described in the beginning of Subsection~\ref{subsec_transfer}.
Suppose that $\Phi_t$ is a $\Gamma$-invariant flow on $\tilde{E}$ which contracts $\nabla$. 
Fix an element $\left[ \varphi \right] \in 
 K_1 \calc^{\Gamma} (\tilde{B} , \calf_{\Gamma c} )$. 
Then there is a constant $\alpha > 0$ such that for every 
$\delta > 0$ there is an $(\alpha, \delta)$-automorphism 
$\psi$ in $\calc^{\Gamma}(\tilde{E}, \calf_{\Gamma c})$ 
satisfying
\[
\tilde{p}_{\ast} [\psi] = \chi(E_q)  \cdot \left[ \varphi \right] \in 
 K_1\calc^{\Gamma} ( \tilde{B}, \calf_{\Gamma c}). 
\] 
Here $\chi(E_q)$ is the Euler-characteristic of the fiber.
\end{proposition}

\begin{proof}
Fix $x \in \tilde{B}$. Again we may assume by 
Lemma~\ref{XversusGamma} that $\varphi : M \to M$ is an automorphism
in $\calc^{\Gamma} (\tilde{B} , \calf_{\Gamma c})$
with $\supp(M) \subset \Gamma x$. 
Let $\alpha_0 > 0$ be such that $|\varphi| \leq  \alpha_0$ and $|\varphi^{-1}| \leq \alpha_0$.
We have to choose suitable triangulations and cellular 
approximations in order to use the transfer $\tr_{\nabla^T}(\varphi)$. 

Fix $\delta_0>0$. By definition of morphisms the set
\[
S = \{\gamma \in \Gamma|
         \varphi_{\gamma x,x} \neq 0
         \mbox{ or }
         (\varphi^{-1})_{\gamma x,x} \neq 0 
      \}
\] 
is finite,
hence there is ${t_0}>0$ such that
\[
(\Phi_{t_0}(\nabla_{\gamma x,x}(e)),\Phi_{t_0}(e))
\]
is $(f(\alpha_0),\delta_0)$-controlled for all 
$\gamma \in S$ and all $e \in \tilde{E}_{x}$.
Let $T_{x}$ be a triangulation of 
$\tilde{E}_{x}$. 
(This gives also triangulations of 
$E_{\pi(x)}$
and
$\tilde{E}_{\gamma x}$.) 
For $\gamma \in \Gamma$
let 
$\nabla^T_{\gamma x,x}$
be a cellular approximation of 
$\nabla_{\gamma x,x}$. By choosing the
triangulation $T_{x}$ sufficiently fine
and the approximation sufficiently close,
we can arrange that 
there is a homotopy 
$H_\gamma : \tilde{E}_{x} \x [0,1] \to
            \tilde{E}_{\gamma x}$
between $\nabla_{\gamma x,x}$ and 
$\nabla^T_{\gamma x,x}$ satisfying
\[
\mathit{dist}_{\tilde{E}}
     (\Phi_{t_0}(H_\gamma(e,s)),
      \Phi_{t_0}(\nabla_{\gamma x,x}(e)))
\leq \delta_0
\]
for all $s \in [0,1]$.
This implies the following:

$(\Phi_{t_0})_*(\tr_\nabla(\varphi))$ and $(\Phi_{t_0})_*(\tr_\nabla(\varphi^{-1}))$ 
are $(f(\alpha_0),\delta_0)$-controlled. Moreover, their products (both ways)
are chain homotopic to the identity  via $(f(\alpha_0),\delta_0)$ controlled chain homotopies.
Remark~\ref{estimateprinciple} implies that $\tau(\tr_\nabla(\varphi))$ can be realized by
an $(k f(\alpha_0),k \delta_0)$-automorphism $\psi$, where $k$ is a constant that depends only on 
the dimension of $\tilde{E}_{x}$. Now the result follows from
\ref{eulermult}.
\end{proof}


\typeout{---------------------foliated ---------------------}


\section{Squeezing Theorems} \label{squeezing}

The Foliated Squeezing Theorem~\ref{foliatedcontrol} below will play an important role in 
the proof of the surjectivity result~\ref{surjectivity}. Roughly such a theorem
says the following: Once an automorphism has sufficient foliated control then we can find 
a representative in the same $K$-theory class with arbitrarily good
ordinary control.  At least this is true away from ``short'' closed geodesics.

Analogous statements for the Whitehead group were proven in \cite{Farrell-Jones(1986)}.
The same geometric arguments (in particular
the existence of a long and thin cell-structure) lead to such a statement in our context once
a suitable version of an ``ordinary'' Squeezing Theorem is available.
Indeed such an ordinary Squeezing Theorem was proven by Quinn.

\subsection{Quinn's Squeezing Theorem}

Roughly speaking an ordinary Squeezing Theorem tells us that once a certain
amount of control ($\epsilon_{g}$-control=''good control'') has been obtained, 
then one can improve the control to arbitrarily good control ($\epsilon_{vg}$-control=''very good control'').
Moreover if one starts out with control only over a part of the space (below this will be $X-S$) 
then the procedure still works but maybe one has to give up a little bit of control near the 
boundary. More precisely:

\label{alpha????}
\begin{theorem}[Squeezing]
\label{quinnsqueezing}  
Let $X$ be a compact Riemannian manifold. 
Let $\pi : \tilde{X} \to X$ denote the universal cover and let $\Gamma$ denote the 
fundamental group. Let $K$ and $S$ be closed subsets of $X$ with $S \cap K = \emptyset$.
Then there is $\e_0 = \e_0(X,K,S)$ and a homeomorphism 
$r = r(X,K,S) : [0,\infty) \to [0,\infty)$ such that the following holds: 
\begin{quote}
Let $\e_0 \geq \e_g \geq \e_{vg} \geq 0$ and $\alpha > r(\e_g)$ and 
let $\phi : M \to M$ in $\calc^{\Gamma}(\tilde{X})$ 
be an $\alpha$-automorphism. Assume moreover that $\phi$ and $\phi^{-1}$ are 
$\e_g$-controlled over $\tilde{X} - \pi^{-1}(S)$.   

Then there is a stabilizing module 
$L$ and an deformation 
$\eta = (e_1,\dots,e_n)$ on $M \oplus L$ 
in $\calc^\Gamma(\tilde{X})$ such that: 
\begin{enumerate}
\item The deformed automorphism $\phi_{new} = (\varphi \oplus id_N) \eta$ is 
      an $\alpha$-automorphism. Moreover, $\phi_{new}$ and $\phi_{new}^{-1}$ 
      are $\e_{vg}$-controlled over $\pi^{-1}(K)$.
\item The deformation $\eta$ is $r(\e_g)$-controlled and each $e_i$ is the identity on
      $(M \oplus N)|_{\pi^{-1}(S)}$. 
\end{enumerate}
\end{quote}
\end{theorem} 
  
\begin{proof}
This is almost \cite[4.5]{Quinn(1982)} and can be proven along the same lines.
Compare also the proof in~\cite{Bartels-Farrell-Jones-Reich(2002)} which is based on
\cite[3.6]{Pedersen(2000)}.
\end{proof}

\subsection{A foliated Squeezing Theorem}

Let $Y$ be a not necessarily compact  Riemannian manifold. 
Let $\pi : \tilde{Y} \to Y$ denote the universal cover and $\Gamma$ the fundamental group.
Let $\Phi: \IR \times Y \to Y$ be a smooth flow satisfying the following condition.
\begin{quote}  
For every $\alpha > 0$ and every compact set $K \subset Y$ there are only finitely many
closed orbits of length less than $\alpha$ which meet $K$.  
\end{quote}

Recall from \ref{alphadeltacontrol} that $E_{\alpha,\delta}$ is a symmetric neighborhood
of the diagonal and hence we have the notion of $E_{\alpha, \delta}$-control
which we abbreviate as $(\alpha,\delta)$-{\em control}. Compare \ref{controlover} and \ref{deformationcontrol}.

\begin{theorem}[Foliated Squeezing]
\label{foliatedcontrol}
There is a constant 
$\mu_1$ which only depends on the dimension of $Y$
for which the following statement is true:

\begin{quote}
Let $\alpha_0 > 0$ be an arbitrarily large number.
Let $g_1, g_2, \dots , g_N$ denote all the closed orbits 
of length less than $\mu_1 \alpha_0$. Let $K$ be an arbitrary compact subset of $Y$
which does not meet $S = \bigcup_{i=1}^{N} g_i$. 
Then there exist
\begin{quote}
numbers $\delta_0  > 0$ and $\mu_2>1$  depending 
on $K$ and $\alpha_0$
\end{quote}
such that:
\begin{quote}
For any  $\epsilon>0$ and $(\alpha , \delta)$ with 
\[
\e \leq \alpha \leq  \delta_0 \mbox{ and }
\e \leq \delta \leq  \delta_0
\]
and every $(\alpha, \delta)$-automorphism 
$\phi:M \to M$ in $ \calc^{\Gamma} ( \tilde{Y} , \calf_{\Gamma c})$ 
there exists 
a stabilizing module $L$ over $\tilde{Y}$ and a deformation 
$\eta$ on $M\oplus L$  in $\calc^{\Gamma}(\tilde{Y}, \calf_{\Gamma c})$
such that:
\begin{enumerate}
\item
The deformed automorphism 
$\phi_{new} = \overline{\eta} (\phi \oplus \id_L )$ and its inverse
$\phi_{new}^{-1}= (\phi^{-1} \oplus \id_L) \overline{\eta}^{-1}$ are both
$\epsilon$-controlled over $\pi^{-1}(Y_2)$.
\item
The deformation $\eta$, $\phi_{new}$ and its inverse $\phi_{new}^{-1}$ are all
everywhere $(\mu_1\alpha,\mu_2\delta)$-controlled. 
\end{enumerate}
\end{quote}
\end{quote}
\end{theorem}

\begin{proof}
For controlled $h$-cobordisms an analogous statement is  \cite[1.6]{Farrell-Jones(1986)}.
For controlled pseudoisotopies see \cite[0.2]{Farrell-Jones(1987a)}. The present statement can be proven 
with similar arguments. A detailed proof will appear in \cite{Bartels-Farrell-Jones-Reich(2002)}.
\end{proof}


\typeout{---------------------geometryII ---------------------}


\section{Geometric preparations needed for proving surjectivity}
\label{sec_geometry}

We now collect the geometric facts, to a large extent following \cite{Farrell-Jones(1986)}, which will
be needed in the proof of surjectivity.

Let $N$ be a simply connected Riemannian manifold whose sectional curvature
$K$ satisfies
\[
-b^2 \leq K \leq -a^2 < 0,
\]
where $a$ and $b$ are positive constants.
Geodesic rays in $N$ are called asymptotic if there distance is bounded. 
The sphere at infinity is by definition
the set of equivalence  classes of unit speed geodesic rays, where two rays 
are equivalent if they are asymptotic (compare \cite{Eberlein-ONeill(1973)}).
Let $SN \subset TN$ be the unit sphere subbundle
of the tangent bundle.
For each point $x \in N$ the fiber
$SN_x$ of the bundle $SN \to N$ 
is naturally homeomorphic to the sphere at infinity by 
the map which sends $v \in SN_x$ to the geodesic ray it determines.
This leads to a family of homeomorphisms 
\[
\nabla_{y,x}: S N_x \to SN_y
\]
which we call the {\em asymptotic fiber transport}. 

Our aim now is to verify that using the geodesic flow we can achieve foliated control:

\begin{proposition} \label{nablacontr}
The geodesic flow on $SN$ contracts the asymptotic fiber transport $\nabla$.
Compare Definition \ref{foliated_contracted}.
\end{proposition}

\begin{proof}
This follows by combining \ref{HeinzImHof} and \ref{Ballmann} below.
\end{proof}

Recall that associated to a point $\theta$ on the sphere at infinity is a Busemann function 
$F : N \to \IR$ which is unique up to an additive constant. 
The level surfaces of $F$ are the horospheres with center $\theta$ and 
$Z = -\mathit{grad} F$ is the unit length vector field on $N$
pointing towards $\theta$. In particular we have $\nabla_{y,x}(Z(x)) = Z(y)$. Clearly $Z$ determines a flow $\Psi_t$
on $N$. If we view $Z$ as a section of the bundle $SN \to N$ and hence as an embedding of $N$ into $SN$ 
then this is just the geodesic flow.
We will need the following consequence of results from \cite{Heintze-ImHof(1978)}.

\begin{proposition} \label{HeinzImHof}
Let $x,y \in N$. Let $\calh$ be the unique horosphere with center $\theta$ that contains $y$.
Let $z = \Psi_s(x)$ be the intersection point of the geodesic determined by $Z(x)$ 
with $\calh$. Assume that $s \geq 0$, i.e. $F(x) \geq F(y)$. Let $\alpha = d(x,y)$,
$c(t) = d(\Psi_t(x),\Psi_t(z))$ and $d(t) = d(\Psi_t(z),\Psi_t(y))$, where $d$ is the distance on $N$.
Then the following holds 
\label{lesenH!!!}
\begin{enumerate}
\item $c(t) \leq  \alpha + \frac{1}{a}$,
\item $d(t) \leq  \frac{2}{b} \sinh(\frac{b}{2} (\alpha + \frac{1}{a}) ) \exp(-at)$.
\end{enumerate}
\end{proposition}

\begin{proof} 
Let $\gamma : [0,\alpha] \to N$ be the geodesic from $x$ to $y$.
Let $s$ be the minimal number in $[0, \alpha ]$ such that $w=\gamma(s)$
is a point on the horosphere $\calh$.
Then \cite[4.8]{Heintze-ImHof(1978)} implies $d(w,z) \leq 1/a$.  
Clearly $d(x,w) \leq \alpha$ and therefore  $c(0) \leq \alpha + 1/a$. 
$\Psi_t$ moves $x$ and $z$ along a common geodesic and
therefore $c(t)$ is constant. This implies the first inequality.
Using $d(0) \leq d(y,w) + d(w,z) \leq \alpha + 1/a$,
the second inequality is a direct consequence of \cite[4.1,4.6]{Heintze-ImHof(1978)}.
\end{proof}

The Riemannian metric on $N$ induces a Riemannian metric on $SN$, cf.~(\ref{tangentmetric}). 
Let $d$ and $d^*$ denote the corresponding distance functions on $N$ and $SN$. 
The following Lemma finishes the proof of \ref{nablacontr}.

\begin{lemma} \label{Ballmann} 
There is a constant $C$ depending only on curvature bounds of $N$ such that
\[
d^*(v,\nabla_{y,x}v) \leq C d(x,y)
\]
for $x,y \in N$ and $v \in SN_x$. 
\end{lemma}

\begin{proof}
This is \cite[1.1]{Ballmann-Brin-Eberlein(1985)}. 
\end{proof}

>From now on let $M$ be our compact Riemannian manifold with strictly negative sectional curvature. 
Let $\Gamma = \pi_1(M)$. Let $\pi : \tilde{M} \to M$ denote the universal cover. 
Due to the Euler characteristic of the fiber which appears in \ref{transfer} we can not 
use the sphere bundle $S\tilde{M} \to \tilde{M}$ itself in our construction.
The problem is resolved by working with the northern-hemisphere bundle over the 
hyperbolic enlargement of $\tilde{M}$.
The hyperbolic enlargement of $\tilde{M}$ is the warped product 
\[
\IH \tilde{M} = \IR \times_{cosh(t)} \tilde{M}.
\]
Compare \cite[Section~7]{Bishop-ONeill(1969)}. 
It is the differentiable manifold $\IR \times \tilde{M}$ equipped with the Riemannian metric 
determined by $dg_{\IH \tilde{M}}^2= dt^2 + cosh(t)^2 dg_{\tilde{M}}^2$. 
 Since $M$ is compact $\tilde{M}$ has strictly negative sectional curvature bounds and the same
follows for $\IH \tilde{M}$, compare \cite[2.1.(vi)]{Farrell-Jones(1986)}.
Note that $\tilde{ \IH M} = \IH (\tilde{M})$.
Further properties of this construction can be found in Section~2 of \cite{Farrell-Jones(1986)}.
Note that the discussion in the beginning of this section applies with $N=\IH\tilde{M}$ to  
the sphere bundle
\[
p_{\IH \tilde{M}}: S \IH \tilde{M} \to \IH \tilde{M}.
\] 
For convenience  we introduce the notation
\[
\IH \tilde{M}_{\left[ -t , t \right]} 
= \left[ -t , t \right] \times \tilde{M}  \subset \IR \times_{cosh(t)} \tilde{M} 
\quad  \mbox{ and } \quad 
S\IH \tilde{M}_{[-t,t]}=p_{\IH \tilde{M}}^{-1}( \IH \tilde{M}_{\left[ -t , t \right]} ).
\]
Note that these spaces are 
$\Gamma$-compact. We use analogous notation for other subsets of $\IR$ instead of $[-t,t]$.

The reason to introduce the hyperbolic enlargement is that the extra distinguished direction
allows to introduce the northern-hemisphere subbundle $S^+ \IH \tilde{M}$ of the sphere-bundle
$S\IH \tilde{M}$. Fix a point $z=(0,x)$ in $\{0\} \times \tilde{M} \subset \IH \tilde{M}$.
Use the identification $T\IH \tilde{M}_{z}= T\IR_0 \times T \tilde{M}_x= \IR \times T \tilde{M}_x$
to define $S^+ \IH \tilde{M}_{z}$ as the subspace of $S \IH \tilde{M}_z$ with positive 
$\IR$ coordinate. Now define $S^+ \IH \tilde{M}_y$ for arbitrary 
$y \in \IH \tilde{M}$ as the image of $S^+ \IH \tilde{M}_{z}$ under the fiber transport. 
This is independent of the chosen point $z$ and defines 
a fiber bundle whose fiber is a disk and hence has Euler characteristic $\chi=1$.
\[
\xymatrix{
S^+\IH \tilde{M} \ar[dr]_{p_{\IH M}}  & \subset &  S\IH \tilde{M} \ar[dl]^{p_{\IH M}}  \\
  & \IH \tilde{M} .
         }             
\]
The northern-hemisphere bundle is invariant under the geodesic flow, i.e.\
we have a flow $\Phi_t : S^+ \IH \tilde{M} \to S^+ \IH \tilde{M}$. The 
asymptotic fiber transport can be restricted to this subbundle, i.e.\
we have a family of homeomorphisms
\[
\nabla_{y,x} : S^+\IH \tilde{M}_{x} 
              \to S^+\IH \tilde{M}_{y}.
\]
The following verifies that Proposition~\ref{transfer}
is applicable in our situation.

\begin{proposition}
\label{nablanice}
The asymptotic fiber transport $\nabla$ on the northern-hemisphere bundle
$p_{\IH \tilde{M}} :S^+ \IH \tilde{M} \to \IH \tilde{M}$ is $\Gamma$-invariant, 
functorial and induces a homotopically trivial action. Compare Subsection~\ref{subsec_transfer}.
The geodesic flow on $S^+ \IH \tilde{M}$ contracts the asymptotic 
fiber transport $\nabla$. Compare Definition \ref{foliated_contracted}.
\end{proposition}

\begin{proof} 
$\Gamma$-invariance and functoriality follow from the definitions. We have a homotopically trivial
action because the fiber is contractible. The last statement follows from \ref{nablacontr}.
\end{proof}

We now want to construct a metric space $X$ together with a contracting map 
\[ 
p_X: S \IH \tilde{M} \to X.
\]
This map will play a crucial role
in the last step in the proof of Theorem~\ref{surjectivity}.
It will be used to deal with the problem that the  foliated squeezing Theorem
only improves control over some $\Gamma$-compact piece like $S \IH \tilde{M}_{\left[ -t , t \right]}$.

Fix a smooth function
$\phi: \IR \to \left[ 0 , 1 \right]$ with $\phi^{-1}(0)= \left( - \infty , -1 \right]
\cup \left[ 1, \infty \right)$ and $\phi^{-1}(1)=\left[-0.5,0.5 \right]$.
Let $T\IH \tilde{M}$ be the total space of the tangent bundle
and let $q: T\IH \tilde{M} \to \IR$ be the composition of the bundle projection
$T \IH \tilde{M} \to \IH \tilde{M}$  with the standard projection
$\IH \tilde{M} = \IR \times \tilde{M} \to \IR$. Define the topological space $X$ as 
\[
X=\{ v | q(v) \in \left[ -1 , 1 \right] \mbox{ and }  |v|=\phi( q(v) ) \} \subset T\IH \tilde{M}
\]
and define $p_X$ as
\begin{eqnarray*}
T \IH \tilde{M} \supset S \IH \tilde{M} & \to & X \\
v & \mapsto & \begin{array}{lcl}
                           \phi(q(v))v          & \mbox{ if } & v \in    T \IH \tilde{M}_{[-1,1]} \\
                            r_{\pm}( 0 \cdot v) & \mbox{ if } & v \notin T \IH \tilde{M}_{[-1,1]} 
                     \end{array} 
\end{eqnarray*}
Here in the second case we identify the zero section of the tangent bundle 
with $\IH \tilde{M}$ and
\[
r_-:\IH \tilde{M}_{\left( - \infty , -1\right]}  = 
\left( - \infty , -1 \right] \times \tilde{M}  \to \IH \tilde{M}_{\{-1\}}= \{-1\} \times \tilde{M}
\]
is the map which sends $(x,t)$ to $(x,-1)$ and similarly $r_+(x,t)=(x,1)$.
Note that restricted to $S \IH \tilde{M}_{[-0.5,0.5]}$ the map $p_X$ is the identity.
\begin{proposition} \label{cooling}
The space $X$ is compact and can be equipped with the structure of a differentiable 
manifold. There exists a metric $d$ on $X$
which is not Riemannian but generates the topology.
With respect to this metric the map $p_X: S \IH \tilde{M} \to X$ 
satisfies 
\begin{enumerate}
\item For $v,w \in S \IH \tilde{M}$ we have 
\[
      d_X(p(v),p(w)) \leq 
        \mathit{dist}_{S\IH\tilde{M}}(v,w).
\]
\item For all $\alpha > 0$ and $\e > 0$ there is $t_0 = t_0(\alpha,\e) > 0$ 
such that for all $v$, $w$ in $S \IH \tilde{M}$
outside of $S \IH \tilde{M}_{ \left[ -t_0 , t_0 \right] }$
we have the following implication:
\[
\mathit{dist}_{S\IH\tilde{M}}(v,w) 
\leq \alpha \Longrightarrow d_X(p(v),p(w)) \leq \e.
\]
\end{enumerate}
\end{proposition}

\begin{proof}
There clearly exists a number $c \in \left[ 1, \infty \right)$ such that
$| \phi^{\prime} (t) | \leq c$ for all $t \in \IR$. 
As a topological space 
\[
X=X^{+} \cap q^{-1} \left( \left[ -1 , 1 \right] \right), \mbox{ where } X^{+}
 = \{ v \in T \IH \tilde{M} | |v|= \phi ( q(v) ) \}.
\]
We can factor the map $p=p_X : S \IH \tilde{M} \to X$ through $X^+$, i.e.\
we consider the commutative triangle 
\[
\xymatrix{ S \IH \tilde{M} \ar[dr]_{p^+} \ar[rr]^{p_X}  & & X  \\
            & X^+  \ar[ur]_{r}& 
         }
\]
Here $p^+ $ is defined by $p^+ ( v) = \phi ( q (v) ) v $ and the retraction $r$ is determined as
follows: If we identify $\IH \tilde{M} = \IR \times \tilde{M}$ with the zero section of 
$T \IH \tilde{M} \to \IH \tilde{M} $ we can write 
\[
X^+ = X \cup ( -\infty , -1 ] \times \tilde{M} \cup [ 1, \infty ) \times \tilde{M}.
\]
We define $r((t,x)) = (-1,x)$ for $t \in (-\infty , -1 ]$, $r((t,x))=(1,x)$ for 
$t \in [1 , \infty )$ and $r(y) = y$ if $ y \in X$.

We proceed to put a metric $d$ on $X$. Recall first that the Riemannian metric on $\IH \tilde{M}$
naturally determines a Riemannian metric on $T \IH \tilde{M}$ as follows: Identify a smooth curve
$A(t)$ in $T \IH \tilde{M}$ with a smooth vector field $J(t)$ along a smooth curve $\alpha(t)$
in $\IH \tilde{M}$, then
\begin{eqnarray} \label{tangentmetric}
|\dot{A} (0) |^2 = | \dot{\alpha} (0) |^2 + |\dot{J} (0)|^2
\end{eqnarray}
Furthermore we denote by $|A|$ the arclength of a piecewise smooth curve 
$A:\left[ a,b \right] \to T \IH \tilde{M}$. Let $x$, $y \in X$, then 
\begin{eqnarray}  \label{ddef}
d(x,y)= \frac{1}{\sqrt{2 + c^2}} 
\inf \left\{ |A|   \left| \rule[-1em]{0cm}{1.5em}  \right.  
\parbox{32ex}{$A$ is a piecewise smooth curve  in $ X^+$ 
connecting $x$ to $ y$} \right\} .
\end{eqnarray}

Let $A(t)$ be a smooth curve in $S \IH \tilde{M}$ given by the unit length vector field $J(t)$
along the curve $\alpha (t)$ in $\IH \tilde{M}$. Then the smooth curve $p^+ \circ A$ in
$X^+$ is given by the scaled vector field $K(t) = \phi( q ( J(t) ) ) J(t)$ along the same curve
$\alpha (t) $ in $\IH \tilde{M} $. 
\begin{claim} 
Let $u$ and $v$ be the tangent vectors (in $T \IH \tilde{M}$ )
at time $t=0$ to $A$ respectively  $p^+ \circ A$.
We have 
\begin{eqnarray} \label{remark1}
|v| \leq \sqrt{ 2 + c^2 } |u|.
\end{eqnarray}
\end{claim}
{\em Proof of the claim.} It clearly suffices to establish the squared inequality. In 
light of (\ref{tangentmetric}) this is equivalent to 
\begin{eqnarray} \label{eleven}
|\dot{K} (0) |^2 \leq (2 + c^2) | \dot{J} (0) |^2 + (1+c^2) |\dot{\alpha} (0) |^2 
\end{eqnarray}
To establish this let $\delta(t)$ and $\beta(t)$ be the components of $\alpha (t)$, i.e.\
$\alpha (t) = ( \delta (t) , \beta (t) ) \in \IR \times \tilde{M} $. 
Hence $\phi( q ( J (t) )) = \phi( \delta(t) )$ and 
\begin{eqnarray} \label{13}
\dot{K} (0)  = 
\phi^{\prime} ( \delta (0) ) \delta^{\prime} (0)  J(0) +
\phi ( \delta (0) ) \dot{J} (0) 
\end{eqnarray}
Consequently $| \dot{K} (0) | \leq c | \dot{\alpha} (0) | + | \dot{J} (0) |$. Squaring
both sides yields 
\begin{eqnarray} \label{15}
| \dot{K} (0) |^2 \leq 
c^2 | \dot{\alpha} (0) |^2 + | \dot{J} (0) |^2 + 
2 c | \dot{\alpha} (0) | |\dot{J} (0) | .
\end{eqnarray}
Therefore to establish inequality (\ref{eleven}) it suffices to show that
\begin{eqnarray}  \label{16}
c^2 x^2 + y^2 + 2cxy \leq (1+c^2) x^2 + (2+c^2) y^2 
\end{eqnarray}
with $x=| \dot{\alpha} (0) | $ and $ y = | \dot{J} (0) |$.
Subtracting the left from the right side we see that this is equivalent to 
\begin{eqnarray}  \label{18}
0 \leq x^2 + (1+c^2) y^2 - 2c xy = (x-cy)^2 + y^2.
\end{eqnarray}
This establishes the claim.

We now prove (ii).
Let $\alpha(t)=(\delta(t),\beta(t))$ be a piecewise smooth curve in $\IH \tilde{M}=\IR \times \tilde{M}$
such that $|\delta(t)| \geq s_0 \geq 1 $ for some constant $s_0$.
Since $\alpha(t)$ is a curve in $X^+$, we can form
the new curve 
\[
A(t)=r(\alpha(t)).
\]
Then $A$ is also a piecewise smooth curve.
>From the definition of $r$ and the definition of the warped product Riemannian metric on
$\IR \times_{cosh(t)} \tilde{M} = \IH \tilde{M} $ it follows that
\begin{eqnarray} \label{remark2}
|A| \leq \frac{cosh(1)}{cosh(s_0)} |\alpha|.
\end{eqnarray}
Let $A:\left[ a , b \right] \to S\IH \tilde{M}$ be a piecewise smooth path given by a vector 
field $V(t)$ along a piecewise smooth curve $\alpha(t)=( \delta(t) , \beta(t))$ in
$\IH \tilde{M}$.

For every $t \in \left[ a , b \right]$ we have 
\begin{eqnarray} \label{remark3}
|\delta(t)| \geq |\delta(a)| - |A|,
\end{eqnarray}
since it is an immediate consequence of the definition of the warped product metric and of 
(\ref{tangentmetric}) that  $|A| \geq |\alpha | \geq |\delta |$.
>From (\ref{remark2}) and (\ref{remark3}) together with the factorization $p^+ = r \circ p$
and the definition of the metric $d$ we conclude that (ii) of the Proposition holds.

The statement (i) in the Proposition is a direct consequence of the following assertion.
Given $\epsilon > 0$ and a piecewise smooth curve $A: \left[ a , b \right] \to X^+$, there
exists a second piecewise smooth curve $\hat{A}: \left[ a , b \right] \to X^+$ satisfying
\begin{eqnarray}
\hat{A} (a) = r ( A ( a ) ), \quad \hat{A} ( b ) = r( A(b) ) \mbox{ and }
|\hat{A} | \leq |A| + 3 \epsilon . 
\end{eqnarray}

To prove this define for $\epsilon \geq 0$ the subspace $X_{\epsilon}$ of $X^+$ as 
\[
X_{\epsilon} = 
X \cup \left[ -1-\epsilon , -1 \right] \times \tilde{M} \cup \left[ 1 , 1 + \epsilon \right] \times \tilde{M}.
\]
Also define a retraction $r_{\epsilon} : X^+ \to X_{\epsilon} $ similar to $r$. 
Note that $X_0=X$ and $r_0=r$. Observe also that $X^+ - X$ is a smooth manifold and that 
$\tilde{M}_{\epsilon} = \{ -(1 + \epsilon) , 1 + \epsilon \} \times \tilde{M} $ is a smooth 
codimension one submanifold when $\epsilon > 0$. Next pick a number $\epsilon_0 \in (0 , \epsilon ) $
such that $A(a)$ and $A(b)$ do not lie in $M_{\epsilon_0}$.
satisfying $q(A(a)) \neq \epsilon_0$ and $q(A(b)) \neq \epsilon_0$. Then perturb $A$
slightly to a new piecewise smooth curve $A_0 : \left[ a , b \right] \to X^+$ satisfying
\begin{eqnarray} 
&& A_0(a) = A(a) \mbox{  and  } A_0 ( b ) = A( b ) .  \\
&& A_0 \mbox{  is transverse to } \tilde{M}_{\epsilon_0} \mbox{  inside of  } X^+ - X.  \label{21} \\
&& |A_0| \leq |A| + \epsilon.
\end{eqnarray}
Now set $A_1= r_{\epsilon_0} \circ A_0$. Then $A_1$ is a piecewise smooth curve in $X^+$
satisfying $|A_1| \leq |A_0|$ because of an
$r_{\epsilon}$-variant of (\ref{remark2}), cf.(\ref{21}).
Also we can connect $r(A(a))$ to $A_1(a)$ and $r(A(b))$ to
$A_1(b)$ and $r(A_1(b))$ by smooth curves in $\IH \tilde{M}$ of length less than $\epsilon$.
Concatenating these two curves with $A_1$ gives the desired curve $\hat{A}$.

It remains to prove that the metric
$d$ in (\ref{ddef}) above induces the subspace topology on $X$ considered
as a subspace of $T \IH \tilde{M}$.
This is an elementary but lengthy argument which we omit.
\end{proof}


\typeout{---------------------surjective ---------------------}


\section{Surjectivity} \label{sec_surjectivity}

Let $\Gamma$ be the fundamental group of a compact Riemannian manifold with 
strictly negative sectional curvature. In this section we prove:

\begin{theorem}
\label{surjectivity}
The map 
\[
A: K_2 \cald^{\Gamma} (E\Gamma ( \Cyc )) \to K_2 \cald^{\Gamma} ( \punkt )
\]
and hence according to \ref{compareassemblymaps} the assembly map
\[
A_{\Cyc \to \All}:H_1^{\Or \Gamma}( E\Gamma( \Cyc )) ; \IK R^{-\infty} ) \to 
H_1^{\Or \Gamma} ( \punkt ; \IK R^{-\infty} ) = K_1( R \Gamma)
\]
is surjective.
\end{theorem}
\begin{proof}
Let $M$ be the compact Riemannian manifold.
We use the notations $\tilde{M}$, $\IH \tilde{M}$, $S \IH \tilde{M}$, $S^{+} \IH \tilde{M}$ etc.\ from
Section~\ref{sec_geometry}. In particular recall the map $p_X: S \IH \tilde{M} \to X$.
Compactness and the negative curvature assumption imply the following:
given any $\alpha>0$ there is only a finite number of closed orbits of length shorter than
$\alpha$. 
Enumerate the closed orbits of the geodesic flow on $S \IH M$ according to their length:
$g_1, g_2, g_3, \ldots$\ , i.e.\ $g_n$ is not longer then $g_{n+1}$.
Note that the closed geodesics  all
lie in $SM \subset S\IH M_{[-0.5,0.5]} \subset S \IH M$ and hence their
preimages $\tilde{g}_i=\pi^{-1}(g_i)$ in $S \IH \tilde{M}$ can be considered as subsets
of $X$.

Let $X(N)$ be the quotient of $X$ obtained by collapsing each of
$\tilde{g}_1, \ldots , \tilde{g}_N$ to a point, i.e.\ we have a pushout
\[
\xymatrix{
\coprod_{i=1}^N  \tilde{g}_i \ar[d] \ar[r] & X \ar[d]^{p_{X(N)}} \\
\pi_0( \coprod_{i=1}^N  \tilde{g}_i) \ar[r] & X(N) .
         }
\]
Consider the following sequence of maps
\[
\xymatrix{
K_2 \cald^{\Gamma} ( X , p_{X(N)} ) \ar[r] &
K_2 \cald^{\Gamma} ( X , \ast ) \ar[r]^-{\cong}_-P &
K_1 \calc^{\Gamma} ( X , \calf_{\Gamma c} ) \ar[r]^-{\cong} &
K_1 R \Gamma
         }
\]
We have the indicated isomorphisms by Lemma~\ref{pedersensmap} and Lemma~\ref{XversusGamma}.
We will show that for a given $a \in K_1 R \Gamma$ there exists an $N$ such that 
$a$ lies in the image of the above composition.
Using \ref{independence-resolution} we can replace the left hand map with the map
$K_2 \cald^{\Gamma} (X (N) ) \to K_2 \cald^{\Gamma} ( \punkt )$
induced from $X(N) \to \punkt$ via the standard resolution.
Since $X(N)$ is a space with 
cyclic isotropy there is a map $X(N) \to E\Gamma(\Cyc)$ by the universal property
of $E\Gamma (\Cyc)$. Hence we can factorize over the assembly map
$A$ and we see that $a$ also lies in the image of the assembly map.

By Remark~\ref{resolvedversion} we can in fact replace 
$\coprod_{i=1}^N  \tilde{g}_i$ by a homotopy equivalent $\Gamma$-invariant  tubular neighborhood
with $T(N) \subset S \IH \tilde{M}_{[-0.5,0.5]}$. 
This tubular neighborhood $T(N)$ will be specified in Step~2 below .
We write $p_{X(T(N))}:X \to X(T(N))$ for the map which collapses each path component 
of $T(N)$ to a point.

Let $a \in K_1 R \Gamma$ be given. According to \ref{smalladdendum} our
task is the following: 
We have to find a number $N$ and a tubular neighborhood $T(N)$ (which will also depend on $a$)
such that for arbitrary 
small $\epsilon > 0$ (note that we do not know the $\epsilon_0$
in \ref{smallisintheimage}) we can find an automorphism
$\phi$ in $\calc(X , \calf_{\Gamma c} )^{\Gamma}$ 
such that 
\begin{enumerate}
\item
The class of $\phi$ maps to $a$ under the isomorphism from Lemma~\ref{XversusGamma}.
\item
$\phi$ and $\phi^{-1}$ are $\epsilon$-controlled over $X-T(N)$ and
\item
$\phi$ is $p_{X(T(N))}$-separating.
\end{enumerate}

We will proceed in three steps. By Lemma~\ref{XversusGamma} all maps in the following
commutative diagram are isomorphisms.
\[
\xymatrix{
K_1 \calc^{\Gamma}(  S^{+} \IH \tilde{M} , \calf_{\Gamma c} ) 
\ar[d]^-{p_{\IH \tilde{M}}} \ar[drr]^-{\cong} \ar[r]^-i &
K_1 \calc^{\Gamma}( S \IH  \tilde{M} , \calf_{\Gamma c} ) \ar[dr]^-{\cong} \ar[r]^{p_X} &
K_1 \calc^{\Gamma}( X , \calf_{\Gamma c} ) \ar[d]^-{\cong}  \\
K_1 \calc^{\Gamma}( \IH \tilde{M} , \calf_{\Gamma c} )  \ar[rr]^-{\cong} &
 & K_1 R \Gamma .
         }
\]
{\bf Step 1. Transfer and geodesic flow.}
Choose an automorphism $\phi$ in $\calc^{\Gamma}( \IH \tilde{M} , \calf_{\Gamma c} )$
such that $\left[ \phi \right]$ is mapped to $a$.
Proposition~\ref{nablanice} verifies
that the asymptotic fiber transport $\nabla$ for the bundle 
$p_{\IH \tilde{M}}: S^{+} \IH \tilde{M} \to \IH \tilde{M}$ constructed in the last section
satisfies the assumptions that are needed to apply Proposition~\ref{transfer}.
Note that the fiber of this bundle is a disk and therefore its Euler-characteristic is $1$.
We achieve the following:
\begin{quote}
There is an $\alpha_0 > 0$ such that for every $\delta > 0$ we find an 
$(\alpha_0,\delta)$-automorphism $\phi_\delta$ in 
$\calc^{\Gamma}(S^+\IH\tilde{M}, \calf_{\Gamma c} )$ with 
$p_{\IH \tilde{M} } (\left[ \phi_{\delta} \right]) = \phi$.
\end{quote}

{\bf Step 2. Applying the Foliated Squeezing Theorem.}

Use the inclusion
$i:S^+\IH\tilde{M} \hra S\IH\tilde{M}$ to consider $\phi_\delta$ as an
$(\alpha_0,\delta)$-automorphism in $\calc^{\Gamma}(S\IH\tilde{M}, \calf_{\Gamma c}  )$.
Let $\mu_1=\mu_1(\dim S \IH \tilde{M})$ be the constant in the 
Foliated Squeezing Theorem~\ref{foliatedcontrol}.

Now choose $N$ such that $g_1, g_2, \dots , g_N$ are the closed orbits of length less than $\mu_1 \alpha_0$.
According to Lemma~\ref{assumptiononT} we can choose a tubular neighborhood $T(N)$ such that there exists
a $\delta_s=\delta_s(\mu_1 \alpha_0 , N)$ with the property that the $(\mu_1 \alpha_0 , \delta_s)$-thickening
of $T(N)$ still consists of $N$ disjoint path components. We also assume that
$T(N) \subset S \IH \tilde{M}_{[-0.5,0.5]}$.
Note that $T(N)$ depends on $\alpha_0$ and hence on the chosen element $a$.

Furthermore pick $t_0=t_0(\mu_1  \alpha_0 +1 , \epsilon ) \geq 1$ as in Proposition~\ref{cooling} and
set $Y_2= S\IH M_{ \left[ -t_0 , t_0 \right]} - \pi( T(N))$.

The Foliated Squeezing Theorem~\ref{foliatedcontrol}
gives us numbers $\delta_0$ and $\mu_2$ depending on $t_0$, $\alpha_0$, 
the flow and the tubular neighborhood $T(N)$.

Choose $\delta \leq \delta_0$ such that also 
$\mu_2 \delta \leq \delta_s$ and $\delta \leq \frac{1}{2}$.

Consider $\phi_{\delta}$.
The Foliated Squeezing Theorem gives us
a new automorphism $\phi_{new}$ representing the same $K$-theory class as $\phi_{\delta}$ such that 
\begin{enumerate}
\item 
$\phi_{new}$ is $\epsilon$-controlled over $S \IH \tilde{M}_{[-t_0 , t_0] } -  T(N)$ and
\item
$\phi_{new}$ is $(\mu_1  \alpha_0 , \mu_2  \delta )$- and hence $(\mu_1 \alpha_0 , \delta_s)$- and 
$\mu_1  \alpha_0 + 1$-controlled over the whole of $S \IH \tilde{M}$.
\end{enumerate}

{\bf Step 3. Collapsing the non-compact part.}

Let $p_X: S\IH\tilde{M} \to X$ be the projection from Section~\ref{sec_geometry} and set
$\phi = (p_X)_*(\phi_{new})$.
Proposition~\ref{cooling} (ii) and the choice of $t_0$ tells us that $p_X$ improves the 
$\mu_1 \alpha_0 + 1$ control over the non-compact part $S \IH \tilde{M} - S \IH \tilde{ M}_{[-t_0, t_0]}$ to
$\epsilon$ control. Moreover \ref{cooling} (i) says that we do not destroy the
$\epsilon$-control we gained in the previous step. Hence $\phi$ is an
$\e$-automorphism over $X-T(N)$. 
Since $\phi_{new}$ is $(\mu_1 \alpha_0 , \delta_s )$-controlled 
Lemma~\ref{assumptiononT} implies that $\phi$ is $p_{X(T(N))}$-separating because 
$\phi_{new}$ and $\phi$ coincide over $T(N) \subset S \IH M_{[-0.5,0.5]}$ where $p_X$ is the identity.
\end{proof}

\begin{lemma} \label{assumptiononT}
Given any $\alpha>0$ and any $N>1$ there exists a constant $\delta_s=\delta_s(\alpha,N)>0$
and a $\Gamma$-invariant tubular neighborhood $T(N)=\coprod_{i=1}^N T_i$ of $\coprod_{i=1}^N \tilde{g}_i$ 
such that the $(\alpha,\delta_s)$-thickening of $T(N)$ in $X$ still
consists of $N$ disjoint path-components.
\end{lemma}

\begin{proof}
Let $S_i$ be pairwise disjoint closed $\Gamma$-invariant tubular neighborhoods of the $\tilde{g}_i$.
Let $Z_i \subset S_i$ consist of all points $x$ such that $x=\Phi_t (y)$ for some
$y$ in the boundary of $S_i$ and some $t$ with $|t|\leq \alpha$. (Recall that $\Phi_t$ is the geodesic flow.)
Since $\tilde{g}_i$ is invariant under the flow it does not intersect $Z_i$. But $Z_i$ and $\tilde{g}_i$
are $\Gamma$-compact. Hence there is a constant $\delta_s$ with $0<\delta_s <\dist (\tilde{g}_i , Z_i )$.
Now let $T_i$ be a tubular neighborhood that is contained in the $\delta_s$-thickening of $\tilde{g}_i$.
\end{proof}


\typeout{----------------------- Literatur -------------------------}
\bibliographystyle{alpha}

\vfill
{\small
\textsc{Westf{\"a}lische Wilhelms-Universit{\"a}t, SFB 478, 48149 M{\"u}nster, Germany\medskip}

\textsc{Department of Mathematics, SUNY, Binghamton, NY 13902, USA\medskip}

\textsc{Department of Mathematics, SUNY, Stony Brook, NY 11794, USA\medskip}

\textsc{Westf{\"a}lische Wilhelms-Universit{\"a}t, SFB 478, 48149 M{\"u}nster, Germany}
}
\vfill




\end{document}